# Random trees and applications*

## Jean-François Le Gall

*DMA-ENS, 45 rue d'Ulm, 75005 PARIS, FRANCE*
*e-mail:* `legall@dma.ens.fr`

**Abstract:** We discuss several connections between discrete and continuous random trees. In the discrete setting, we focus on Galton-Watson trees under various conditionings. In particular, we present a simple approach to Aldous' theorem giving the convergence in distribution of the contour process of conditioned Galton-Watson trees towards the normalized Brownian excursion. We also briefly discuss applications to combinatorial trees. In the continuous setting, we use the formalism of real trees, which yields an elegant formulation of the convergence of rescaled discrete trees towards continuous objects. We explain the coding of real trees by functions, which is a continuous version of the well-known coding of discrete trees by Dyck paths. We pay special attention to random real trees coded by Brownian excursions, and in a particular we provide a simple derivation of the marginal distributions of the CRT. The last section is an introduction to the theory of the Brownian snake, which combines the genealogical structure of random real trees with independent spatial motions. We introduce exit measures for the Brownian snake and we present some applications to a class of semilinear partial differential equations.

**AMS 2000 subject classifications:** Primary 60J80; secondary 05C05, 35J65, 60C05, 60J65.

**Keywords and phrases:** random tree, contour process, conditioned tree, Brownian motion, Brownian excursion, real tree, coding of trees, CRT, Brownian snake, exit measure, partial differential equation.

Received September 2005.

## Introduction

The main purposes of these notes are to present some recent work about continuous genealogical structures and to point at some of their applications. The interest for these continuous branching structures first arose from their connections with the measure-valued branching processes called superprocesses, which have been studied extensively since the end of the eighties. Independently of the theory of superprocesses, Aldous [1],[2] discussed scaling limits of various classes of discrete trees conditioned to be large. In the case of a Galton-Watson tree with a finite variance critical offspring distribution and conditioned to have a large number of vertices, he proved that the scaling limit is a continuous random tree called the Brownian CRT. Moreover, this limiting continuous object can be coded by a normalized Brownian excursion, a fact that is reminiscent

---

of the Brownian snake construction of superprocesses [25]. In the recent years, these ideas have been extended to much more general continuous trees: See in particular [12] for a discussion of Lévy trees, which are the possible scaling limits of critical or subcritical Galton-Watson trees.

Section 1 below is concerned with scaling limits of discrete trees. In fact, we do not really discuss limits of trees, but consider rather certain coding functions of these trees, namely the height function and the contour function (see Fig.1 below). Our main result (Theorem 1.8) states that the rescaled height function associated with a forest of independent (critical, finite variance) Galton-Watson trees converges in distribution towards reflected Brownian motion on the positive half-line. From this, one can derive a variety of analogous limit theorems for a single Galton-Watson tree conditioned to be large. The derivation is quite simple if the conditioning is "non-degenerate": For instance if the tree is conditioned to have height greater than $n$ (resp. total progeny greater than $n$), the scaling limit of the height function will be a Brownian excursion conditioned to have height greater than 1 (resp. duration greater than 1). For degenerate conditionings, things become a little more complicated: The case of a Galton-Watson tree conditioned to have exactly $n$ vertices, corresponding to Aldous theorem [2], is treated under an exponential moment assumption using an idea of Marckert and Mokkadem [31]. We briefly discuss some applications to various classes of "combinatorial" trees.

Although the limit theorems of Section 1 give a lot of useful information about asymptotics of discrete random trees, it is a bit unsatisfactory that they only discuss continuous limits for the coding functions of trees and not for the trees themselves. The formalism of real trees, which is briefly presented in Section 2, provides an elegant way of restating the limit theorems of Section 1 in terms of convergence of trees. The use of real trees for probabilistic purposes seems to be quite recent: See in particular [17] and [12]. In Section 2, we first discuss the coding of a (compact) real tree by a continuous function. This is of course a continuous analogue of the correspondence between a discrete tree and its contour function. This coding makes it possible to get a simple and efficient construction of the CRT and related random trees as trees coded by various kinds of Brownian excursions. As an application, we use some tools of Brownian excursion theory to derive the finite-dimensional marginals of the CRT, which had been computed by Aldous with a very different method.

Section 3 gives an introduction to the path-valued process called the Brownian snake and its connections with certain semilinear partial differential equations. The Brownian snakes combines the genealogical structure of the random real trees studied in Section 2 with spatial motions governed by a general Markov process $\xi$. Informally, each Brownian snake path corresponds to the spatial positions along the ancestral line of a vertex in the tree. The precise definition of the Brownian snake is thus motivated by the coding of real trees that is discussed in Section 2. In view of applications to PDE, we introduce the exit measure from a domain $D$, which is in a sense uniformly spread over the set of exit points of the Brownian snake paths from $D$. We then derive the key integral equation (Theorem 3.11) for the Laplace functional of the exit measure. In the particular case when the underlying spatial motion $\xi$ is $d$-dimensional Brownian motion, this quickly leads to the connection between the Brownian snake and the semilinear PDE $\Delta u = u^2$. Up to some point, this connection can be viewed as a reformulation of (a special case of) Dynkin's work [13] about superprocesses. Although we content ourselves with a simple application to solutions with boundary blow-up, the results of Section 3 provide most of the background that is necessary to understand the recent deep applications of the Brownian snake to the classification and probabilistic representation of solutions of $\Delta u = u^2$ in a domain [33].

The concepts and results that are presented here have many other recent applications. Let us list a few of these. Random continuous trees can be used to model the genealogy of self-similar fragmentations [20]. The Brownian snake has turned out to be a powerful tool in the study of super-Brownian motion: See the monograph [27] and the references therein. The random measure known as ISE (see [3] and Section 3 below), which is easily obtained from the Brownian snake driven by a normalized Brownian excursion, has appeared in asymptotics for various models of statistical mechanics [7],[21],[23]. There are similar limit theorems for interacting particle systems in $\mathbb{Z}^d$: See [5] for the voter model and coalescing random walks and the recent paper [22] for the contact process. Another promising area of application is the asymptotic study of planar maps. Using a bijection between quadrangulations and well-labelled discrete trees, Chassaing and Schaeffer [6] were able to derive precise asymptotics for random quadrangulations in terms of a one-dimensional Brownian snake conditioned to remain positive (see [30] for the definition of this object). The Chassaing-Schaeffer asymptotics have been extended to more general planar maps by Marckert and Miermont [32]. In fact one expects the existence of a universal continuous limit of planar maps that should be described by random real trees of the type considered here.

## 1. From Discrete to Continuous Trees

In this section, we first explain how discrete random trees can be coded by discrete paths called the height function and the contour function of the tree. We then prove that the rescaled height function associated with a forest of independent Galton-Watson trees converges in distribution towards reflecting Brownian motion on the positive half-line. This has several interesting consequences for the asymptotic behavior of various functionals of Galton-Watson forests or trees. We also discuss analogous results for a single Galton-Watson tree conditioned to be large. In particular we recover a famous theorem of Aldous showing that the suitably rescaled contour function of a Galton-Watson tree conditioned to have $n$ vertices converges in distribution towards the normalized Brownian excursion as $n \to \infty$. Consequences for various classes of "combinatorial trees" are outlined in subsection 1.5.

### *1.1. Discrete trees*

We will be interested in (finite) rooted ordered trees, which are also called plane trees in combinatorics (see e.g. [41]). We first introduce the set of labels

$$\mathcal{U} = \bigcup_{n=0}^{\infty} \mathbb{N}^n$$

where $\mathbb{N} = \{1, 2, \dots\}$ and by convention $\mathbb{N}^0 = \{\varnothing\}$. An element of $\mathcal{U}$ is thus a sequence $u = (u^1, \dots, u^n)$ of elements of $\mathbb{N}$, and we set $|u| = n$, so that $|u|$ represents the "generation" of $u$. If $u = (u^1, \dots u^m)$ and $v = (v^1, \dots, v^n)$ belong to $\mathcal{U}$, we write $uv = (u^1, \dots u^m, v^1, \dots, v^n)$ for the concatenation of $u$ and $v$. In particular $u\varnothing = \varnothing u = u$.

The mapping $\pi : \mathcal{U} \backslash \{\varnothing\} \longrightarrow \mathcal{U}$ is defined by $\pi(u^1 \dots u^n) = u^1 \dots u^{n-1}$ ($\pi(u)$ is the "father" of $u$).

A (finite) rooted ordered tree $\mathbf{t}$ is a finite subset of $\mathcal{U}$ such that:

(i) $\varnothing \in \mathbf{t}$.
(ii) $u \in \mathbf{t} \backslash \{\varnothing\} \Rightarrow \pi(u) \in \mathbf{t}$.
(iii) For every $u \in \mathbf{t}$, there exists an integer $k_u(\mathbf{t}) \geq 0$ such that, for every $j \in \mathbb{N}$, $uj \in \mathbf{t}$ if and only if $1 \leq j \leq k_u(\mathbf{t})$

The number $k_u(\mathbf{t})$ is interpreted as the "number of children" of $u$ in $\mathbf{t}$.

We denote by $\mathbf{A}$ the set of all rooted ordered trees. In what follows, we see each vertex of the tree $\mathbf{t}$ as an individual of a population whose $\mathbf{t}$ is the family tree. The cardinality $\#(\mathbf{t})$ of $\mathbf{t}$ is the total progeny.

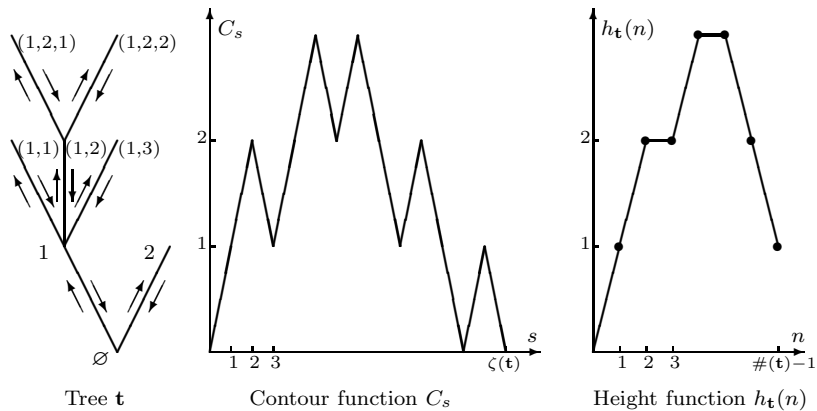

Tree $\mathbf{t}$     Contour function $C_s$     Height function $h_{\mathbf{t}}(n)$

Figure 1

We will now explain how trees can be coded by discrete functions. We first introduce the (discrete) **height function** associated with a tree $\mathbf{t}$. Let us denote by $u_0 = \varnothing, u_1, u_2, \dots, u_{\#(\mathbf{t})-1}$ the elements of $\mathbf{t}$ listed in lexicographical order. The height function $(h_{\mathbf{t}}(n); 0 \leq n < \#(\mathbf{t}))$ is defined by

$$h_{\mathbf{t}}(n) = |u_n|, \quad 0 \leq n < \#(\mathbf{t}).$$

The height function is thus the sequence of the generations of the individuals of **t**, when these individuals are listed in the lexicographical order (see Fig.1 for an example). It is easy to check that $h_{\mathbf{t}}$ characterizes the tree **t**.

The **contour function** (or Dyck path in the terminology of [41]) gives another way of characterizing the tree, which is easier to visualize on a picture (see Fig.1). Suppose that the tree is embedded in the half-plane in such a way that edges have length one. Informally, we imagine the motion of a particle that starts at time $t = 0$ from the root of the tree and then explores the tree from the left to the right, moving continuously along the edges at unit speed (in the way explained by the arrows of Fig.1), until all edges have been explored and the particle has come back to the root. Since it is clear that each edge will be crossed twice in this evolution, the total time needed to explore the tree is $\zeta(\mathbf{t}) := 2(\#(\mathbf{t}) - 1)$. The value $C_s$ of the contour function at time $s \in [0, \zeta(\mathbf{t})]$ is the distance (on the tree) between the position of the particle at time $s$ and the root. By convention $C_s = 0$ if $s \geq \zeta(\mathbf{t})$. Fig.1 explains the construction of the contour function better than a formal definition.

We will introduce still another way of coding the tree. We denote by $\mathcal{S}$ the set of all finite sequences of nonnegative integers $m_1, \ldots, m_p$ (with $p \geq 1$) such that

- $m_1 + m_2 + \cdots + m_i \geq i \quad , \; \forall i \in \{1, \ldots, p-1\}$;
- $m_1 + m_2 + \cdots + m_p = p - 1$.

Recall that $u_0 = \varnothing, u_1, u_2, \ldots, u_{\#(\mathbf{t})-1}$ are the elements of **t** listed in lexicographical order.

**Proposition 1.1** *The mapping*

$$\Phi : \mathbf{t} \longrightarrow (k_{u_0}(\mathbf{t}), k_{u_1}(\mathbf{t}), \ldots, k_{u_{\#(\mathbf{t})-1}}(\mathbf{t}))$$

*defines a bijection from* **A** *onto* $\mathcal{S}$.

**Proof.** We note that if $\#(\mathbf{t}) = p$, the sum $k_{u_0}(\mathbf{t}) + k_{u_1}(\mathbf{t}) + \cdots + k_{u_{\#(\mathbf{t})-1}}(\mathbf{t})$ counts the total number of children of all individuals in the tree and is thus equal to $p-1$ (because $\varnothing$ is not counted !). Furthermore, if $i \in \{0, 1, \ldots, p-2\}$, $k_{u_0} + \cdots + k_{u_i}$ is the number of children of $u_0, \ldots, u_i$ and thus greater than or equal to $i$, because $u_1, \ldots, u_i$ are counted among these children (in the lexicographical order, an individual is visited before his children). There is even a strict inequality because the father of $u_{i+1}$ belongs to $\{u_0, \ldots, u_i\}$. It follows that $\Phi$ maps **A** into $\mathcal{S}$. We leave the rest of the proof to the reader. $\qquad\square$

Let $\mathbf{t} \in \mathbf{A}$ and $p = \#(\mathbf{t})$. Rather than the sequence $(m_1, \ldots, m_p) = \Phi(\mathbf{t})$, we will often consider the finite sequence of integers

$$x_n = \sum_{i=1}^{n} (m_i - 1) , \qquad 0 \leq n \leq p$$

which satisfies the following properties

- $x_0 = 0$ and $x_p = -1$.
- $x_n \geq 0$ for every $0 \leq n \leq p - 1$.
- $x_i - x_{i-1} \geq -1$ for every $1 \leq i \leq p$.

Such a sequence is called a Lukasiewicz path. Obviously the mapping $\Phi$ of the proposition induces a bijection between trees (in **A**) and Lukasiewicz paths.

We now observe that there is a simple relation between the height function of a tree and its Lukasiewicz path.

**Proposition 1.2** *The height function $h_{\mathbf{t}}$ of a tree $\mathbf{t}$ is related to the Lukasiewicz path of $\mathbf{t}$ by the formula*

$$h_{\mathbf{t}}(n) = \#\{j \in \{0, 1, \ldots, n-1\} : x_j = \inf_{j \leq \ell \leq n} x_\ell\},$$

*for every $n \in \{0, 1, \ldots, \#(\mathbf{t}) - 1\}$.*

**Proof.** Obviously,

$$h_{\mathbf{t}}(n) = \#\{j \in \{0, 1, \ldots, n-1\} : u_j \prec u_n\}$$

where $\prec$ stands for the genealogical order on the tree ($u \prec v$ if $v$ is a descendant of $u$). Thus it is enough to prove that for $j \in \{0, 1, \ldots, n-1\}$, $u_j \prec u_n$ holds iff

$$x_j = \inf_{j \leq \ell \leq n} x_\ell.$$

To this end, it suffices to verify that

$$\inf\{k \geq j : x_k < x_j\}$$

is equal either to $\#(\mathbf{t})$, in the case when all $u_k$ with $k > j$ are descendants of $u_j$, or to the first index $k > j$ such that $u_k$ is not a descendant of $u_j$.

However, writing

$$x_k - x_j = \sum_{i=j+1}^{k} (m_i - 1)$$

and using the same arguments as in the proof of Proposition 1.1 (to prove that $\Phi$ takes values in $\mathcal{S}$), we see that for every $\ell > j$ such that $u_\ell$ is a descendant of $u_j$, we have $x_\ell - x_j \geq 0$, whereas on the other hand $x_k - x_j = -1$ if $k$ is the first $\ell > j$ such that $u_\ell$ is not a descendant of $j$ (or $k = p$ if there are no such $\ell$). This completes the proof. $\qquad\square$

### *1.2. Galton-Watson trees*

Let $\mu$ be a critical or subcritical offspring distribution. This means that $\mu$ is a probability measure on $\mathbb{Z}_+$ such that

$$\sum_{k=0}^{\infty} k\mu(k) \leq 1.$$

We exclude the trivial case where $\mu(1) = 1$.

We will make use of the following explicit construction of Galton-Watson trees: Let $(K_u, u \in \mathcal{U})$ be a collection of independent random variables with law $\mu$, indexed by the label set $\mathcal{U}$. Denote by $\theta$ the random subset of $\mathcal{U}$ defined by

$$\theta = \{u = u^1 \ldots u^n \in \mathcal{U} : u^j \leq K_{u^1 \ldots u^{j-1}} \text{ for every } 1 \leq j \leq n\}.$$

**Proposition 1.3** *$\theta$ is a.s. a tree. Moreover, if*

$$Z_n = \#\{u \in \theta : |u| = n\},$$

*$(Z_n, n \geq 0)$ is a Galton-Watson process with offspring distribution $\mu$ and initial value $Z_0 = 1$.*

**Remark.** Clearly $k_u(\theta) = K_u$ for every $u \in \theta$.

The tree $\theta$, or any random tree with the same distribution, will be called a Galton-Watson tree with offspring distribution $\mu$, or in short a $\mu$-Galton-Watson tree. We also write $\Pi_\mu$ for the distribution of $\theta$ on the space $\mathbf{A}$.

We leave the easy proof of the proposition to the reader. The finiteness of the tree $\theta$ comes from the fact that the Galton-Watson process with offspring distribution $\mu$ becomes extinct a.s., so that $Z_n = 0$ for $n$ large.

If $\mathbf{t}$ is a tree and $1 \leq j \leq k_\varnothing(\mathbf{t})$, we write $T_j\mathbf{t}$ for the tree $\mathbf{t}$ shifted at $j$:

$$T_j\mathbf{t} = \{u \in \mathcal{U} : ju \in \mathbf{t}\}.$$

Note that $T_j\mathbf{t}$ is a tree.

Then $\Pi_\mu$ may be characterized by the following two properties (see e.g. [34] for more general statements):

(i) $\Pi_\mu(k_\varnothing = j) = \mu(j), \quad j \in \mathbb{Z}_+$.
(ii) For every $j \geq 1$ with $\mu(j) > 0$, the shifted trees $T_1\mathbf{t}, \ldots, T_j\mathbf{t}$ are independent under the conditional probability $\Pi_\mu(d\mathbf{t} \mid k_\varnothing = j)$ and their conditional distribution is $\Pi_\mu$.

Property (ii) is often called the branching property of the Galton-Watson tree.

We now give an explicit formula for $\Pi_\mu$.

**Proposition 1.4** *For every $\mathbf{t} \in \mathbf{A}$,*

$$\Pi_\mu(\mathbf{t}) = \prod_{u \in \mathbf{t}} \mu(k_u(\mathbf{t})).$$

**Proof.** We can easily check that

$$\{\theta = \mathbf{t}\} = \bigcap_{u \in \mathbf{t}} \{K_u = k_u(\mathbf{t})\},$$

so that

$$\Pi_\mu(\mathbf{t}) = P(\theta = \mathbf{t}) = \prod_{u \in \mathbf{t}} P(K_u = k_u(\mathbf{t})) = \prod_{u \in \mathbf{t}} \mu(k_u(\mathbf{t})).$$

$\square$

Recall from Proposition 1.1 the definition of the mapping $\Phi$.

**Proposition 1.5** *Let $\theta$ be a $\mu$-Galton-Watson tree. Then*

$$\Phi(\theta) \overset{\text{(d)}}{=} (M_1, M_2, \ldots, M_T),$$

*where the random variables $M_1, M_2, \ldots$ are independent with distribution $\mu$, and*

$$T = \inf\{n \geq 1 : M_1 + \cdots + M_n < n\}.$$

**Remark.** The fact that $T < \infty$ a.s. is indeed a consequence of our approach, but is also easy to prove directly by a martingale argument.

**Proof.** We may assume that $\theta$ is given by the preceding explicit construction. Write $U_0 = \varnothing, U_1, \ldots, U_{\#(\theta)-1}$ for the elements of $\theta$ listed in lexicographical order, in such a way that

$$\Phi(\theta) = (K_{U_0}, K_{U_1}, \ldots, K_{U_{\#(\theta)-1}}).$$

We already know that $K_{U_0} + \cdots + K_{U_n} \geq n+1$ for every $n \in \{0, 1, \ldots, \#(\theta)-2\}$, and $K_{U_0} + \cdots + K_{U_{\#(\theta)-1}} = \#(\theta) - 1$.

It will be convenient to also define $U_p$ for $p \geq \#(\theta)$, for instance by setting

$$U_p = U_{\#(\theta)-1} 1 \ldots 1$$

where in the right-hand side we have added $p - \#(\theta) + 1$ labels 1. Then the proof of the proposition reduces to checking that, for every $p \geq 0$, $K_{U_0}, \ldots, K_{U_p}$ are independent with distribution $\mu$. The point is that the labels $U_j$ are random (they depend on the collection $(K_u, u \in \mathcal{U})$) and so we cannot just use the fact that the variables $K_u, u \in \mathcal{U}$ are i.i.d. with distribution $\mu$. We argue by induction on $p$. For $p = 0$ or $p = 1$, the result is obvious since $U_0 = \varnothing$ and $U_1 = 1$ are deterministic.

Fix $p \geq 2$ and assume that the desired result holds at order $p-1$. Use the notation $u \leq v$ for the lexicographical order on $\mathcal{U}$ (in contrast with $u \prec v$ for the genealogical order !). As usual $u < v$ if $u \leq v$ and $u \neq v$. The point is to observe that, for every fixed $u \in \mathcal{U}$, the random set

$$\theta \cap \{v \in \mathcal{U} : v \leq u\}$$

is measurable with respect to the $\sigma$-field $\sigma(K_v, v < u)$. This readily follows from the construction of $\theta$. As a consequence, the event

$$\{U_p = u\} \cap \{\#(\theta) > p\}$$

is measurable with respect to $\sigma(K_v, v < u)$. It is also easy to see that the same measurability property holds for the event

$$\{U_p = u\} \cap \{\#(\theta) \leq p\}.$$

Hence $\{U_p = u\}$ is measurable with respect to $\sigma(K_v, v < u)$.

Finally, if $g_0, g_1, \ldots, g_p$ are nonnegative functions on $\{0, 1, \ldots\}$,

$$
\begin{aligned}
&E[g_0(K_{U_0}) g_1(K_{U_1}) \ldots g_p(K_{U_p})] \\
&= \sum_{u_0 < u_1 < \cdots < u_p} E\Big[1_{\{U_0 = u_0, \ldots, U_p = u_p\}}\, g_0(K_{u_0}) \ldots g_p(K_{u_p})\Big] \\
&= \sum_{u_0 < u_1 < \cdots < u_p} E\Big[1_{\{U_0 = u_0, \ldots, U_p = u_p\}}\, g_0(K_{u_0}) \ldots g_{p-1}(K_{u_{p-1}})\Big]\, E[g_p(K_{u_p})]
\end{aligned}
$$

because $K_{u_p}$ is independent of $\sigma(K_v, v < u_p)$, and we use the preceding measurability property.

Then $E[g_p(K_{u_p})] = \mu(g_p)$ does not depend on $u_p$, and taking $g_p = 1$ in the preceding formula we see that

$$E[g_0(K_{U_0}) g_1(K_{U_1}) \ldots g_p(K_{U_p})] = E[g_0(K_{U_0}) g_1(K_{U_1}) \ldots g_{p-1}(K_{U_{p-1}})]\, \mu(g_p).$$

An application of the induction assumption completes the proof. □

**Corollary 1.6** *Let $(S_n, n \geq 0)$ be a random walk on $\mathbb{Z}$ with initial value $S_0$ and jump distribution $\nu(k) = \mu(k + 1)$ for every $k \geq -1$. Set*

$$T = \inf\{n \geq 1 : S_n = -1\}.$$

*Then the Lukasiewicz path of a $\mu$-Galton-Watson tree $\theta$ has the same distribution as $(S_0, S_1, \ldots, S_T)$. In particular, $\#(\theta)$ and $T$ have the same distribution.*

This is an immediate consequence of the preceding proposition.

### 1.3. Convergence to Brownian motion

Our goal is to show that the height functions (or contour functions) of Galton-Watson trees (resp. of Galton-Watson forests) converge in distribution, modulo a suitable rescaling, towards Brownian excursions (resp. reflected Brownian motions).

We fix a critical offspring distribution $\mu$ with finite variance $\sigma^2 > 0$. Note that the criticality means that we now have

$$\sum_{k=0}^{\infty} k\mu(k) = 1.$$

Let $\theta_1, \theta_2, \ldots$ be a sequence of independent $\mu$-Galton-Watson trees. With each $\theta_i$ we can associate its height function $(h_{\theta_i}(n), 0 \leq n \leq \#(\theta_i) - 1)$. We then define the height process $(H_n, n \geq 0)$ of the forest by concatenating the functions $h_{\theta_1}, h_{\theta_2}, \ldots$:

$$H_n = h_{\theta_i}(n - (\#(\theta_1) + \cdots + \#(\theta_{i-1})))$$

if $\#(\theta_1) + \cdots + \#(\theta_{i-1}) \leq n < \#(\theta_1) + \cdots + \#(\theta_i)$. Clearly, the function $(H_n, n \geq 0)$ determines the sequence of trees. To be specific, the "$k$-th excursion" of $H$ from 0 (more precisely, the values of $H$ between its $k$-th zero and the next one) is the height function of the $k$-th tree in the sequence.

By combining Corollary 1.6 with Proposition 1.2, we arrive at the following result (cf Corollary 2.2 in [29]).

**Proposition 1.7** *We have for every $n \geq 0$*

$$H_n = \#\{k \in \{0, 1, \ldots, n-1\} : S_k = \inf_{k \leq j \leq n} S_j\}. \tag{1}$$

*where $(S_n, n \geq 0)$ is a random walk with the distribution described in Corollary 1.6.*

This is the main ingredient for the proof of the following theorem. By definition, a reflected Brownian motion (started at the origin) is the absolute value of a standard linear Brownian motion started at the origin. The notation $[x]$ refers to the integer part of $x$.

**Theorem 1.8** *Let $\theta_1, \theta_2, \ldots$ be a sequence of independent $\mu$-Galton-Watson trees, and let $(H_n, n \geq 0)$ be the associated height process. Then*

$$\left(\frac{1}{\sqrt{p}} H_{[pt]}, t \geq 0\right) \xrightarrow[p \to \infty]{\text{(d)}} \left(\frac{2}{\sigma} \gamma_t, t \geq 0\right)$$

*where $\gamma$ is a reflected Brownian motion. The convergence holds in the sense of weak convergence on the Skorokhod space $\mathbb{D}(\mathbb{R}_+, \mathbb{R}_+)$.*

Let us establish the weak convergence of finite-dimensional marginals in the theorem.

Let $S = (S_n, n \geq 0)$ be as in Proposition 1.7. Note that the jump distribution $\nu$ has mean 0 and finite variance $\sigma^2$, and thus the random walk $S$ is recurrent. We also introduce the notation

$$M_n = \sup_{0 \leq k \leq n} S_k, \quad I_n = \inf_{0 \leq k \leq n} S_k.$$

Donsker's invariance theorem gives

$$\left(\frac{1}{\sqrt{p}} S_{[pt]}, t \geq 0\right) \xrightarrow[p \to \infty]{\text{(d)}} (\sigma B_t, t \geq 0) \tag{2}$$

where $B$ is a standard linear Brownian motion started at the origin.

For every $n \geq 0$, introduce the time-reversed random walk $\widehat{S}^n$ defined by

$$\widehat{S}_k^n = S_n - S_{(n-k)^+}$$

and note that $(\widehat{S}_k^n, 0 \leq k \leq n)$ has the same distribution as $(S_n, 0 \leq k \leq n)$. From formula (1), we have

$$H_n = \#\{k \in \{0, 1, \ldots, n-1\} : S_k = \inf_{k \leq j \leq n} S_j\} = \Phi_n(\widehat{S}^n),$$

where for any discrete trajectory $\omega = (\omega(0), \omega(1), \ldots)$, we have set

$$\Phi_n(\omega) = \#\{k \in \{1, \ldots, n\} : \omega(k) = \sup_{0 \leq j \leq k} \omega(j)\}.$$

We also set

$$K_n = \Phi_n(S) = \#\{k \in \{1, \ldots, n\} : S_k = M_k\}.$$

**Lemma 1.9** *Define a sequence of stopping times $T_j$, $j = 0, 1, \ldots$ inductively by setting $T_0 = 0$ and for every $j \geq 1$,*

$$T_j = \inf\{n > T_{j-1} : S_n = M_n\}.$$

*Then the random variables $S_{T_j} - S_{T_{j-1}}$, $j = 1, 2, \ldots$ are independent and identically distributed, with distribution*

$$P[S_{T_1} = k] = \nu([k, \infty[) , \quad k \geq 0.$$

**Proof.** The fact that the random variables $S_{T_j} - S_{T_{j-1}}$, $j = 1, 2, \ldots$ are independent and identically distributed is a straightforward consequence of the strong Markov property. It remains to compute the distribution of $S_{T_1}$.

The invariant measure of the recurrent random walk $S$ is the counting measure on $\mathbb{Z}$. By a standard result, if $R_0 = \inf\{n \geq 1 : S_n = 0\}$, we have for every $i \in \mathbb{Z}$,

$$E\Big[\sum_{n=0}^{R_0-1} 1_{\{S_n=i\}}\Big] = 1.$$

Notice that $T_1 \leq R_0$ and that the random walk takes positive values on $]T_1, R_0[$. It easily follows that for every $i \leq 0$

$$E\Big[\sum_{n=0}^{T_1-1} 1_{\{S_n=i\}}\Big] = 1.$$

Therefore, for any function $g : \mathbb{Z} \longrightarrow \mathbb{Z}_+$,

$$E\Big[\sum_{n=0}^{T_1-1} g(S_n)\Big] = \sum_{i=0}^{-\infty} g(i). \tag{3}$$

Then, for any function $f : \mathbb{Z} \longrightarrow \mathbb{Z}_+$,

$$E[f(S_{T_1})] = E\Big[\sum_{k=0}^{\infty} 1_{\{k<T_1\}} f(S_{k+1}) 1_{\{S_{k+1} \geq 0\}}\Big]$$

$$= \sum_{k=0}^{\infty} E\Big[ 1_{\{k < T_1\}} f(S_{k+1}) 1_{\{S_{k+1} \geq 0\}} \Big]$$

$$= \sum_{k=0}^{\infty} E\Big[ 1_{\{k < T_1\}} \sum_{j=0}^{\infty} \nu(j) f(S_k + j) 1_{\{S_k + j \geq 0\}} \Big]$$

$$= \sum_{i=0}^{-\infty} \sum_{j=0}^{\infty} \nu(j) f(i + j) 1_{\{i+j \geq 0\}}$$

$$= \sum_{m=0}^{\infty} f(m) \sum_{j=m}^{\infty} \nu(j),$$

which gives the desired formula. In the third equality we used the Markov property at time $k$ and in the fourth one we applied (3). □

Note that the distribution of $S_{T_1}$ has a finite first moment:

$$E[S_{T_1}] = \sum_{k=0}^{\infty} k \, \nu([k, \infty)) = \sum_{j=0}^{\infty} \frac{j(j+1)}{2} \nu(j) = \frac{\sigma^2}{2}.$$

The next lemma is the key to the first part of the proof.

**Lemma 1.10** *We have*

$$\frac{H_n}{S_n - I_n} \xrightarrow[n \to \infty]{(P)} \frac{2}{\sigma^2},$$

*where the notation $\xrightarrow{(P)}$ means convergence in probability.*

**Proof.** From our definitions, we have

$$M_n = \sum_{T_k \leq n} (S_{T_k} - S_{T_{k-1}}) = \sum_{k=1}^{K_n} (S_{T_k} - S_{T_{k-1}}). \qquad (4)$$

Using Lemma 1.9 and the law of large numbers (note that $K_n \longrightarrow \infty$), we get

$$\frac{M_n}{K_n} \xrightarrow[n \to \infty]{(\text{a.s.})} E[S_{T_1}] = \frac{\sigma^2}{2}.$$

By replacing $S$ with the time-reversed walk $\widehat{S}^n$ we see that for every $n$, the pair $(M_n, K_n)$ has the same distribution as $(S_n - I_n, H_n)$. Hence the previous convergence entails

$$\frac{S_n - I_n}{H_n} \xrightarrow[n \to \infty]{(P)} \frac{\sigma^2}{2},$$

and the lemma follows. □

From (2), we have for every choice of $0 \leq t_1 \leq t_2 \leq \cdots \leq t_m$,

$$\frac{1}{\sqrt{p}} \Big( S_{[pt_1]} - I_{[pt_1]}, \ldots, S_{[pt_m]} - I_{[pt_m]} \Big) \xrightarrow[p \to \infty]{(d)} \sigma \Big( B_{t_1} - \inf_{0 \leq s \leq t_1} B_s, \ldots, B_{t_m} - \inf_{0 \leq s \leq t_m} B_s \Big).$$

Therefore it follows from Lemma 1.10 that

$$\frac{1}{\sqrt{p}}\Big(H_{[pt_1]},\ldots,H_{[pt_m]}\Big) \xrightarrow[p\to\infty]{\text{(d)}} \frac{2}{\sigma}\Big(B_{t_1} - \inf_{0\le s\le t_1}B_s,\ldots,B_{t_m} - \inf_{0\le s\le t_m}B_s\Big).$$

However, a famous theorem of Lévy states that the process

$$\gamma_t = B_t - \inf_{0\le s\le t}B_s$$

is a reflected Brownian motion. This completes the proof of the convergence of finite-dimensional marginals in Theorem 1.8.

We will now discuss the functional convergence in Theorem 1.8. To this end, we will need more precise estimates. We will give details of the argument in the case when $\mu$ has small exponential moments, that is there exists $\lambda > 0$ such that

$$\sum_{k=0}^{\infty} e^{\lambda k}\,\mu(k) < \infty.$$

Our approach in that case is inspired from [31]. See [28] for a proof in the general case. We first state a lemma.

**Lemma 1.11** *Let $\varepsilon \in (0,\frac{1}{4})$. We can find $\varepsilon' > 0$ and an integer $N \ge 1$ such that, for every $n \ge N$ and $\ell \in \{0,1,\ldots,n\}$,*

$$P[|M_\ell - \frac{\sigma^2}{2}K_\ell| > n^{\frac{1}{4}+\varepsilon}] < \exp(-n^{\varepsilon'}).$$

We postpone the proof of the lemma and complete the proof of Theorem 1.8. We apply Lemma 1.11 with $\varepsilon = 1/8$. Since for every $n$ the pair $(M_n, K_n)$ has the same distribution as $(S_n - I_n, H_n)$, we get that, for every sufficiently large $n$, and $\ell \in \{0,1,\ldots,n\}$,

$$P[|S_\ell - I_\ell - \frac{\sigma^2}{2}H_\ell| > n^{\frac{3}{8}}] < \exp(-n^{\varepsilon'}).$$

Hence

$$P\Big[\sup_{0\le\ell\le n}|S_\ell - I_\ell - \frac{\sigma^2}{2}H_\ell| > n^{\frac{3}{8}}\Big] < n\,\exp(-n^{\varepsilon'}).$$

Let $A \ge 1$ be a fixed integer. We deduce from the preceding bound that, for every $p$ sufficiently large,

$$P\Big[\sup_{0\le t\le A}|S_{[pt]} - I_{[pt]} - \frac{\sigma^2}{2}H_{[pt]}| > (Ap)^{\frac{3}{8}}\Big] < Ap\,\exp(-(Ap)^{\varepsilon'}). \qquad (5)$$

A simple application of the Borel-Cantelli lemma gives

$$\sup_{0\le t\le A}\Big|\frac{S_{[pt]} - I_{[pt]}}{\sqrt{p}} - \frac{\sigma^2}{2}\frac{H_{[pt]}}{\sqrt{p}}\Big| \xrightarrow[p\to\infty]{} 0\,, \qquad \text{a.s.}$$

It is now clear that the theorem follows from the convergence

$$(\frac{1}{\sqrt{p}}(S_{[pt]} - I_{[pt]}), t \geq 0) \xrightarrow[p \to \infty]{(d)} (\sigma(B_t - \inf_{0 \leq s \leq t} B_s), t \geq 0),$$

which is an immediate consequence of (2). $\qquad\square$

We still have to prove Lemma 1.11. We first state a very simple "moderate deviations" lemma for sums of independent random variables.

**Lemma 1.12** *Let $Y_1, Y_2, \ldots$ be a sequence of i.i.d. real random variables. We assume that there exists a number $\lambda > 0$ such that $E[\exp(\lambda|Y_1|)] < \infty$, and that $E[Y_1] = 0$. Then, for every $\alpha > 0$, we can choose $N$ sufficiently large so that for every $n \geq N$ and $\ell \in \{1, 2, \ldots, n\}$,*

$$P[|Y_1 + \cdots + Y_\ell| > n^{\frac{1}{2}+\alpha}] \leq \exp(-n^{\alpha/2}).$$

**Proof.** The assumption implies that $E[e^{\lambda Y_1}] = 1 + c\lambda^2 + o(\lambda^2)$ as $\lambda \to 0$, where $c = \frac{1}{2}\text{var}(Y_1)$. Hence we can find a constant $C$ such that for every sufficiently small $\lambda > 0$,

$$E[e^{\lambda Y_1}] \leq e^{C\lambda^2}.$$

It follows that, for every sufficiently small $\lambda > 0$,

$$P[Y_1 + \cdots + Y_\ell > n^{\frac{1}{2}+\alpha}] \leq e^{-\lambda n^{\frac{1}{2}+\alpha}} E[e^{\lambda(Y_1 + \ldots + Y_\ell)}] \leq e^{-\lambda n^{\frac{1}{2}+\alpha}} e^{Cn\lambda^2}.$$

If $n$ is sufficiently large we can take $\lambda = n^{-1/2}$ and the desired result follows (after also replacing $Y_i$ with $-Y_i$). $\qquad\square$

Let us now prove Lemma 1.11. We choose $\alpha \in (0, \varepsilon/2)$ and to simplify notation we put $m_n = [n^{\frac{1}{2}+\alpha}]$. Then, for every $\ell \in \{0, 1, \ldots, n\}$,

$$P[|M_\ell - \frac{\sigma^2}{2}K_\ell| > n^{\frac{1}{4}+\varepsilon}] \leq P[K_\ell > m_n] + P[|M_\ell - \frac{\sigma^2}{2}K_\ell| > n^{\frac{1}{4}+\varepsilon}; K_\ell \leq m_n]. \tag{6}$$

Recalling (4), we have first

$$P[|M_\ell - \frac{\sigma^2}{2}K_\ell| > n^{\frac{1}{4}+\varepsilon}; K_\ell \leq m_n]$$

$$\leq P\Big[\sup_{0 \leq k \leq m_n} |\sum_{j=1}^{k}((S_{T_j} - S_{T_{j-1}}) - \frac{\sigma^2}{2})| > n^{\frac{1}{4}+\varepsilon}\Big]$$

$$\leq P\Big[\sup_{0 \leq k \leq m_n} |\sum_{j=1}^{k}((S_{T_j} - S_{T_{j-1}}) - \frac{\sigma^2}{2})| > m_n^{\frac{1}{2}+\varepsilon}\Big]$$

$$\leq m_n \exp(-m_n^{\varepsilon/2}),$$

where the last bound holds for $n$ large by Lemma 1.12. Note that we are assuming that $\mu$ has small exponential moments, and the same holds for the law of $S_{T_1}$ by Lemma 1.9, which allows us to apply Lemma 1.12.

We still need to bound the first term in the right-hand side of (6). Plainly,

$$P[K_\ell > m_n] \le P[K_n > m_n] \le P[S_{T_{m_n}} \le M_n],$$

and so

$$P[K_\ell > m_n] \le P[S_{T_{m_n}} \le n^{\frac{1}{2}+\frac{\alpha}{2}}] + P[M_n > n^{\frac{1}{2}+\frac{\alpha}{2}}].$$

Applying again Lemma 1.12, we get that for $n$ large,

$$P[M_n > n^{\frac{1}{2}+\frac{\alpha}{2}}] \le n \sup_{1 \le \ell \le n} P[S_\ell > n^{\frac{1}{2}+\frac{\alpha}{2}}] \le n \exp(-n^{\alpha/4}).$$

Finally,

$$P[S_{T_{m_n}} \le n^{\frac{1}{2}+\frac{\alpha}{2}}] = P[S_{T_{m_n}} - \frac{\sigma^2}{2}m_n \le n^{\frac{1}{2}+\frac{\alpha}{2}} - \frac{\sigma^2}{2}m_n]$$

and since $S_{T_{m_n}} - \frac{\sigma^2}{2}m_n$ is the sum of $m_n$ i.i.d. centered random variables having small exponential moments, we can again apply Lemma 1.12 (or a classical large deviations estimate) to get the needed bound. This completes the proof of Lemma 1.11. $\qquad\square$

### *1.4. Some applications*

Let us first recall some important properties of linear Brownian motion. Let $\beta$ be a standard linear Brownian motion started at 0. Then there exists a continuous increasing process $L_t^0 = L_t^0(\beta)$ called the local time of $\beta$ at 0 such that if $N_\varepsilon(t)$ denotes the number of positive excursions of $\beta$ away from 0 with height greater than $\varepsilon$ and completed before time $t$, one has

$$\lim_{\varepsilon \to 0} 2\varepsilon \, N_\varepsilon(t) = L_t^0$$

for every $t \ge 0$, a.s. The topological support of the measure $dL_t^0$ coincides a.s. with the zero set $\{t \ge 0 : \beta_t = 0\}$. Moreover, the above-mentioned Lévy theorem can be strengthened in the form

$$(B_t - \underline{B}_t, -\underline{B}_t; t \ge 0) \stackrel{(d)}{=} (|\beta_t|, L_t^0(\beta); t \ge 0)$$

where $\underline{B}_t = \inf_{0 \le s \le t} B_s$. See e.g. [38] Chapter VI, Theorem VI.2.3.

Keeping the notation of subsection 1.3, we set for every $n \ge 0$,

$$\Lambda_n = k \quad \text{iff } \#(\theta_1) + \cdots + \#(\theta_{k-1}) \le n < \#(\theta_1) + \cdots + \#(\theta_k)$$

in such a way that $k$ is the index of the tree to which the $n^{\text{th}}$-visited vertex belongs.

The convergence stated in Theorem 1.8 can now be strengthened in the following form:

$$\left(\frac{1}{\sqrt{p}}H_{[pt]}, \frac{1}{\sqrt{p}}\Lambda_{[pt]}; t \geq 0\right) \xrightarrow[p\to\infty]{(\mathrm{d})} \left(\frac{2}{\sigma}|\beta_t|, \sigma L_t^0(\beta); t \geq 0\right). \tag{7}$$

This is a simple consequence of our arguments: It is easily seen that

$$\Lambda_n = 1 - \inf_{j \leq n} S_j = 1 - I_n.$$

On the other hand, we saw that, for every $A > 0$,

$$\sup_{0 \leq t \leq A}\left|\frac{S_{[pt]} - I_{[pt]}}{\sqrt{p}} - \frac{\sigma^2}{2}\frac{H_{[pt]}}{\sqrt{p}}\right| \xrightarrow[p\to\infty]{} 0, \qquad \text{a.s.}$$

Combining with Donsker's theorem, we get

$$\left(\frac{1}{\sqrt{p}}H_{[pt]}, \frac{1}{\sqrt{p}}\Lambda_{[pt]}; t \geq 0\right) \xrightarrow[p\to\infty]{(\mathrm{d})} \left(\frac{2}{\sigma}(B_t - \underline{B}_t), -\sigma\underline{B}_t; t \geq 0\right)$$

and an application of Lévy's theorem in the form recalled above yields (7).

We will now apply (7) to study the asymptotics of a single Galton-Watson tree conditioned to be large. We write $h(\theta) = \sup\{|v| : v \in \theta\}$ for the height of the tree $\theta$. Let us fix $x > 0$ and for every integer $p \geq 1$ denote by $\theta^{\{x\sqrt{p}\}}$ a random tree with distribution

$$\Pi_\mu(da \mid h(a) \geq x\sqrt{p})$$

where we recall that $\Pi_\mu$ is the law of the Galton-Watson tree with offspring distribution $\mu$.

We denote by $H^{\{x\sqrt{p}\}}$ the height function of $\theta^{\{x\sqrt{p}\}}$. By convention, $H_n^{\{x\sqrt{p}\}} = 0$ if $n \geq \#(\theta^{\{x\sqrt{p}\}})$.

**Corollary 1.13** *We have*

$$\left(\frac{1}{\sqrt{p}}H_{[pt]}^{\{x\sqrt{p}\}}, t \geq 0\right) \xrightarrow[p\to\infty]{(\mathrm{d})} \left(\frac{2}{\sigma}e_t^{\sigma x/2}, t \geq 0\right)$$

*where $e^{\sigma x/2}$ is a Brownian excursion conditioned to have height greater than $\sigma x/2$.*

The excursion $e^{\sigma x/2}$ can be constructed explicitly in the following way. Set

$$T = \inf\{t \geq 0 : |\beta_t| \geq \sigma x/2\}$$
$$G = \sup\{t \leq T : \beta_t = 0\}$$
$$D = \inf\{t \geq T : \beta_t = 0\}.$$

Then we may take

$$e_t^{\sigma x/2} = |\beta_{(G+t)\wedge D}|.$$

**Proof.** We rely on (7). From the Skorokhod representation theorem, we may for every $p \geq 1$ find a process $(H_t^{(p)}, \Lambda_t^{(p)})$ such that

$$(H_t^{(p)}, \Lambda_t^{(p)})_{t \geq 0} \overset{\text{(d)}}{=} (\frac{1}{\sqrt{p}} H_{[pt]}, \frac{1}{\sqrt{p}} \Lambda_{[pt]})_{t \geq 0}$$

and

$$(H_t^{(p)}, \Lambda_t^{(p)})_{t \geq 0} \underset{p \to \infty}{\longrightarrow} (\frac{2}{\sigma} |\beta_t|, \sigma L_t^0(\beta))_{t \geq 0}. \tag{8}$$

uniformly on every compact set, a.s.

Set

$$T^{(p)} = \inf\{t \geq 0 : |H_t^{(p)}| \geq x\}$$
$$G^{(p)} = \sup\{t \leq T^{(p)} : H_t^{(p)} = 0\} - p^{-1}$$
$$D^{(p)} = \inf\{t \geq T^{(p)} : H_t^{(p)} = 0\}.$$

The reason for the term $-p^{-1}$ in the formula for $G^{(p)}$ comes from the integer part in $\frac{1}{\sqrt{p}} H_{[pt]}$: We want the process $H_{G^{(p)}+t}^{(p)}$ to stay at 0 during the interval $[0, p^{-1})$.

By construction, $(H_{(G^{(p)}+t) \wedge D^{(p)}}^{(p)}, t \geq 0)$ has the distribution of the (rescaled) height process of the first tree in the sequence $\theta_1, \theta_2, \ldots$ with height greater than $x\sqrt{p}$, which is distributed as $\theta^{\{x\sqrt{p}\}}$. Therefore,

$$(H_{(G^{(p)}+t) \wedge D^{(p)}}^{(p)}, t \geq 0) \overset{\text{(d)}}{=} (\frac{1}{\sqrt{p}} H_{[pt]}^{\{x\sqrt{p}\}}, t \geq 0).$$

The corollary will thus follow from (8) if we can prove that

$$G^{(p)} \overset{\text{a.s.}}{\longrightarrow} G , \quad D^{(p)} \overset{\text{a.s.}}{\longrightarrow} D.$$

Using the fact that immediately after time $T$ the process $|\beta|$ hits levels strictly larger than $\sigma x/2$, we easily get from (8) that

$$T^{(p)} \underset{p \to \infty}{\overset{\text{a.s.}}{\longrightarrow}} T.$$

From this and (8) again it follows that

$$\liminf D^{(p)} \geq D \quad \text{a.s.}$$
$$\limsup G^{(p)} \geq G \quad \text{a.s.}$$

Let us prove that we have also $\limsup D^{(p)} \leq D$ a.s. (the same argument works for $\liminf G^{(p)}$). Let us fix $t > 0$. From the support property of local time, we have

$$L_t^0 > L_D^0 \quad \text{a.s. on } \{D < t\}.$$

Thanks to (8), we get

$$\Lambda_t^{(p)} > \sigma L_D^0 \quad \text{for } p \text{ large, a.s. on } \{D < t\}.$$

Now note that $L_D^0 = L_T^0$. Since $T^{(p)}$ converges to $T$, (8) also shows that $\Lambda_{T^{(p)}}^{(p)}$ converges to $\sigma L_T^0$, and it follows that

$$\Lambda_t^{(p)} > \Lambda_{T^{(p)}}^{(p)} \quad \text{for } p \text{ large, a.s. on } \{D < t\}.$$

Observing that $\Lambda^{(p)}$ stays constant on the interval $[T^{(p)}, D^{(p)})$, we conclude that

$$t \geq D^{(p)} \quad \text{for } p \text{ large, a.s. on } \{D < t\}.$$

This is enough to prove that $\limsup D^{(p)} \leq D$ a.s. $\qquad \square$

Replacing $p$ by $p^2$ and taking $x = 1$, we deduce from the corollary (in fact rather from its proof) that

$$\frac{1}{p^2} \#(\theta^{\{p\}}) \xrightarrow[p \to \infty]{\text{(d)}} \zeta_{\sigma/2}$$

where $\zeta_{\sigma/2} = D - G$ is the length of excursion $e^{\sigma/2}$. Indeed this immediately follows from the convergence of $D^{(p)} - G^{(p)}$ towards $D - G$ and the fact that, by construction

$$\#(\theta^{\{p\}}) \overset{\text{(d)}}{=} p^2(D^{(p^2)} - G^{(p^2)}) + 1$$

in the notation of the preceding proof.

Notice that the Laplace transform of the limiting law is known explicitly :

$$E[\exp(-\lambda \zeta_{\sigma/2})] = \frac{\sigma\sqrt{2\lambda}/2}{\sinh(\sigma\sqrt{2\lambda}/2)} \exp(-\sigma\sqrt{2\lambda}/2).$$

This basically follows from the Williams decomposition of Itô's excursion measure (Theorem XII.4.5 in [38]) and the known formulas for the hitting time of $\sigma/2$ by a three-dimensional Bessel process or a linear Brownian motion started at 0 (see e.g. [38]).

**Exercise.** Show the convergence in distribution of $p^{-1}h(\theta^{\{p\}})$ and identify the limiting law.

We will now discuss "occupation measures". Rather than considering a single tree as above, we will be interested in a finite forest whose size will tend to $\infty$ with $p$. Precisely, we fix $b > 0$, and we set

$$H_n^p = \begin{cases} H_n & \text{if } \Lambda_n \leq bp, \\ 0 & \text{if } \Lambda_n > bp, \end{cases}$$

in such a way that $H^p$ is the height process for a collection of $[bp]$ independent Galton-Watson trees. Then it easily follows from (7) that

$$(\frac{1}{p}H_{[p^2t]}^p, t \geq 0) \xrightarrow[p \to \infty]{\text{(d)}} (\frac{2}{\sigma}|\beta|_{t \wedge \tau_{b/\sigma}}, t \geq 0), \tag{9}$$

where, for every $r > 0$,

$$\tau_r := \inf\{t \geq 0 : L_t^0 > r\}.$$

Indeed, we can write

$$\frac{1}{p}H_{[p^2t]}^p = \frac{1}{p}H_{[p^2(t \wedge \tau_b^{(p)})]}$$

where $\tau_b^{(p)} = \frac{1}{p^2}\inf\{n \geq 0 : \Lambda_n > bp\} = \inf\{t \geq 0 : \frac{1}{p}\Lambda_{[p^2t]} > b\}$. Then we observe from (7) that

$$\big((\frac{1}{p}H_{[p^2t]})_{t \geq 0} \ , \ \tau_b^{(p)}\big) \xrightarrow[p \to \infty]{(\mathrm{d})} \big((\frac{2}{\sigma}\,|\beta_t|)_{t \geq 0} \ , \ \tau_{b/\sigma}\big)$$

and (9) follows.

Taking $b = 1$, we deduce from (9) that, for every $x > 0$,

$$P\Big[\sup_{1 \leq i \leq p} h(\theta_i) > px\Big] \xrightarrow[p \to \infty]{} P\Big[\sup_{t \leq \tau_{1/\sigma}} \frac{2}{\sigma}|\beta_t| > x\Big] = 1 - \exp(-\frac{2}{\sigma^2 x}).$$

The last equality is a simple consequence of excursion theory for linear Brownian motion (see e.g. Chapter XII in [38]). Now obviously

$$P\Big[\sup_{1 \leq i \leq p} h(\theta_i) > px\Big] = 1 - (1 - P[h(\theta) > px])^p$$

and we recover the classical fact in the theory of branching processes

$$P[h(\theta) \geq n] \underset{n \to \infty}{\sim} \frac{2}{\sigma^2 n}.$$

We now set $Z_0^p = p$ and, for every $n \geq 1$,

$$Z_n^p = \sum_{i=1}^p \#\{u \in \theta_i : |u| = n\} = \#\{k \geq 0 : H_k^p = n\}.$$

From Proposition 1.3, we know that $(Z_n^p, n \geq 0)$ is a Galton-Watson branching process with offspring distribution $\mu$. We can thus apply the classical diffusion approximation to this process.

**Theorem 1.14**

$$(\frac{1}{p}Z_{[pt]}^p, t \geq 0) \xrightarrow[p \to \infty]{(d)} (X_t, t \geq 0),$$

*where the limiting process $X$ is a diffusion process with infinitesimal generator $\frac{1}{2}\sigma^2 x\frac{d^2}{dx^2}$, which can be obtained as the unique solution of the stochastic differential equation*

$$dX_t = \sigma\sqrt{X_t}\,dB_t$$
$$X_0 = 1.$$

For a proof, see e.g. Theorem 9.1.3 in Ethier and Kurtz [16]. It is easy to see that, for every $p \geq 1$,

- $Z_n^p$ is a martingale,
- $(Z_n^p)^2 - \sigma^2 \sum_{k=0}^{n-1} Z_k^p$ is a martingale,

which strongly suggests that the limiting process $X$ is of the form stated in the theorem.

The process $X$ is called Feller's branching diffusion. When $\sigma = 2$, this is also the zero-dimensional squared Bessel process in the terminology of [38]. Note that $X$ hits 0 in finite time and is absorbed at 0.

To simplify notation, let us fix $\mu$ with $\sigma = 2$. Let $f_1, \ldots, f_q$ be $q$ continuous functions with compact support from $\mathbb{R}_+$ into $\mathbb{R}_+$. As a consequence of (9) we have

$$\Big( \int_0^{\tau_1^{(p)}} f_i(\frac{1}{p} H_{[p^2 t]}^p) \, dt \Big)_{1 \leq i \leq q} \xrightarrow[p \to \infty]{(d)} \Big( \int_0^{\tau_{1/2}} f_i(|\beta_t|) \, dt \Big)_{1 \leq i \leq q}.$$

On the other hand,

$$\int_0^{\tau_1^{(p)}} f_i(\frac{1}{p} H_{[p^2 t]}^p) \, dt = \frac{1}{p^2} \sum_{n=0}^{\infty} Z_n^p f_i(\frac{n}{p}) = \int_0^{\infty} da \, f_i(\frac{[pa]}{p}) \frac{1}{p} Z_{[pa]}^p.$$

By using Theorem 1.14, we see that

$$\Big( \int_0^{\tau_{1/2}} f_i(|\beta_t|) \, dt \Big)_{1 \leq i \leq q} \overset{(d)}{=} \Big( \int_0^{\infty} da \, f_i(a) \, X_a \Big)_{1 \leq i \leq q}.$$

In other words, the occupation measure of $|\beta|$ over the time interval $[0, \tau_{1/2}]$, that is the measure

$$f \longrightarrow \int_0^{\tau_{1/2}} f(|\beta_t|) \, dt$$

has the same distribution as the measure $X_a da$. We have recovered one of the celebrated Ray-Knight theorems for Brownian local times (see e.g. Theorem XI.2.3 in [38]).

### 1.5. Galton-Watson trees with a fixed progeny

We can also use Theorem 1.8 to recover a famous result of Aldous concerning Galton-Watson trees conditioned to have a large (fixed) number of vertices (see [2], Aldous dealt with the contour function rather than the height function, but this is more or less equivalent as we will see in the next subsection). We will follow an idea of [31]. We assume as in the end of subsection 1.3 that $\mu$ has a small exponential moment. Our results hold without this assumption, but it will simplify the proof.

For every $p \geq 1$ we denote by $\theta^{(p)}$ a $\mu$-Galton-Watson tree conditioned to have $\#(\theta) = p$. For this to make sense we need $P(\#(\theta) = p) > 0$ for every $p \geq 1$, which holds if $\mu(1) > 0$ (in fact, we only need $P(\#(\theta) = p) > 0$ for $p$ large, which holds under an aperiodicity condition on $\mu$). In the notation of subsection 1.2, the distribution of $\theta^{(p)}$ is $\Pi_\mu(da \mid \#(a) = p)$.

We denote by $(H_k^{(p)})_{0 \leq k \leq p}$ the height process of $\theta^{(p)}$, with the convention $H_p^{(p)} = 0$.

We also need to introduce the normalized Brownian excursion $(\mathbf{e}_t)_{0 \leq t \leq 1}$. This is simply the Brownian excursion conditioned to have length 1. For instance, we may look at the first positive excursion of $\beta$ (away from 0) with length greater than 1, write $[G, D]$ for the corresponding time interval, and set

$$E_t = \beta_{(G+t) \wedge D} \quad , \qquad t \geq 0$$

and

$$\mathbf{e}_t = \frac{1}{\sqrt{D - G}} \, E_{(D-G)t} \quad , \qquad 0 \leq t \leq 1.$$

A more intrinsic construction of the normalized Brownian excursion will be presented in the next subsection.

**Theorem 1.15** *We have*

$$(\frac{1}{\sqrt{p}} H_{[pt]}^{(p)}, 0 \leq t \leq 1) \xrightarrow[p \to \infty]{(\mathrm{d})} (\frac{2}{\sigma} \, \mathbf{e}_t, 0 \leq t \leq 1).$$

This is obviously very similar to our previous results Theorem 1.8 and Corollary 1.13. However, because the present conditioning "degenerates in the limit" $p \to \infty$ (there is no Brownian excursion with length exactly equal to 1), we cannot use the same strategy of proof as in Corollary 1.13.

**Proof.** Let $(H_n, n \geq 0)$ be as in Theorem 1.8 the height process associated with a sequence of independent $\mu$-Galton-Watson trees. We may and will assume that $H$ is given in terms of the random walk $S$ as in (1).

Denote by $T_1$ the number of vertices of the first tree in the sequence, or equivalently

$$T_1 = \inf\{n \geq 1 : H_n = 0\} = \inf\{n \geq 0 : S_n = -1\}.$$

A simple combinatorial argument (consider all circular permutations of the $p$ increments of the random walk $S$ over the interval $[0, p]$) shows that, for every $p \geq 1$,

$$P(T_1 = p) = \frac{1}{p} P(S_p = -1).$$

On the other hand classical results for random walk (see e.g. P9 in Chapter II of [40]) give

$$\lim_{p \to \infty} \sqrt{p} \, P(S_p = -1) = \frac{1}{\sigma \sqrt{2\pi}},$$

and it follows that

$$P(T_1 = p) \underset{p \to \infty}{\sim} \frac{1}{\sigma \sqrt{2\pi p^3}}. \tag{10}$$

Recall from the end of the proof of Theorem 1.8 (see (5)) that we can find $\varepsilon > 0$ so that, for $p$ large enough

$$P\Big[\sup_{0 \le t \le 1} |\frac{H_{[pt]}}{\sqrt{p}} - \frac{2}{\sigma^2} \frac{S_{[pt]} - I_{[pt]}}{\sqrt{p}}| > p^{-1/8}\Big] < \exp(-p^\varepsilon).$$

By comparing with (10), we see that we have also for $p$ large

$$P\Big[\sup_{0 \le t \le 1} |\frac{H_{[pt]}}{\sqrt{p}} - \frac{2}{\sigma^2} \frac{S_{[pt]} - I_{[pt]}}{\sqrt{p}}| > p^{-1/8} \,\Big|\, T_1 = p\Big] < \exp(-p^{\varepsilon'}),$$

for any $\varepsilon' < \varepsilon$. Since $I_n = 0$ for $0 \le n < T_1$, we have also for $p$ large

$$P\Big[\sup_{0 \le t \le 1} |\frac{H_{[pt]}}{\sqrt{p}} - \frac{2}{\sigma^2} \frac{S_{[pt]}}{\sqrt{p}}| > p^{-1/8} \,\Big|\, T_1 = p\Big] < \exp(-p^{\varepsilon'}).$$

Now obviously $(H_k^{(p)}, 0 \le k \le p)$ has the same distribution as $(H_k, 0 \le k \le p)$ under $P(\cdot \mid T_1 = p)$. Therefore Theorem 1.15 is a consequence of the last bound and the following lemma which relates the normalized Brownian excursion to the random walk excursion with a fixed long duration.

**Lemma 1.16** *The distribution of the process $(\frac{1}{\sigma\sqrt{p}}S_{[pt]}, 0 \le t \le 1)$ under the conditional probability $P(\cdot \mid T_1 = p)$ converges as $p$ tends to $\infty$ to the law of the normalized Brownian excursion.*

We omit the proof of this lemma, which can be viewed as a conditional version of Donsker's theorem. See Kaigh [24].

See also Duquesne [10] for generalisations of Theorem 1.15.

**Application.** We immediately deduce from Theorem 1.15 that, for every $x > 0$,

$$\lim_{p \to \infty} P(h(\theta) > x\sqrt{p} \mid \#(\theta) = p) = P(\sup_{0 \le t \le 1} \mathbf{e}_t > \frac{\sigma x}{2}). \tag{11}$$

There is an explicit (complicated) formula for the right-hand side of (11).

**Combinatorial consequences.** For several particular choices of $\mu$, the measure $\Pi_\mu(da \mid \#(a) = p)$ coincides with the uniform probability measure on a class of "combinatorial trees" with $p$ vertices, and Theorem 1.15 gives information about the proportion of trees in this class that satisfy certain properties. To make this more explicit, consider first the case when $\mu$ is the geometric distribution with parameter $\frac{1}{2}$ ($\mu(k) = 2^{-k-1}$). Then $\Pi_\mu(da \mid \#(a) = p)$ is the uniform distribution on the set $\mathbf{A}_p$ of all rooted ordered trees with $p$ vertices (this follows from Proposition 1.4). Thus (11) shows that the height of a tree chosen at random in $\mathbf{A}_p$ is of order $\sqrt{p}$, and more precisely it gives the asymptotic proportion of those trees in $\mathbf{A}_p$ with height greater than $x\sqrt{p}$. Similar arguments apply to other functionals than the height: For instance, if $0 \le a < b$ are given, we could derive asymptotics for the number of vertices in the tree between generations $a\sqrt{p}$ and $b\sqrt{p}$, for the number of vertices at generation $[a\sqrt{p}]$ that have descendants at generation $[b\sqrt{p}]$, etc. The limiting distributions are obtained in terms of the normalized Brownian excursion, and are not always easy to compute explicitly !

Similarly, if $\mu$ is the Poisson distribution, that is $\mu(k) = e^{-1}/k!$, then $\Pi_\mu(da \mid \#(a) = p)$ yields the uniform distribution on the set of all rooted Cayley trees with $p$ vertices. Let us explain this more in detail. Recall that a Cayley tree (with $p$ vertices) is an unordered tree on vertices labelled $1, 2, \ldots, p$. The root can then be any of these vertices. By a famous formula due to Cayley, there are $p^{p-2}$ Cayley trees on $p$ vertices, and so $p^{p-1}$ rooted Cayley trees on $p$ vertices. If we start from a rooted ordered tree distributed according to $\Pi_\mu(da \mid \#(a) = p)$, then assign labels $1, 2, \ldots, p$ uniformly at random to vertices, and finally "forget" the ordering of the tree we started from, we get a random tree that is uniformly distributed over the set of all rooted Cayley trees with $p$ vertices. Hence Theorem 1.15, and in particular (11) also give information about the properties of large Cayley trees.

As a last example, we can take $\mu = \frac{1}{2}(\delta_0 + \delta_2)$, and it follows from Proposition 1.4 that, provided $p$ is odd, $\Pi_\mu(da \mid \#(a) = p)$ is the uniform distribution over the set of (complete) binary trees with $p$ vertices. Strictly speaking, Theorem 1.15 does not include this case, since we assumed that $\mu(1) > 0$. It is however not hard to check that the convergence of Theorem 1.15 still holds, provided we restrict our attention to odd values of $p$.

It is maybe unexpected that these different classes of combinatorial trees give rise to the same scaling limit, and in particular that the limiting law appearing in (11) is the same in each case. Note however that the constant $\sigma$ varies: $\sigma^2 = 1$ for (complete) binary trees or for Cayley trees, whereas $\sigma^2 = 2$ for rooted ordered trees.

As a final remark, let us observe that the convergence in distribution of Theorem 1.15 is often not strong enough to deduce rigorously the desired combinatorial asymptotics (this is the case for instance if one looks at the height profile of the tree, that is the number of vertices at every level in the tree). Still Theorem 1.15 allows one to guess what the limit should be in terms of the normalized Brownian excursion. See in particular [9] for asymptotics of the profile that confirmed a conjecture of Aldous.

### 1.6. Convergence of contour functions

In this subsection, we briefly explain how the preceding results can be stated as well in terms of the contour processes of the trees rather than the height processes as discussed above. The contour function of a tree was discussed in subsection 1.1 (see Fig.1). Notice that in contrast to the height process it is convenient to have the contour function indexed by a real parameter.

We will give the result corresponding to Theorem 1.8. So we consider again a sequence $\theta_1, \theta_2, \ldots$ of independent $\mu$-Galton-Watson trees and we denote by $(C_t, t \geq 0)$ the process obtained by concatenating the contour functions of $\theta_1, \theta_2, \ldots$ Here we need to define precisely what we mean by concatenation.

In Fig.1 and the discussion of subsection 1.1, the contour function of a tree $\theta$ is naturally defined on the time interval $[0, \zeta(\theta)]$, where $\zeta(\theta) = 2(\#(\theta) - 1)$. This has the unpleasant consequence that the contour function of the tree consisting only of the root is trivial. For this reason we make the slightly artificial convention that the contour function $C_t(\theta)$ is defined for $0 \leq t \leq \xi(\theta) = 2\#(\theta) - 1$, by taking $C_t = 0$ if $\zeta(\theta) \leq t \leq \xi(\theta)$. We then obtain $(C_t, t \geq 0)$ by concatenating the functions $(C_t(\theta_1), 0 \leq t \leq \xi(\theta_1))$, $(C_t(\theta_2), 0 \leq t \leq \xi(\theta_2))$, etc.

For every $n \geq 0$, we set

$$J_n = 2n - H_n + I_n.$$

Note that the sequence $J_n$ is strictly increasing and $J_n \geq n$.

Recall that the value at time $n$ of the height process corresponds to the generation of the individual visited at step $n$, assuming that individuals are visited in lexicographical order one tree after another. It is easily checked by induction on $n$ that $[J_n, J_{n+1}]$ is exactly the time interval during which the contour process goes from the individual $n$ to the individual $n + 1$. From this observation, we get

$$\sup_{t \in [J_n, J_{n+1}]} |C_t - H_n| \leq |H_{n+1} - H_n| + 1.$$

A more precise argument for this bound follows from the explicit formula for $C_t$ in terms of the height process: For $t \in [J_n, J_{n+1}]$,

$$
\begin{aligned}
C_t &= H_n - (t - J_n) && \text{if } t \in [J_n, J_{n+1} - 1], \\
C_t &= (H_{n+1} - (J_{n+1} - t))^+ && \text{if } t \in [J_{n+1} - 1, J_{n+1}].
\end{aligned}
$$

These formulas are easily checked by induction on $n$.

Define a random function $\varphi : \mathbb{R}_+ \longrightarrow \{0, 1, \ldots\}$ by setting $\varphi(t) = n$ iff $t \in [J_n, J_{n+1})$. From the previous bound, we get for every integer $m \geq 1$,

$$\sup_{t \in [0,m]} |C_t - H_{\varphi(t)}| \leq \sup_{t \in [0, J_m]} |C_t - H_{\varphi(t)}| \leq 1 + \sup_{n \leq m} |H_{n+1} - H_n|. \tag{12}$$

Similarly, it follows from the definition of $J_n$ that

$$\sup_{t \in [0,m]} |\varphi(t) - \frac{t}{2}| \leq \sup_{t \in [0, J_m]} |\varphi(t) - \frac{t}{2}| \leq \frac{1}{2} \sup_{n \leq m} H_n + \frac{1}{2}|I_m| + 1. \tag{13}$$

**Theorem 1.17** *We have*

$$\left(\frac{1}{\sqrt{p}} C_{2pt}, \, t \geq 0\right) \xrightarrow[p \to \infty]{(\mathrm{d})} \left(\frac{2}{\sigma}|\beta_t|, t \geq 0\right). \tag{14}$$

*where $\beta$ is a standard linear Brownian motion.*

**Proof.** For every $p \geq 1$, set $\varphi_p(t) = p^{-1}\varphi(pt)$. By (12), we have for every $m \geq 1$,

$$\sup_{t \leq m} \left| \frac{1}{\sqrt{p}} C_{2pt} - \frac{1}{\sqrt{p}} H_{p\varphi_p(2t)} \right| \leq \frac{1}{\sqrt{p}} + \frac{1}{\sqrt{p}} \sup_{t \leq 2m} |H_{[pt]+1} - H_{[pt]}| \xrightarrow[p \to \infty]{} 0 \tag{15}$$

in probability, by Theorem 1.8.

On the other hand, the convergence (2) implies that, for every $m \geq 1$,

$$\frac{1}{\sqrt{p}} I_{mp} \xrightarrow[p \to \infty]{(d)} \sigma \inf_{t \leq m} B_t. \tag{16}$$

Then, we get from (13)

$$\sup_{t \leq m} |\varphi_p(2t) - t| \leq \frac{1}{p} \sup_{k \leq 2mp} H_k + \frac{1}{p} |I_{2mp}| + \frac{2}{p} \xrightarrow[p \to \infty]{} 0 \tag{17}$$

in probability, by Theorem 1.8 and (16).

The statement of the theorem now follows from Theorem 1.8, (15) and (17).
$\square$

**Remark.** There is one special case where Theorem 1.17 is easy, without any reference to Theorem 1.8. This is the case where $\mu$ is the geometric distribution $\mu(k) = 2^{-k-1}$, which satisfies our assumptions with $\sigma^2 = 2$. In that case, it is not hard to see that away from the origin the contour process $(C_n, n \geq 0)$ behaves like simple random walk (indeed, by the properties of the geometric distribution, the probability for an individual to have at least $n+1$ children knowing that he has at least $n$ is $1/2$ independently of $n$). A simple argument then shows that the statement of Theorem 1.17 follows from Donsker's invariance theorem.

Clearly, Corollary 1.13 and Theorem 1.15 can also be restated in terms of the contour process of the respective trees. Simply replace $H_{[pt]}^{\{p\}}$ by $C_{[2pt]}^{\{p\}}$ (with an obvious notation) in Corollary 1.13 and $H_{[pt]}^{(p)}$ by $C_{[2pt]}^{(p)}$ in Theorem 1.15.

### *1.7. Conclusion*

The various results of this section show that the rescaled height processes (or contour processes) of large Galton-Watson trees converge in distribution towards Brownian excursions. Still we did not assert that the trees themselves converge. In fact, a precise mathematical formulation of this fact requires a formal definition of what the limiting random trees are and what the convergence means. In the next section, we will give a precise definition of continuous trees and discuss a topology on the space of continuous trees. This will make it possible to reinterpret the results of this section as convergence theorems for random trees.

*Bibibliographical notes.* The coding of discrete trees by contour functions (Dyck paths) or Lukasiewicz words is well known: See e.g. [41]. Theorem 1.8 can be viewed as a variant of Aldous' theorem about the scaling limit of the contour function of Galton-Watson trees [2]. The method that is presented here is taken from [28], with an additional idea from [31]. More general statements can be found in Chapter 2 of the monograph [11]. See Chapters 5 and 6 of Pitman [36], and the references therein, for more results about the connections between trees, random walks and Brownian motion.

## 2. Real Trees and their Coding by Brownian Excursions

In this section, we first describe the formalism of real trees, which can be used to give a precise mathematical meaning to the convergence of rescaled discrete trees towards continuous objects. We then show how a real tree can be coded by a continuous function in a way similar to the coding of discrete trees by their contour functions. Aldous' Continuum Random Tree (the CRT) can be defined as the random real tree coded by a normalized Brownian excursion. For every integer $p \geq 1$, we then compute the $p$-dimensional marginal distribution (that is the law of the reduced tree consisting of the ancestral lines of $p$ individuals chosen uniformly at random) of the tree coded by a Brownian excursion under the Itô excursion measure. Via a conditioning argument, this leads to a simple derivation of the marginal distributions of the CRT.

### 2.1. Real trees

We start with a formal definition. In the present work, we consider only *compact* real trees, and so we include this compactness property in the definition.

**Definition 2.1** *A compact metric space $(\mathcal{T}, d)$ is a real tree if the following two properties hold for every $a, b \in \mathcal{T}$.*

(i) *There is a unique isometric map $f_{a,b}$ from $[0, d(a,b)]$ into $\mathcal{T}$ such that $f_{a,b}(0) = a$ and $f_{a,b}(d(a,b)) = b$.*

(ii) *If $q$ is a continuous injective map from $[0,1]$ into $\mathcal{T}$, such that $q(0) = a$ and $q(1) = b$, we have*

$$q([0,1]) = f_{a,b}([0, d(a,b)]).$$

*A rooted real tree is a real tree $(\mathcal{T}, d)$ with a distinguished vertex $\rho = \rho(\mathcal{T})$ called the root. In what follows, real trees will always be rooted, even if this is not mentioned explicitly.*

Let us consider a rooted real tree $(\mathcal{T}, d)$. The range of the mapping $f_{a,b}$ in (i) is denoted by $[\![a, b]\!]$ (this is the line segment between $a$ and $b$ in the tree). In particular, $[\![\rho, a]\!]$ is the path going from the root to $a$, which we will interpret as the ancestral line of vertex $a$. More precisely we define a partial order on the tree by setting $a \preccurlyeq b$ ($a$ is an ancestor of $b$) if and only if $a \in [\![\rho, b]\!]$.

If $a, b \in \mathcal{T}$, there is a unique $c \in \mathcal{T}$ such that $[\![\rho, a]\!] \cap [\![\rho, b]\!] = [\![\rho, c]\!]$. We write $c = a \wedge b$ and call $c$ the most recent common ancestor to $a$ and $b$.

By definition, the multiplicity of a vertex $a \in \mathcal{T}$ is the number of connected components of $\mathcal{T} \backslash \{a\}$. Vertices of $\mathcal{T} \backslash \{\rho\}$ which have multiplicity 1 are called leaves.

Our goal is to study the convergence of random real trees. To this end, it is of course necessary to have a notion of distance between two real trees. We will use the Gromov-Hausdorff distance between compact metric spaces, which has been introduced by Gromov (see e.g. [19]) in view of geometric applications.

If $(E, \delta)$ is a metric space, we use the notation $\delta_{Haus}(K, K')$ for the usual Hausdorff metric between compact subsets of $E$ :

$$\delta_{Haus}(K, K') = \inf\{\varepsilon > 0 : K \subset U_\varepsilon(K') \text{ and } K' \subset U_\varepsilon(K)\},$$

where $U_\varepsilon(K) := \{x \in E : \delta(x, K) \leq \varepsilon\}$.

Then, if $\mathcal{T}$ and $\mathcal{T}'$ are two rooted compact metric spaces, with respective roots $\rho$ and $\rho'$, we define the distance $d_{GH}(\mathcal{T}, \mathcal{T}')$ by

$$d_{GH}(\mathcal{T}, \mathcal{T}') = \inf\{\delta_{Haus}(\varphi(\mathcal{T}), \varphi'(\mathcal{T}')) \vee \delta(\varphi(\rho), \varphi'(\rho'))\}$$

where the infimum is over all choices of a metric space $(E, \delta)$ and all isometric embeddings $\varphi : \mathcal{T} \longrightarrow E$ and $\varphi' : \mathcal{T}' \longrightarrow E$ of $\mathcal{T}$ and $\mathcal{T}'$ into $(E, \delta)$.

Two rooted compact metric spaces $\mathcal{T}_1$ and $\mathcal{T}_2$ are called equivalent if there is a root-preserving isometry that maps $\mathcal{T}_1$ onto $\mathcal{T}_2$. Obviously $d_{GH}(\mathcal{T}, \mathcal{T}')$ only depends on the equivalence classes of $\mathcal{T}$ and $\mathcal{T}'$. Then $d_{GH}$ defines a metric on the set of all equivalent classes of rooted compact metric spaces (cf [19] and [17]). We denote by $\mathbb{T}$ the set of all (equivalence classes of) rooted real trees.

**Theorem 2.1** *The metric space* $(\mathbb{T}, d_{GH})$ *is complete and separable.*

We will not really use this theorem (see however the remarks after Lemma 2.4). So we refer the reader to [17], Theorem 1 for a detailed proof.

We will use the following alternative definition of $d_{GH}$. First recall that if $(\mathcal{T}_1, d_1)$ and $(\mathcal{T}_2, d_2)$ are two compact metric spaces, a correspondence between $\mathcal{T}_1$ and $\mathcal{T}_2$ is a subset $\mathcal{R}$ of $\mathcal{T}_1 \times \mathcal{T}_2$ such that for every $x_1 \in \mathcal{T}_1$ there exists at least one $x_2 \in \mathcal{T}_2$ such that $(x_1, x_2) \in \mathcal{R}$ and conversely for every $y_2 \in \mathcal{T}_2$ there exists at least one $y_1 \in \mathcal{T}_1$ such that $(y_1, y_2) \in \mathcal{R}$. The distortion of the correspondence $\mathcal{R}$ is defined by

$$\mathrm{dis}(\mathcal{R}) = \sup\{|d_1(x_1, y_1) - d_2(x_2, y_2)| : (x_1, x_2), (y_1, y_2) \in \mathcal{R}\}.$$

Then, if $\mathcal{T}$ and $\mathcal{T}'$ are two rooted compact metric spaces with respective roots $\rho$ and $\rho'$, we have

$$d_{GH}(\mathcal{T}, \mathcal{T}') = \frac{1}{2} \inf_{\mathcal{R} \in \mathcal{C}(\mathcal{T}, \mathcal{T}'), (\rho, \rho') \in \mathcal{R}} \mathrm{dis}(\mathcal{R}), \tag{18}$$

where $\mathcal{C}(\mathcal{T}, \mathcal{T}')$ denotes the set of all correspondences between $\mathcal{T}$ and $\mathcal{T}'$ (see Lemma 2.3 in [17] – actually this lemma is stated for trees but the proof applies as well to compact metric spaces).

### 2.2. *Coding real trees*

In this subsection, we describe a method for constructing real trees, which is particularly well-suited to our forthcoming applications to random trees. We consider a (deterministic) continuous function $g : [0, \infty[ \longrightarrow [0, \infty[$ with compact support and such that $g(0) = 0$. To avoid trivialities, we will also assume that $g$ is not identically zero. For every $s, t \geq 0$, we set

$$m_g(s, t) = \inf_{r \in [s \wedge t, s \vee t]} g(r),$$

and

$$d_g(s, t) = g(s) + g(t) - 2m_g(s, t).$$

Clearly $d_g(s, t) = d_g(t, s)$ and it is also easy to verify the triangle inequality

$$d_g(s, u) \leq d_g(s, t) + d_g(t, u)$$

for every $s, t, u \geq 0$. We then introduce the equivalence relation $s \sim t$ iff $d_g(s, t) = 0$ (or equivalently iff $g(s) = g(t) = m_g(s, t)$). Let $\mathcal{T}_g$ be the quotient space

$$\mathcal{T}_g = [0, \infty[ / \sim .$$

Obviously the function $d_g$ induces a distance on $\mathcal{T}_g$, and we keep the notation $d_g$ for this distance. We denote by $p_g : [0, \infty[ \longrightarrow \mathcal{T}_g$ the canonical projection. Clearly $p_g$ is continuous (when $[0, \infty[$ is equipped with the Euclidean metric and $\mathcal{T}_g$ with the metric $d_g$).

We set $\rho = p_g(0)$. If $\zeta > 0$ is the supremum of the support of $g$, we have $p_g(t) = \rho$ for every $t \geq \zeta$. In particular, $\mathcal{T}_g = p_g([0, \zeta])$ is compact.

**Theorem 2.2** *The metric space $(\mathcal{T}_g, d_g)$ is a real tree. We will view $(\mathcal{T}_g, d_g)$ as a rooted tree with root $\rho = p_g(0)$.*

**Remark.** It is also possible to prove that any (rooted) real tree can be represented in the form $\mathcal{T}_g$. We will leave this as an exercise for the reader.

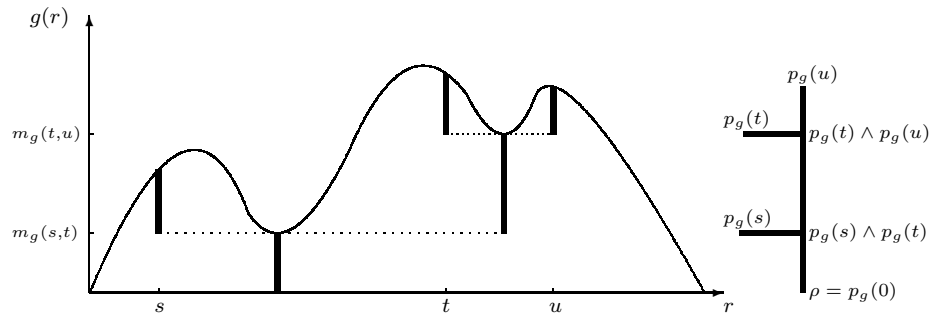

Figure 2

To get an intuitive understanding of Theorem 2.2, the reader should have a look at Figure 2. This figure shows how to construct a simple subtree of $\mathcal{T}_g$, namely the "reduced tree" consisting of the union of the ancestral lines in $\mathcal{T}_g$ of three vertices $p_g(s), p_g(t), p_g(u)$ corresponding to three (given) times $s, t, u \in [0, \zeta]$. This reduced tree is the union of the five bold line segments that are constructed from the graph of $g$ in the way explained on the left part of the figure. Notice that the lengths of the horizontal dotted lines play no role in the construction, and that the reduced tree should be viewed as pictured on the right part of Figure 2. The ancestral line of $p_g(s)$ (resp. $p_g(t), p_g(u)$) is a line segment of length $g(s)$ (resp. $g(t), g(u)$). The ancestral lines of $p_g(s)$ and $p_g(t)$ share a common part, which has length $m_g(s,t)$ (the line segment at the bottom in the left or the right part of Figure 2), and of course a similar property holds for the ancestral lines of $p_g(s)$ and $p_g(u)$, or of $p_g(t)$ and $p_g(u)$.

We present below an elementary proof of Theorem 2.2, which uses only the definition of a real tree, and also helps to understand the notions of ancestral line and most recent common ancestor in $\mathcal{T}_g$. Another argument depending on Theorem 2.1 is presented at the end of the subsection. See also [18] for a short proof using the characterization of real trees via the so-called four-point condition.

Before proceeding to the proof of the theorem, we state and prove the following root change lemma, which is of independent interest.

**Lemma 2.3** *Let $s_0 \in [0, \zeta[$. For any real $r \geq 0$, denote by $\overline{r}$ the unique element of $[0, \zeta[$ such that $r - \overline{r}$ is an integer multiple of $\zeta$. Set*

$$g'(s) = g(s_0) + g(\overline{s_0 + s}) - 2m_g(s_0, \overline{s_0 + s}),$$

*for every $s \in [0, \zeta]$, and $g'(s) = 0$ for $s > \zeta$. Then, the function $g'$ is continuous with compact support and satisfies $g'(0) = 0$, so that we can define $\mathcal{T}_{g'}$. Furthermore, for every $s, t \in [0, \zeta]$, we have*

$$d_{g'}(s,t) = d_g(\overline{s_0 + s}, \overline{s_0 + t}) \tag{19}$$

*and there exists a unique isometry $R$ from $\mathcal{T}_{g'}$ onto $\mathcal{T}_g$ such that, for every $s \in [0, \zeta]$,*

$$R(p_{g'}(s)) = p_g(\overline{s_0 + s}). \tag{20}$$

Assuming that Theorem 2.2 is proved, we see that $\mathcal{T}_{g'}$ coincides with the real tree $\mathcal{T}_g$ re-rooted at $p_g(s_0)$. Thus the lemma tells us which function codes the tree $\mathcal{T}_g$ re-rooted at an arbitrary vertex.

**Proof.** It is immediately checked that $g'$ satisfies the same assumptions as $g$, so that we can make sense of $\mathcal{T}_{g'}$. Then the key step is to verify the relation (19). Consider first the case where $s, t \in [0, \zeta - s_0[$. Then two possibilities may occur.

If $m_g(s_0 + s, s_0 + t) \geq m_g(s_0, s_0 + s)$, then $m_g(s_0, s_0 + r) = m_g(s_0, s_0 + s) = m_g(s_0, s_0 + t)$ for every $r \in [s, t]$, and so

$$m_{g'}(s,t) = g(s_0) + m_g(s_0 + s, s_0 + t) - 2m_g(s_0, s_0 + s).$$

It follows that

$$
\begin{aligned}
d_{g'}(s,t) &= g'(s) + g'(t) - 2m_{g'}(s,t) \\
&= g(s_0 + s) - 2m_g(s_0, s_0 + s) + g(s_0 + t) - 2m_g(s_0, s_0 + t) \\
&\quad -2(m_g(s_0 + s, s_0 + t) - 2m_g(s_0, s_0 + s)) \\
&= g(s_0 + s) + g(s_0 + t) - 2m_g(s_0 + s, s_0 + t) \\
&= d_g(s_0 + s, s_0 + t).
\end{aligned}
$$

If $m_g(s_0 + s, s_0 + t) < m_g(s_0, s_0 + s)$, then the minimum in the definition of $m_{g'}(s, t)$ is attained at $r_1$ defined as the first $r \in [s, t]$ such that $g(s_0 + r) = m_g(s_0, s_0 + s)$ (because for $r \in [r_1, t]$ we will have $g(s_0 + r) - 2m_g(s_0, s_0 + r) \geq -m_g(s_0, s_0 + r) \geq -m_g(s_0, s_0 + r_1)$). Therefore,

$$m_{g'}(s, t) = g(s_0) - m_g(s_0, s_0 + s),$$

and

$$
\begin{aligned}
d_{g'}(s, t) &= g(s_0 + s) - 2m_g(s_0, s_0 + s) + g(s_0 + t) \\
&\quad - 2m_g(s_0, s_0 + t) + 2m_g(s_0, s_0 + s) \\
&= d_g(s_0 + s, s_0 + t).
\end{aligned}
$$

The other cases are treated in a similar way and are left to the reader.

By (19), if $s, t \in [0, \zeta]$ are such that $d_{g'}(s, t) = 0$, we have $d_g(\overline{s_0 + s}, \overline{s_0 + t}) = 0$ so that $p_g(\overline{s_0 + s}) = p_g(\overline{s_0 + t})$. Noting that $\mathcal{T}_{g'} = p_{g'}([0, \zeta])$ (the supremum of the support of $g'$ is less than or equal to $\zeta$), we can define $R$ in a unique way by the relation (20). From (19), $R$ is an isometry, and it is also immediate that $R$ takes $\mathcal{T}_{g'}$ onto $\mathcal{T}_g$. □

**Proof of Theorem 2.2.** Let us start with some preliminaries. For $\sigma, \sigma' \in \mathcal{T}_g$, we set $\sigma \preccurlyeq \sigma'$ if and only if $d_g(\sigma, \sigma') = d_g(\rho, \sigma') - d_g(\rho, \sigma)$. If $\sigma = p_g(s)$ and $\sigma' = p_g(t)$, it follows from our definitions that $\sigma \preccurlyeq \sigma'$ iff $m_g(s, t) = g(s)$. It is immediate to verify that this defines a partial order on $\mathcal{T}_g$.

For any $\sigma_0, \sigma \in \mathcal{T}_g$, we set

$$[\![\sigma_0, \sigma]\!] = \{\sigma' \in \mathcal{T}_g : d_g(\sigma_0, \sigma) = d_g(\sigma_0, \sigma') + d_g(\sigma', \sigma)\}.$$

If $\sigma = p_g(s)$ and $\sigma' = p_g(t)$, then it is easy to verify that $[\![\rho, \sigma]\!] \cap [\![\rho, \sigma']\!] = [\![\rho, \gamma]\!]$, where $\gamma = p_g(r)$, if $r$ is any time which achieves the minimum of $g$ between $s$ and $t$. We then put $\gamma = \sigma \wedge \sigma'$.

We set $\mathcal{T}_g[\sigma] := \{\sigma' \in \mathcal{T}_g : \sigma \preccurlyeq \sigma'\}$. If $\mathcal{T}_g[\sigma] \neq \{\sigma\}$ and $\sigma \neq \rho$, then $\mathcal{T}_g \backslash \mathcal{T}_g[\sigma]$ and $\mathcal{T}_g[\sigma] \backslash \{\sigma\}$ are two nonempty disjoint open sets. To see that $\mathcal{T}_g \backslash \mathcal{T}_g[\sigma]$ is open, let $s$ be such that $p_g(s) = \sigma$ and note that $\mathcal{T}_g[\sigma]$ is the image under $p_g$ of the compact set $\{u \in [0, \zeta] : m_g(s, u) = g(s)\}$. The set $\mathcal{T}_g[\sigma] \backslash \{\sigma\}$ is open because if $\sigma' \in \mathcal{T}_g[\sigma]$ and $\sigma' \neq \sigma$, it easily follows from our definitions that the open ball centered at $\sigma'$ with radius $d_g(\sigma, \sigma')$ is contained in $\mathcal{T}_g[\sigma] \backslash \{\sigma\}$.

We now prove property (i) of the definition of a real tree. So we fix $\sigma_1$ and $\sigma_2$ in $\mathcal{T}_g$ and we have to prove existence and uniqueness of the mapping $f_{\sigma_1, \sigma_2}$. By using Lemma 2.3 with $s_0$ such that $p_g(s_0) = \sigma_1$, we may assume that $\sigma_1 = \rho$. If $\sigma \in \mathcal{T}_g$ is fixed, we have to prove that there exists a unique isometric mapping $f = f_{\rho, \sigma}$ from $[0, d_g(\rho, \sigma)]$ into $\mathcal{T}_g$ such that $f(0) = \rho$ and $f(d_g(\rho, \sigma)) = \sigma$. Let $s \in p_g^{-1}(\{\sigma\})$, so that $g(s) = d_g(\rho, \sigma)$. Then, for every $a \in [0, d_g(\rho, \sigma)]$, we set

$$v(a) = \inf\{r \in [0, s] : m_g(r, s) = a\}.$$

Note that $g(v(a)) = a$. We put $f(a) = p_g(v(a))$. We have $f(0) = \rho$ and $f(d_g(\rho, \sigma)) = \sigma$, the latter because $m_g(v(g(s)), s) = g(s)$ implies $p_g(v(g(s))) =$ $p_g(s) = \sigma$. It is also easy to verify that $f$ is an isometry: If $a, b \in [0, d_g(\rho, \sigma)]$ with $a \le b$, it is immediate that $m_g(v(a), v(b)) = a$, and so

$$d_g(f(a), f(b)) = g(v(a)) + g(v(b)) - 2a = b - a.$$

To get uniqueness, suppose that $\tilde{f}$ is an isometric mapping satisfying the same properties as $f$. Then, if $a \in [0, d_g(\rho, \sigma)]$,

$$d_g(\tilde{f}(a), \sigma) = d_g(\rho, \sigma) - a = d_g(\rho, \sigma) - d_g(\rho, \tilde{f}(a)).$$

Therefore, $\tilde{f}(a) \preccurlyeq \sigma$. Recall that $\sigma = p_g(s)$, and choose $t$ such that $p_g(t) = \tilde{f}(a)$. Note that $g(t) = d_g(\rho, p_g(t)) = a$. Since $\tilde{f}(a) \preccurlyeq \sigma$ we have $g(t) = m_g(t, s)$. On the other hand, we also know that $a = g(v(a)) = m_g(v(a), s)$. It follows that we have $a = g(t) = g(v(a)) = m_g(v(a), t)$ and thus $d_g(t, v(a)) = 0$, so that $\tilde{f}(a) = p_g(t) = p_g(v(a)) = f(a)$. This completes the proof of (i).

As a by-product of the preceding argument, we see that $f([0, d_g(\rho, \sigma)]) = [[\rho, \sigma]]$: Indeed, we have seen that for every $a \in [0, d_g(\rho, \sigma)]$, we have $f(a) \preccurlyeq \sigma$ and, on the other hand, if $\eta \preccurlyeq \sigma$, the end of the proof of (i) just shows that $\eta = f(d_g(\rho, \eta))$.

We turn to the proof of (ii). We let $q$ be a continuous injective mapping from $[0, 1]$ into $\mathcal{T}_g$, and we aim at proving that $q([0, 1]) = f_{q(0), q(1)}([0, d_g(q(0), q(1))])$. From Lemma 2.3 again, we may assume that $q(0) = \rho$, and we set $\sigma = q(1)$. Then we have just noticed that $f_{0, \sigma}([0, d_g(\rho, \sigma)]) = [[\rho, \sigma]]$.

We first argue by contradiction to prove that $[[\rho, \sigma]] \subset q([0, 1])$. Suppose that $\eta \in [[\rho, \sigma]] \backslash q([0, 1])$, and in particular, $\eta \ne \rho, \sigma$. Then $q([0, 1])$ is contained in the union of the two disjoint open sets $\mathcal{T}_g \backslash \mathcal{T}_g[\eta]$ and $\mathcal{T}_g[\eta] \backslash \{\eta\}$, with $q(0) = \rho \in \mathcal{T}_g[\eta] \backslash \{\eta\}$ and $q(1) = \sigma \in \mathcal{T}_g \backslash \mathcal{T}_g[\eta]$. This contradicts the fact that $q([0, 1])$ is connected.

Conversely, suppose that there exists $a \in (0, 1)$ such that $q(a) \notin [[\rho, \sigma]]$. Set $\eta = q(a)$ and let $\gamma = \sigma \wedge \eta$. Note that $\gamma \in [[\rho, \eta]] \cap [[\rho, \sigma]]$ (from the definition of $\sigma \wedge \eta$, it is immediate to verify that $d_g(\eta, \sigma) = d_g(\eta, \gamma) + d_g(\gamma, \sigma)$). From the first part of the proof of (ii), $\gamma \in q([0, a])$ and, via a root change argument, $\gamma \in q([a, 1])$. Since $q$ is injective, this is only possible if $\gamma = q(a) = \eta$, which contradicts the fact that $\eta \notin [[\rho, \sigma]]$. $\square$

Once we know that $(\mathcal{T}_g, d_g)$ is a real tree, it is straightforward to verify that the notation $\sigma \preccurlyeq \sigma'$, $[[\sigma, \sigma']]$, $\sigma \wedge \sigma'$ introduced in the preceding proof is consistent with the definitions of subsection 2.1 stated for a general real tree.

Let us briefly discuss multiplicities of vertices in the tree $\mathcal{T}_g$. If $\sigma \in \mathcal{T}_g$ is not a leaf then we must have $\ell(\sigma) < r(\sigma)$, where

$$\ell(\sigma) := \inf p_g^{-1}(\{\sigma\}), \quad r(\sigma) := \sup p_g^{-1}(\{\sigma\})$$

are respectively the smallest and the largest element in the equivalence class of $\sigma$ in $[0, \zeta]$. Note that $m_g(\ell(\sigma), r(\sigma)) = g(\ell(\sigma)) = g(r(\sigma)) = d_g(\rho, \sigma)$. Denote by $(a_i, b_i), i \in \mathcal{I}$ the connected components of the open set $(\ell(\sigma), r(\sigma)) \cap \{t \in$ $[0, \infty[: g(t) > d_g(\rho, \sigma)\}$ (the index set $\mathcal{I}$ is empty if $\sigma$ is a leaf). Then we claim that the connected components of the open set $\mathcal{T}_g \backslash \{\sigma\}$ are the sets $p_g((a_i, b_i))$, $i \in \mathcal{I}$ and $\mathcal{T}_g \backslash \mathcal{T}_g[\sigma]$ (the latter only if $\sigma$ is not the root). We have already noticed that $\mathcal{T}_g \backslash \mathcal{T}_g[\sigma]$ is open, and the argument used above for $\mathcal{T}_g[\sigma] \backslash \{\sigma\}$ also shows that the sets $p_g((a_i, b_i))$, $i \in \mathcal{I}$ are open. Finally the sets $p_g((a_i, b_i))$ are connected as continuous images of intervals, and $\mathcal{T}_g \backslash \mathcal{T}_g[\sigma]$ is also connected because if $\sigma', \sigma'' \in \mathcal{T}_g \backslash \mathcal{T}_g[\sigma]$, $[\![\rho, \sigma']\!] \cup [\![\rho, \sigma'']\!]$ is a connected closed set contained in $\mathcal{T}_g \backslash \mathcal{T}_g[\sigma]$.

We conclude this subsection with a lemma comparing the trees coded by two different functions $g$ and $g'$.

**Lemma 2.4** *Let $g$ and $g'$ be two continuous functions with compact support from $[0, \infty[$ into $[0, \infty[$, such that $g(0) = g'(0) = 0$. Then,*

$$d_{GH}(\mathcal{T}_g, \mathcal{T}_{g'}) \leq 2\|g - g'\|,$$

*where $\|g - g'\|$ stands for the uniform norm of $g - g'$.*

**Proof.** We rely on formula (18) for the Gromov-Hausdorff distance. We can construct a correspondence between $\mathcal{T}_g$ and $\mathcal{T}_{g'}$ by setting

$$\mathcal{R} = \{(\sigma, \sigma') : \sigma = p_g(t) \text{ and } \sigma' = p_{g'}(t) \text{ for some } t \geq 0\}.$$

In order to bound the distortion of $\mathcal{R}$, let $(\sigma, \sigma') \in \mathcal{R}$ and $(\eta, \eta') \in \mathcal{R}$. By our definition of $\mathcal{R}$ we can find $s, t \geq 0$ such that $p_g(s) = \sigma$, $p_{g'}(s) = \sigma'$ and $p_g(t) = \eta$, $p_{g'}(t) = \eta'$. Now recall that

$$\begin{aligned} d_g(\sigma, \eta) &= g(s) + g(t) - 2m_g(s, t), \\ d_{g'}(\sigma', \eta') &= g'(s) + g'(t) - 2m_{g'}(s, t), \end{aligned}$$

so that

$$|d_g(\sigma, \eta) - d_{g'}(\sigma', \eta')| \leq 4\|g - g'\|.$$

Thus we have $\text{dis}(\mathcal{R}) \leq 4\|g - g'\|$ and the desired result follows from (18). $\square$

Lemma 2.4 suggests the following alternative proof of Theorem 2.2. Denote by $\mathcal{C}_{00}$ the set of all functions $g : [0, \infty[ \longrightarrow [0, \infty[$ that satisfy the assumptions stated at the beginning of this subsection, and such that the following holds: There exist $\varepsilon > 0$ and $\rho > 0$ such that, for every $i \in \mathbb{N}$, the function $g$ is linear with slope $\rho$ or $-\rho$ over the interval $[(i-1)\varepsilon, i\varepsilon]$. Then it is easy to see that $\mathcal{T}_g$ is a real tree if $g \in \mathcal{C}_{00}$. Indeed, up to a simple time-space rescaling, $g$ will be the contour function of a discrete tree $\mathbf{t} \in \mathbf{A}$, and $\mathcal{T}_g$ coincides (up to rescaling) with the real tree that can be constructed from $\mathbf{t}$ in an obvious way. Then, a general function $g$ can be written as the uniform limit of a sequence $(g_n)$ in $\mathcal{C}_{00}$, and Lemma 2.4 implies that $\mathcal{T}_g$ is the limit in the Gromov-Hausdorff distance of the sequence $\mathcal{T}_{g_n}$. Since each $\mathcal{T}_{g_n}$ is a real tree, $\mathcal{T}_g$ also must be a real tree, by Theorem 2.1 (we do not really need Theorem 2.1, but only the fact that the set of all real trees is closed in the set of all compact metric spaces equipped with the Gromov-Hausdorff metric, cf Lemma 2.1 in [17]).

**Remark.** Recalling that any rooted real tree can be represented in the form $\mathcal{T}_g$, the separability in Theorem 2.1 can be obtained as a consequence of Lemma 2.4 and the separability of the space of continuous functions with compact support on $\mathbb{R}_+$.

### 2.3. The continuum random tree

We recall from Section 1 the notation $\mathbf{e} = (\mathbf{e}_t, 0 \le t \le 1)$ for the normalized Brownian excursion. By convention, we take $\mathbf{e}_t = 0$ if $t > 1$.

**Definition 2.2** *The continuum random tree (CRT) is the random real tree $\mathcal{T}_{\mathbf{e}}$ coded by the normalized Brownian excursion.*

The CRT $\mathcal{T}_{\mathbf{e}}$ is thus a random variable taking values in the set $\mathbb{T}$. Note that the measurability of this random variable follows from Lemma 2.4.

**Remark.** Aldous [1],[2] uses a different method to define the CRT. The preceding definition then corresponds to Corollary 22 in [2]. Note that our normalization differs by an unimportant scaling factor 2 from the one in Aldous' papers: The CRT there is the tree $\mathcal{T}_{2\mathbf{e}}$ instead of $\mathcal{T}_{\mathbf{e}}$.

We can restate many of the results of Section 1 in terms of weak convergence in the space $\mathbb{T}$. Rather than doing this in an exhaustive manner, we will give a typical example showing that the CRT is the limit of rescaled "combinatorial" trees.

Recall from subsection 1.1 the notation $\mathbf{A}$ for the set of all (finite) rooted ordered trees, and denote by $\mathbf{A}_n$ the subset of $\mathbf{A}$ consisting of trees with $n$ vertices. We may and will view each element $\mathbf{t}$ of $\mathbf{A}$ as a rooted real tree: Simply view $\mathbf{t}$ as a union of line segments of length 1 in the plane, in the way represented in the left part of Figure 1, equipped with the obvious distance (the distance between $\sigma$ and $\sigma'$ is the length of the shortest path from $\sigma$ to $\sigma'$ in the tree). Alternatively, if $(C_t, t \ge 0)$ is the contour function of the tree, this means that we identify $\mathbf{t} = \mathcal{T}_C$ (this is not really an identification, because the tree $\mathbf{t}$ has an order structure which disappears when we consider it as a real tree).

For any $\lambda > 0$ and a tree $\mathcal{T} \in \mathbb{T}$, the tree $\lambda\mathcal{T}$ is the "same" tree with all distances multiplied by the factor $\lambda$ (if the tree is embedded in the plane as suggested above, this corresponds to replacing the set $\mathcal{T}$ by $\lambda\mathcal{T}$).

**Theorem 2.5** *For every $n \ge 1$, let $\mathcal{T}_{(n)}$ be a random tree distributed uniformly over $\mathbf{A}_n$. Then $(2n)^{-1/2}\mathcal{T}_{(n)}$ converges in distribution to the CRT $\mathcal{T}_{\mathbf{e}}$, in the space $\mathbb{T}$.*

**Proof.** Let $\theta$ be a Galton-Watson tree with geometric offspring distribution $\mu(k) = 2^{-k-1}$, and for every $n \ge 1$ let $\theta_n$ be distributed as $\theta$ conditioned to have $n$ vertices. Then it is easy to verify that $\theta_n$ has the same distribution as $\mathcal{T}_{(n)}$. On the other hand, let $(C_t^n, t \ge 0)$ be the contour function of $\theta_n$, and let

$$\widetilde{C}_t^n = (2n)^{-1/2} C_{2nt}^n , \quad t \ge 0.$$

From Theorem 1.15 (restated in terms of the contour function as explained in subsection 1.6), we have

$$(\widetilde{C}_t^n, t \geq 0) \xrightarrow[n \to \infty]{\text{(d)}} (\mathbf{e}_t, t \geq 0).$$

On the other hand, from the observations preceding the theorem, and the fact that $\theta_n$ has the same distribution as $\mathcal{T}_{(n)}$, it is immediate that the tree $\mathcal{T}_{\widetilde{C}^n}$ coded by $\widetilde{C}^n$ has the same distribution as $(2n)^{-1/2}\mathcal{T}_{(n)}$. The statement of the theorem now follows from the previous convergence and Lemma 2.4. $\qquad \square$

We could state analogues of Theorem 2.5 for several other classes of combinatorial trees, such as the ones considered at the end of subsection 1.5. For instance, if $\tau_n$ is distributed uniformly among all rooted Cayley trees with $n$ vertices , then $(4n)^{-1/2}\tau_n$ converges in distribution to the CRT $\mathcal{T}_\mathbf{e}$, in the space $\mathbb{T}$. Notice that Cayley trees are not ordered, but that they can be obtained from (ordered) Galton-Watson trees with Poisson offspring distribution by "forgetting" the order, as was explained at the end of subsection 1.5. By applying the same argument as in the preceding proof to these (conditioned) Galton-Watson trees, we get the desired convergence for rescaled Cayley trees.

### 2.4. The Itô excursion measure

Our goal is to derive certain explicit distributions for the CRT, and more specifically its so-called finite-dimensional marginal distributions. For these calculations, we will need some basic properties of Brownian excursions. Before dealing with the normalized Brownian excursion, we will consider Itô's measure.

We denote by $(B_t, t \geq 0)$ a linear Brownian motion, which starts at $x$ under the probability measure $P_x$. We set

$$S_t = \sup_{s \leq t} B_s \ , \qquad I_t = \inf_{s \leq t} B_s$$

and, for every $a \in \mathbb{R}$, $T_a = \inf\{t \geq 0 : B_t = a\}$. The reflection principle gives the law of the pair $(S_t, B_t)$: If $a \geq 0$ and $b \in (-\infty, a]$,

$$P_0[S_t \geq a, B_t \leq b] = P_0[B_t \geq 2a - b].$$

It follows that, for every $t > 0$, the density under $P_0$ of the law of the pair $(S_t, B_t)$ is

$$\gamma_t(a, b) = \frac{2(2a - b)}{\sqrt{2\pi t^3}} \exp\left(-\frac{(2a - b)^2}{2t}\right) 1_{\{a \geq 0, b \leq a\}}. \tag{21}$$

The reflection principle also implies that $S_t$ and $|B_t|$ have the same distribution. Let $a > 0$. Observing that $\{T_a \leq t\} = \{S_t \geq a\}$, $P_0$ a.s., we obtain that the density of $T_a$ under $P_0$ is the function

$$q_a(t) = \frac{a}{\sqrt{2\pi t^3}} \exp\left(-\frac{a^2}{2t}\right).$$

Notice the relation $\gamma_t(a, b) = 2\, q_{2a-b}(t)$ for $a \geq 0$ and $b < a$.

For every $\varepsilon > 0$, denote by $\nu_\varepsilon$ the law of the first excursion of $B$ away from 0 that hits level $\varepsilon$. More specifically, if $G_\varepsilon = \sup\{t < T_\varepsilon : B_t = 0\}$ and $D_\varepsilon = \inf\{t > T_\varepsilon : B_t = 0\}$, $\nu_\varepsilon$ is the law of $(B_{(G_\varepsilon + t) \wedge D_\varepsilon}, t \geq 0)$. The measure $\nu_\varepsilon$ is thus a probability measure on the set $\mathcal{C} = C(\mathbb{R}_+, \mathbb{R}_+)$ of all continuous functions from $\mathbb{R}_+$ into $\mathbb{R}_+$, and is supported on $\mathcal{C}_\varepsilon = \{e \in \mathcal{C} : \sup e(s) \geq \varepsilon\}$. If $0 < \varepsilon < \varepsilon'$, we have

$$\nu_\varepsilon(\mathcal{C}_{\varepsilon'}) = P_\varepsilon[T_{\varepsilon'} < T_0] = \frac{\varepsilon}{\varepsilon'}$$

and

$$\nu_{\varepsilon'} = \nu_\varepsilon(\cdot \mid \mathcal{C}_{\varepsilon'}) = \frac{\varepsilon'}{\varepsilon} \nu_\varepsilon(\cdot \cap \mathcal{C}_{\varepsilon'}).$$

For every $\varepsilon > 0$, set

$$n_\varepsilon = \frac{1}{2\varepsilon}\, \nu_\varepsilon.$$

Then $n_{\varepsilon'} = n_\varepsilon(\cdot \cap \mathcal{C}_{\varepsilon'})$ for every $0 < \varepsilon < \varepsilon'$. This leads to the following definition.

**Definition 2.3** *The $\sigma$-finite measure $n$ on $\mathcal{C}$ defined by*

$$n = \lim_{\varepsilon \downarrow 0} \uparrow n_\varepsilon$$

*is called the Itô measure of positive excursions of linear Brownian motion.*

Let us briefly state some simple properties of the Itô measure. First $n(de)$ is supported on the set $\mathcal{E}$ consisting of all elements $e \in \mathcal{C}$ which have the property that there exists $\sigma = \sigma(e) > 0$ such that $e(t) > 0$ if and only if $0 < t < \sigma$ (the number $\sigma(e)$ is called the length, or the duration of excursion $e$). By construction, $n_\varepsilon$ is the restriction of $n$ to $\mathcal{C}_\varepsilon$, and in particular $n(\mathcal{C}_\varepsilon) = (2\varepsilon)^{-1}$. Finally, if $T_\varepsilon(e) = \inf\{t \geq 0 : e(t) = \varepsilon\}$, the law of $(e(T_\varepsilon(e) + t), t \geq 0)$ under $n(\cdot \mid T_\varepsilon < \infty) = \nu_\varepsilon$ is the law of $(B_{t \wedge T_0}, t \geq 0)$ under $P_\varepsilon$. The last property follows from the construction of the measure $\nu_\varepsilon$ and the strong Markov property of Brownian motion at time $T_\varepsilon$.

**Proposition 2.6** (i) *For every $t > 0$, and every measurable function $g : \mathbb{R}_+ \longrightarrow \mathbb{R}_+$ such that $g(0) = 0$,*

$$\int n(de)\, g(e(t)) = \int_0^\infty dx\, q_x(t)\, g(x). \tag{22}$$

*In particular, $n(\sigma > t) = n(e(t) > 0) = (2\pi t)^{-1/2} < \infty$. Moreover,*

$$n\Big(\int_0^\infty dt\, g(e(t))\Big) = \int_0^\infty dx\, g(x). \tag{23}$$

(ii) *Let $t > 0$ and let $\Phi$ and $\Psi$ be two nonnegative measurable functions defined respectively on $C([0, t], \mathbb{R}_+)$ and $\mathcal{C}$. Then,*

$$\int n(de)\, \Phi(e(r), 0 \leq r \leq t) \Psi(e(t + r), r \geq 0)$$

$$= \int n(de)\, \Phi(e(r), 0 \leq r \leq t)\, E_{e(t)}\big[\Psi(B_{r \wedge T_0}, r \geq 0)\big].$$

**Proof.** (i) We may assume that $g$ is bounded and continuous and that there exists $\alpha > 0$ such that $g(x) = 0$ if $x \leq \alpha$. Then, by dominated convergence,

$$\int n(de)\, g(e(t)) = \lim_{\varepsilon \downarrow 0} \int n(de)\, g(T_\varepsilon(e) + t)\, 1_{\{T_\varepsilon(e) < \infty\}} = \lim_{\varepsilon \downarrow 0} \frac{1}{2\varepsilon} E_\varepsilon[g(B_{t \wedge T_0})],$$

using the property stated just before the proposition. From formula (21) we get

$$\begin{aligned} E_\varepsilon[g(B_{t \wedge T_0})] = E_\varepsilon[g(B_t)\, 1_{\{t < T_0\}}] &= E_0[g(\varepsilon - B_t)\, 1_{\{S_t < \varepsilon\}}] \\ &= \int_0^\varepsilon da \int_{-\infty}^a db\, g(\varepsilon - b)\, \gamma_t(a, b). \end{aligned}$$

The first assertion in (i) now follows, observing that $q_x(t) = \frac{1}{2}\gamma_t(0, -x)$. The identity $n(e(t) > 0) = (2\pi t)^{-1/2} < \infty$ is obtained by taking $g(x) = 1_{\{x > 0\}}$. The last assertion in (i) follows from (22), recalling that the function $t \to q_x(t)$ is a probability density.

(ii) We may assume that $\Phi$ and $\Psi$ are bounded and continuous and that there exists $\alpha \in (0, t)$ such that $\Phi(\omega(r), 0 \leq r \leq t) = 0$ if $\omega(\alpha) = 0$. The proof then reduces to an application of the Markov property of Brownian motion stopped at time $T_0$, by writing

$$\int n(de)\, \Phi(e(r), 0 \leq r \leq t)\Psi(e(t + r), r \geq 0)$$

$$= \lim_{\varepsilon \to 0} \int n(de)\, 1_{\{T_\varepsilon(e) < \infty\}}\, \Phi(e(T_\varepsilon(e) + r), 0 \leq r \leq t)\Psi(e(T_\varepsilon(e) + t + r), r \geq 0)$$

$$= \lim_{\varepsilon \to 0} \frac{1}{2\varepsilon} E_\varepsilon\Big[\Phi(B_{r \wedge T_0}, 0 \leq r \leq t)\, \Psi(B_{(t+s) \wedge T_0}, s \geq 0)\Big]$$

$$= \lim_{\varepsilon \to 0} \frac{1}{2\varepsilon} E_\varepsilon\Big[\Phi(B_{r \wedge T_0}, 0 \leq r \leq t)\, E_{B_{t \wedge T_0}}[\Psi(B_{s \wedge T_0}, s \geq 0)]\Big]$$

$$= \lim_{\varepsilon \to 0} \int n(de)\, 1_{\{T_\varepsilon(e) < \infty\}}\, \Phi(e(T_\varepsilon(e) + r), 0 \leq r \leq t)E_{e(T_\varepsilon(e)+t)}[\Psi(B_{s \wedge T_0}, s \geq 0)]$$

$$= \int n(de)\, \Phi(e(r), 0 \leq r \leq t)\, E_{e(t)}\big[\Psi(B_{r \wedge T_0}, r \geq 0)\big].$$

In the first and in the last equality, dominated convergence is justified by the fact that $1_{\{T_\varepsilon(e) < \infty\}}\Phi(e(T_\varepsilon(e) + r), 0 \leq r \leq t) = 0$ if $\sigma(e) \leq \alpha$, and the property $n(\sigma > \alpha) < \infty$. $\qquad\square$

### 2.5. *Finite-dimensional marginals under the Itô measure*

If $(\mathcal{T}, d)$ is a real tree with root $\rho$, and if $x_1, \ldots, x_p \in \mathcal{T}$, the subtree spanned by $x_1, \ldots, x_p$ is simply the set

$$\mathcal{T}(x_1, \ldots, x_p) = \bigcup_{i=1}^p [\![\rho, x_i]\!].$$

It is easy to see that $\mathcal{T}(x_1, \ldots, x_p)$, equipped with the distance $d$, is again a real tree, which has a discrete structure: More precisely, $\mathcal{T}(x_1, \ldots, x_p)$ can be represented by a discrete skeleton (which is a discrete rooted tree with $p$ labelled leaves) and the collection, indexed by vertices of the skeleton, of lengths of "branches".

Rather than giving formal definitions for a general real tree, we will concentrate on the case of the tree $\mathcal{T}_g$ coded by $g$ in the sense of subsection 2.2.

Recall that $\mathcal{C}$ denotes the set of all continuous functions from $\mathbb{R}_+$ into $\mathbb{R}_+$. We consider a general continuous function $g \in \mathcal{C}$ (in contrast with subsection 2.2, we do not assume that $g(0) = 0$ and $g$ has compact support). If $0 \le t_1 \le t_2 \le \cdots \le t_p$, we will define a "marked tree"

$$\theta(g; t_1, \ldots, t_p) = \Big( \tau(g; t_1, \ldots, t_p), (h_v)_{v \in \tau(g; t_1, \ldots, t_p)} \Big)$$

where $\tau(g; t_1, \ldots, t_p) \in \mathbf{A}$ ($\mathbf{A}$ is the set of all rooted ordered trees as in Section 1) and $h_v \ge 0$ for every $v \in \tau(g; t_1, \ldots, t_p)$. An example with $p = 3$ is given in Figure 3 below. In this example, $\tau(g; t_1, t_2, t_3) = \{\varnothing, 1, 2, (2, 1), (2, 2)\}$ as pictured in the right part of Figure 3, and the numbers $h_v, v \in \tau(g; t_1, t_2, t_3)$ are the lengths of the bold segments as indicated on the left part of the figure.

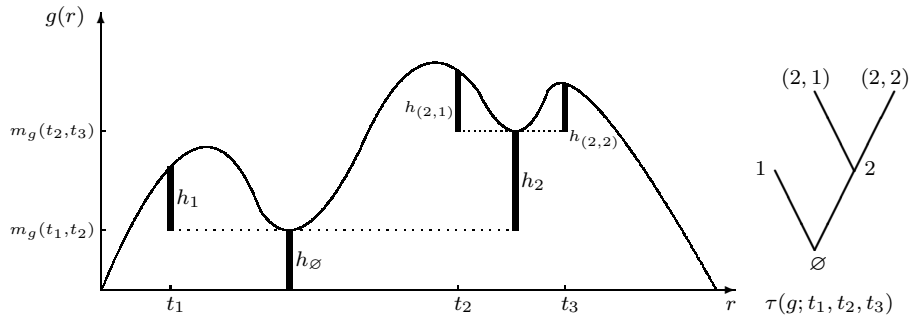

Figure 3

To give a precise definition of $\theta(g; t_1, \ldots, t_p)$, we proceed by induction on $p$.

If $p = 1$, $\tau(g; t_1) = \{\varnothing\}$ and $h_\varnothing(g; t_1) = g(t_1)$.

Let $p \ge 2$ and suppose that the "marked tree" $\theta(g; t_1, \ldots, t_j)$ has been constructed up to order $p - 1$. Then there exists an integer $k \in \{1, \ldots, p - 1\}$ and $k$ integers $1 \le i_1 < i_2 < \cdots < i_k \le p - 1$ such that $m_g(t_i, t_{i+1}) = m_g(t_1, t_p)$ iff $i \in \{i_1, \ldots, i_k\}$. Set $i_0 = 0$ and $i_{k+1} = p$ by convention. For every $\ell \in \{1, 2, \ldots, k + 1\}$, define $g^\ell \in \mathcal{C}$ by the formula

$$g^\ell(t) = g((t \vee t_{i_{\ell-1}+1}) \wedge t_{i_\ell}) - m_g(t_1, t_p).$$

We then let $\tau(g; t_1, \ldots, t_p)$ be the concatenation of the trees $\tau(g^\ell; t_{i_{\ell-1}+1}, \ldots, t_{i_\ell})$ for $1 \le \ell \le k + 1$: Precisely,

$$\tau(g; t_1, \ldots, t_p) = \{\varnothing\} \cup \bigcup_{\ell=1}^{k+1} \{\ell u : u \in \tau(g^\ell; t_{i_{\ell-1}+1}, \ldots, t_{i_\ell})\}.$$

Furthermore, if for $\ell \in \{1, \dots, k+1\}$

$$\theta(g^\ell; t_{i_{\ell-1}+1}, \dots, t_{i_\ell}) = (\tau(g^\ell; t_{i_{\ell-1}+1}, \dots, t_{i_\ell}), (h_v^\ell)_{v \in \tau(g^\ell; t_{i_{\ell-1}+1}, \dots, t_{i_\ell})}),$$

we define the marks $(h_v)_{v \in \tau(g; t_1, \dots, t_p)}$ by setting

$$h_v = h_u^\ell, \text{ if } v = \ell u, \ u \in \tau(g^\ell; t_{i_{\ell-1}+1}, \dots, t_{i_\ell}),$$

and $h_\varnothing = m_g(t_1, t_p)$

This completes the construction of the tree by induction. Note that $k+1$ is the number of children of $\varnothing$ in the tree $\theta(g; t_1, \dots, t_p)$, and $m(t_1, t_p)$ is the mark of $\varnothing$.

If now $g$ satisfies the conditions in subsection 2.2, it is easy to see that $\theta(g; t_1, \dots, t_p)$ corresponds to the tree $\mathcal{T}_g(p_g(t_1), \dots, p_g(t_p))$ spanned by the vertices $p_g(t_1), \dots, p_g(t_p)$ in the tree $\mathcal{T}_g$. More precisely, if we attach to every $v \in \tau(g; t_1, \dots, t_p)$ a line segment in the plane with length $h_v$, in such a way that the line segments attached to $v$ and to its children share a common end (the same for all children of $v$) and that the line segments otherwise do not intersect, the union of the resulting line segments will give a representative of the equivalence class of $\mathcal{T}_g(p_g(t_1), \dots, p_g(t_p))$. (Note that the order structure of $\tau(g; t_1, \dots, t_p)$ plays no role in this construction.)

We let $\mathbf{A}_{(p)}$ be the set of all rooted ordered trees with $p$ leaves (a leaf of a tree $\tau \in \mathbf{A}$ is a vertex $u \in \tau$ with no child, i.e. such that $k_u(\tau) = 0$, with the notation of Section 1). We denote by $\Theta_{(p)}$ the set of all marked trees with $p$ leaves: Elements of $\Theta_{(p)}$ are of the form $\theta = (\tau, (h_v)_{v \in \tau})$ where $\tau \in \mathbf{A}_{(p)}$ and $h_v \geq 0$ for every $v \in \tau$. The set $\Theta_{(p)}$ is equipped with the obvious topology and the associated Borel $\sigma$-field. We also consider the corresponding sets of binary trees: $\mathbf{A}_{(p)}^{\mathrm{bin}}$ is the set of all binary rooted trees with $p$ leaves (and hence $2p - 1$ vertices), and $\Theta_{(p)}^{\mathrm{bin}}$ is the set of marked trees $\theta = (\tau, (h_v)_{v \in \tau})$ whose skeleton $\tau$ belongs to $\mathbf{A}_{(p)}^{\mathrm{bin}}$. Recall that

$$\#(\mathbf{A}_{(p)}^{\mathrm{bin}}) = \frac{(2p-2)!}{(p-1)! \, p!} =: c_p$$

is the Catalan number of order $p - 1$.

**Theorem 2.7** *The law of the tree $\theta(e; t_1, \dots, t_p)$ under the measure*

$$n(de) \, 1_{\{0 \leq t_1 \leq \cdots \leq t_p \leq \sigma(e)\}} dt_1 \dots dt_p$$

*is $2^{p-1} \Lambda_p$, where $\Lambda_p$ is the uniform measure on $\Theta_{(p)}^{\mathrm{bin}}$, defined by*

$$\int \Lambda_p(d\theta) \, F(\theta) = \sum_{\tau \in \mathbf{A}_{(p)}^{\mathrm{bin}}} \int \prod_{v \in \tau} dh_v \, F(\tau, \{h_v, v \in \tau\}).$$

It should be clear from our construction that the tree $\theta(e; t_1, \dots, t_p)$ only depends on the values of $e(t_1), \dots, e(t_p)$ and of the successive minima $m_e(t_1, t_2)$, $m_e(t_2, t_3), \dots, m_e(t_{p-1}, t_p)$. The key tool in the proof of Theorem 2.7 will thus be the following proposition. To simplify notation we write $m(s, t) = m_e(s, t)$.

**Proposition 2.8** *Let $F$ be nonnegative and measurable on $\mathbb{R}_+^{2p-1}$. Then*

$$n\Big(\int_{\{0 \leq t_1 \leq \cdots \leq t_p \leq \sigma\}} dt_1 \ldots dt_p \, f(m(t_1, t_2), \ldots, m(t_{p-1}, t_p), e(t_1), \ldots, e(t_p))\Big)$$

$$= 2^{p-1} \int_{\mathbb{R}_+^{2p-1}} d\alpha_1 .. d\alpha_{p-1} d\beta_1 .. d\beta_p \Big(\prod_{i=1}^{p-1} 1_{\{\alpha_i \leq \beta_i \wedge \beta_{i+1}\}}\Big) f(\alpha_1, .., \alpha_{p-1}, \beta_1, .., \beta_p).$$

Before proving this proposition, we state a lemma which is an immediate consequence of (21). Recall that $B$ is a Brownian motion that starts from $x$ under the probability $P_x$, and that $I = (I_t, t \geq 0)$ is the associated minimum process.

**Lemma 2.9** *If $g$ is a nonnegative measurable function on $\mathbb{R}^3$ and $x \geq 0$,*

$$E_x\Big(\int_0^{T_0} dt \, g(t, I_t, B_t)\Big) = 2 \int_0^x dy \int_y^\infty dz \int_0^\infty dt \, q_{x+z-2y}(t) \, g(t, y, z) \quad (24)$$

*In particular, if $h$ is a nonnegative measurable function on $\mathbb{R}^2$,*

$$E_x\Big(\int_0^{T_0} dt \, h(I_t, B_t)\Big) = 2 \int_0^x dy \int_y^\infty dz \, h(y, z). \quad (25)$$

**Proof of Proposition 2.8.** This is a simple consequence of Lemma 2.9. For $p = 1$, the result is exactly formula (23) in Proposition 2.6. We proceed by induction on $p$ using the Markov property under $n$ (property (ii) in Proposition 2.6) and then (25):

$$n\Big(\int_{\{0 \leq t_1 \leq \cdots \leq t_p \leq \sigma\}} dt_1 \ldots dt_p \, f(m(t_1, t_2), \ldots, m(t_{p-1}, t_p), e(t_1), \ldots, e(t_p))\Big)$$

$$= n\Big(\int_{\{0 \leq t_1 \leq \cdots \leq t_{p-1} \leq \sigma\}} dt_1 \ldots dt_{p-1}$$

$$E_{e(t_{p-1})}\Big[\int_0^{T_0} dt \, f(m(t_1, t_2), \ldots, m(t_{p-2}, t_{p-1}), I_t, e(t_1), \ldots, e(t_{p-1}), B_t)\Big]\Big)$$

$$= 2 \, n\Big(\int_{\{0 \leq t_1 \leq \cdots \leq t_{p-1} \leq \sigma\}} dt_1 \ldots dt_{p-1}$$

$$\int_0^{e(t_{p-1})} d\alpha_{p-1} \int_{\alpha_{p-1}}^\infty d\beta_p \, f(m(t_1, t_2), .., m(t_{p-2}, t_{p-1}), \alpha_{p-1}, e(t_1), .., e(t_{p-1}), \beta_p)\Big).$$

The proof is then completed by using the induction hypothesis. □

**Proof of Theorem 2.7.** Let $\Gamma_p$ be the measurable function from $\mathbb{R}_+^{2p-1}$ into $\Theta_{(p)}$ such that

$$\theta(e; t_1, \ldots, t_p) = \Gamma_p(m(t_1, t_2), \ldots, m(t_{p-1}, t_p), e(t_1), \ldots, e(t_p)).$$

The existence of this function easily follows from our construction by induction of the marked tree $\theta(e; t_1, \ldots, t_p)$.

Denote by $\Delta_p$ the measure on $\mathbb{R}_+^{2p-1}$ defined by

$$\Delta_p(d\alpha_1 \dots d\alpha_{p-1}d\beta_1 \dots d\beta_p) = \Big(\prod_{i=1}^{p-1} 1_{[0,\beta_i \wedge \beta_{i+1}]}(\alpha_i)\Big)d\alpha_1 \dots d\alpha_{p-1}d\beta_1 \dots d\beta_p.$$

In view of Proposition 2.8, the proof of Theorem 2.7 reduces to checking that $\Gamma_p(\Delta_p) = \Lambda_p$. For $p = 1$, this is obvious.

Let $p \geq 2$ and suppose that the result holds up to order $p - 1$. For every $j \in \{1, \dots, p-1\}$, let $H_j$ be the subset of $\mathbb{R}_+^{2p-1}$ defined by

$$H_j = \{(\alpha_1, \dots, \alpha_{p-1}, \beta_1, \dots, \beta_p); \alpha_i > \alpha_j \text{ for every } i \neq j\}.$$

Then,

$$\Delta_p = \sum_{j=1}^{p-1} 1_{H_j} \cdot \Delta_p.$$

On the other hand, it is immediate to verify that $1_{H_j} \cdot \Delta_p$ is the image of the measure

$$\Delta_j(d\alpha_1' \dots d\beta_j') \otimes 1_{(0,\infty)}(h)dh \otimes \Delta_{p-j}(d\alpha_1'' \dots d\beta_{p-j}'')$$

under the mapping $\Phi : (\alpha_1', \dots, \beta_j', h, \alpha_1'' \dots, \beta_{p-j}'') \longrightarrow (\alpha_1, \dots, \beta_p)$ defined by

$$\begin{aligned}
\alpha_j &= h, \\
\alpha_i &= \alpha_i' + h &&\text{for } 1 \leq i \leq j-1, \\
\beta_i &= \beta_i' + h &&\text{for } 1 \leq i \leq j, \\
\alpha_i &= \alpha_{i-j}'' + h &&\text{for } j+1 \leq i \leq p-1, \\
\beta_i &= \beta_{i-j}'' + h &&\text{for } j+1 \leq i \leq p.
\end{aligned}$$

The construction by induction of the tree $\theta(e; t_1, \dots, t_p)$ exactly shows that, a.e. for the measure $\Delta_j(d\alpha_1' \dots d\beta_j') \otimes 1_{(0,\infty)}(h)dh \otimes \Delta_{p-j}(d\alpha_1'' \dots d\beta_{p-j}'')$, we have

$$\Gamma_p \circ \Phi(\alpha_1', \dots, \beta_j', h, \alpha_1'' \dots, \beta_{p-j}'') = \Gamma_j(\alpha_1', \dots, \beta_j') \underset{h}{*} \Gamma_{p-j}(\alpha_1'' \dots, \beta_{p-j}'').$$

where if $\theta \in \Theta_{(j)}$ and $\theta' \in \Theta_{(p-j)}$, the tree $\theta \underset{h}{*} \theta'$ is obtained by concatenating the discrete skeletons of $\theta$ and $\theta'$ (as above in the construction by induction of $\theta(g; t_1, \dots, t_p)$) and assigning the mark $h$ to the root $\varnothing$.

Together with the induction hypothesis, the previous observations imply that for any nonnegative measurable function $f$ on $\Theta_{(p)}$,

$$\begin{aligned}
\int \Delta_p(du) \, 1_{H_j}(u) \, f(\Gamma_p(u)) &= \int_0^\infty dh \iint \Delta_j(du')\Delta_{p-j}(du'')f\big(\Gamma_p(\Phi(u', h, u''))\big) \\
&= \int_0^\infty dh \iint \Delta_j(du')\Delta_{p-j}(du'')f\big(\Gamma_j(u') \underset{h}{*} \Gamma_{p-j}(u'')\big) \\
&= \int_0^\infty dh \int \Lambda_j \underset{h}{*} \Lambda_{p-j}(d\theta) \, f(\theta)
\end{aligned}$$

where we write $\Lambda_j \underset{h}{*} \Lambda_{p-j}$ for the image of $\Lambda_j(d\theta)\Lambda_{p-j}(d\theta')$ under the mapping $(\theta, \theta') \longrightarrow \theta \underset{h}{*} \theta'$. To complete the proof, simply note that

$$\Lambda_p = \sum_{j=1}^{p-1} \int_0^\infty dh \, \Lambda_j \underset{h}{*} \Lambda_{p-j}.$$

$\square$

**Remark.** The fact that we get only binary trees in Theorem 2.7 corresponds to the property that local minima of Brownian motion are distinct: If $0 < t_1 < \cdots < t_p$ and if the local minima $m_g(t_i, t_{i+1})$ are distinct, the tree $\theta(g; t_1, \ldots, t_p)$ is clearly binary.

### 2.6. Finite-dimensional marginals of the CRT

In this subsection, we propose to calculate the law of the tree $\theta(e; t_1, \ldots, t_p)$ when $e$ is a normalized Brownian excursion. This corresponds to choosing $p$ vertices independently uniformly on the CRT (the uniform measure on the CRT $\mathcal{T}_e$ is by definition the image of Lebesgue measure on $[0, 1]$ under the mapping $p_e$) and determining the law of the tree spanned by these vertices. In contrast with the measure $\Lambda_p$ of Theorem 2.7, we will get for every $p$ a *probability* measure on $\Theta_{(p)}^{\text{bin}}$, which may be called the $p$-dimensional marginal distribution of the CRT (cf Aldous [1], [2]).

We first recall the connection between the Itô measure and the normalized Brownian excursion. Informally, the law of the normalized Brownian excursion ($\mathbf{e}$ in the notation of Section 1) is $n(de \mid \sigma(e) = 1)$. More precisely, using a standard desintegration theorem for measures, together with the Brownian scaling property, one easily shows that there exists a unique collection of probability measures $(n_{(s)}, s > 0)$ on the set $\mathcal{E}$ of excursions, such that the following properties hold:

(i) For every $s > 0$, $n_{(s)}(\sigma = s) = 1$.

(ii) For every $\lambda > 0$ and $s > 0$, the law under $n_{(s)}(de)$ of $e_\lambda(t) = \sqrt{\lambda} \, e(t/\lambda)$ is $n_{(\lambda s)}$.

(iii) For every Borel subset $A$ of $\mathcal{E}$,

$$n(A) = \frac{1}{2}(2\pi)^{-1/2} \int_0^\infty s^{-3/2} \, n_{(s)}(A) \, ds.$$

The measure $n_{(1)}$ is the law of the normalized Brownian excursion $\mathbf{e}$ which was considered in Section 1 and in subsection 2.3 above.

Our first goal is to get a statement more precise than Theorem 2.7 by considering the pair $(\theta(e; t_1, \ldots, t_p), \sigma)$ instead of $\theta(e; t_1, \ldots, t_p)$. If $\theta = (\tau, \{h_v, v \in \tau\})$ is a marked tree, the length of $\theta$ is defined in the obvious way by

$$L(\theta) = \sum_{v \in \tau} h_v.$$

**Proposition 2.10** *The law of the pair* $(\theta(e; t_1, \ldots, t_p), \sigma)$ *under the measure*

$$n(de) \, 1_{\{0 \leq t_1 \leq \cdots \leq t_p \leq \sigma(e)\}} dt_1 \ldots dt_p$$

*is*

$$2^{p-1} \, \Lambda_p(d\theta) \, q_{2L(\theta)}(s) ds.$$

**Proof.** Recall the notation of the proof of Theorem 2.7. We will verify that, for any nonnegative measurable function $f$ on $\mathbb{R}_+^{3p}$,

$$n\Big(\int_{\{0 \leq t_1 \leq \cdots \leq t_p \leq \sigma\}} dt_1 \ldots dt_p$$
$$f\big(m(t_1, t_2), \ldots, m(t_{p-1}, t_p), e(t_1), \ldots, e(t_p), t_1, t_2 - t_1, \ldots, \sigma - t_p\big)\Big)$$
$$= 2^{p-1} \int \Delta_p(d\alpha_1 \ldots d\alpha_{p-1} d\beta_1 \ldots d\beta_p) \int_{\mathbb{R}_+^{p+1}} ds_1 \ldots ds_{p+1} \, q_{\beta_1}(s_1) q_{\beta_1 + \beta_2 - 2\alpha_1}(s_2) \ldots$$
$$\ldots q_{\beta_{p-1} + \beta_p - 2\alpha_{p-1}}(s_p) q_{\beta_p}(s_{p+1}) \, f(\alpha_1, \ldots, \alpha_{p-1}, \beta_1, \ldots, \beta_p, s_1, \ldots, s_{p+1}). \quad (26)$$

Suppose that (26) holds. It is easy to check (for instance by induction on $p$) that

$$2 \, L(\Gamma_p(\alpha_1, \ldots, \alpha_{p-1}, \beta_1, \ldots, \beta_p)) = \beta_1 + \sum_{i=1}^{p-1} (\beta_i + \beta_{i-1} - 2\alpha_i) + \beta_p.$$

Using the convolution identity $q_x * q_y = q_{x+y}$ (which is immediate from the interpretation of $q_x$ as the density of $T_x$ under $P_0$), we get from (26) that

$$n\Big(\int_{\{0 \leq t_1 \leq \cdots \leq t_p \leq \sigma\}} dt_1 \ldots dt_p \, f\big(m(t_1, t_2), \ldots, m(t_{p-1}, t_p), e(t_1), \ldots, e(t_p), \sigma\big)\Big)$$
$$= 2^{p-1} \int \Delta_p(d\alpha_1 \ldots d\alpha_{p-1} d\beta_1 \ldots d\beta_p) \int_0^\infty dt \, q_{2L(\Gamma_p(\alpha_1, \ldots, \beta_p))}(t) \, f(\alpha_1, \ldots, \beta_p, t).$$

As in the proof of Theorem 2.7, the statement of Proposition 2.10 follows from this last identity and the equality $\Gamma_p(\Delta_p) = \Lambda_p$.

It remains to prove (26). The case $p = 1$ is easy: By using the Markov property under the Itô measure (Proposition 2.6 (ii)), then the definition of the function $q_x$ and finally (22), we get

$$\int n(de) \int_0^\sigma dt \, f(e(t), t, \sigma - t) = \int n(de) \int_0^\sigma dt \, E_{e(t)}\big[f(e(t), t, T_0)\big]$$
$$= \int n(de) \int_0^\sigma dt \int_0^\infty dt' \, q_{e(t)}(t') f(e(t), t, t')$$
$$= \int_0^\infty dx \int_0^\infty dt \, q_x(t) \int_0^\infty dt' \, q_x(t') \, f(x, t, t').$$

Let $p \geq 2$. Applying the Markov property under $n$ successively at $t_p$ and at $t_{p-1}$, and then using (24), we obtain

$$n\Big(\int_{\{0 \leq t_1 \leq \cdots \leq t_p \leq \sigma\}} dt_1 \ldots dt_p$$

$$\times f\big(m(t_1,t_2),\ldots,m(t_{p-1},t_p),e(t_1),\ldots,e(t_p),t_1,t_2-t_1,\ldots,\sigma-t_p\big)\Big)$$

$$= n\Big(\int_{\{0 \leq t_1 \leq \cdots \leq t_{p-1} \leq \sigma\}} dt_1 \ldots dt_{p-1} \, E_{e(t_{p-1})}\Big[\int_0^{T_0} dt \int_0^\infty ds \, q_{B_t}(s)$$

$$\times f\big(m(t_1,t_2),..,m(t_{p-2},t_{p-1}),I_t,e(t_1),..,e(t_{p-1}),B_t,t_1,..,t_{p-1}-t_{p-2},t,s\big)\Big]\Big)$$

$$= 2n\Big(\int_{\{0 \leq t_1 \leq \cdots \leq t_{p-1} \leq \sigma\}} dt_1 \ldots dt_{p-1}$$

$$\int_0^{e(t_{p-1})} dy \int_y^\infty dz \int_0^\infty dt \int_0^\infty ds \, q_{e(t_{p-1})+z-2y}(t) q_z(s)$$

$$\times f\big(m(t_1,t_2),..,m(t_{p-2},t_{p-1}),y,e(t_1),..,e(t_{p-1}),z,t_1,..,t_{p-1}-t_{p-2},t,s\big)\Big).$$

It is then straightforward to complete the proof by induction on $p$. □

We can now state and prove the main result of this subsection.

**Theorem 2.11** *The law of the tree* $\theta(e;t_1,\ldots,t_p)$ *under the probability measure*

$$p! \, 1_{\{0 \leq t_1 \leq \cdots \leq t_p \leq 1\}} dt_1 \ldots dt_p \, n_{(1)}(de)$$

*is*

$$p! \, 2^{p+1} \, L(\theta) \exp\big(-2\,L(\theta)^2\big) \, \Lambda_p(d\theta).$$

**Proof.** Let $F$ be a nonnegative bounded continuous function on $\Theta_{(p)}$ and let $h$ be bounded, nonnegative and measurable on $\mathbb{R}_+$. By Proposition 2.10,

$$\int n(de) \, h(\sigma) \int_{\{0 \leq t_1 \leq \cdots \leq t_p \leq \sigma\}} dt_1 \ldots dt_p \, F\big(\theta(e;t_1,\ldots,t_p)\big)$$

$$= 2^{p-1} \int_0^\infty ds \, h(s) \int \Lambda_p(d\theta) \, q_{2L(\theta)}(s) \, F(\theta).$$

On the other hand, using the properties of the definition of the measures $n_{(s)}$, we have also

$$\int n(de) \, h(\sigma) \int_{\{0 \leq t_1 \leq \cdots \leq t_p \leq \sigma\}} dt_1 \ldots dt_p \, F\big(\theta(e;t_1,\ldots,t_p)\big)$$

$$= \frac{1}{2}(2\pi)^{-1/2} \int_0^\infty ds \, s^{-3/2} \, h(s) \int n_{(s)}(de)$$

$$\int_{\{0 \leq t_1 \leq \cdots \leq t_p \leq s\}} dt_1 \ldots dt_p \, F(\theta(e;t_1,\ldots,t_p)).$$

By comparing with the previous identity, we get for a.a. $s > 0$,

$$\int n_{(s)}(de) \int_{\{0 \leq t_1 \leq \cdots \leq t_p \leq s\}} dt_1 \ldots dt_p \, F(\theta(e; t_1, \ldots, t_p))$$

$$= 2^{p+1} \int \Lambda_p(d\theta) \, L(\theta) \, \exp\left(-\frac{2L(\theta)^2}{s}\right) F(\theta).$$

Both sides of the previous equality are continuous functions of $s$ (use the scaling property of $n_{(s)}$ for the left side). Thus the equality holds for every $s > 0$, and in particular for $s = 1$. This completes the proof. $\square$

**Concluding remarks.** If we pick $t_1, \ldots, t_p$ in $[0, 1]$ we can consider the increasing rearrangement $t'_1 \leq t'_2 \leq \cdots \leq t'_p$ of $t_1, \ldots, t_p$ and define $\theta(e; t_1, \ldots, t_p) = \theta(e; t'_1, \ldots, t'_p)$. We can also keep track of the initial ordering and consider the tree $\tilde{\theta}(e; t_1, \ldots, t_p)$ defined as the tree $\theta(e; t_1, \ldots, t_p)$ where leaves are labelled $1, \ldots, p$, the leaf corresponding to time $t_i$ receiving the label $i$. (This labelling has nothing to do with the ordering of the tree.) Theorem 2.11 implies that the law of the tree $\tilde{\theta}(e; t_1, \ldots, t_p)$ under the probability measure

$$1_{[0,1]^p}(t_1, \ldots, t_p) dt_1 \ldots dt_p \, n_{(1)}(de)$$

has density

$$2^{p+1} L(\theta) \exp(-2L(\theta)^2)$$

with respect to $\widetilde{\Theta}_{(p)}^{\mathrm{bin}}(d\theta)$, the uniform measure on the set of labelled marked (ordered) binary trees with $p$ leaves.

We can then "forget" the ordering. Let $\bar{\theta}(e; t_1, .., t_p)$ be the tree $\tilde{\theta}(e; t_1, .., t_p)$ without the order structure. Since there are $2^{p-1}$ possible orderings for a given labelled binary tree with $p$ leaves, we get that the law (under the same measure) of the tree $\bar{\theta}(e; t_1, \ldots, t_p)$ has density

$$2^{2p} L(\theta) \exp(-2L(\theta)^2)$$

with respect to $\bar{\Theta}_{(p)}^{\mathrm{bin}}(d\theta)$, the uniform measure on the set of labelled marked (unordered) binary trees with $p$ leaves.

In agreement with Aldous' normalization of the CRT, replace the excursion $e$ by $2e$ (this simply means that all marks $h_v$ are multiplied by 2). We obtain that the law of the tree $\bar{\theta}(2e; t_1, \ldots, t_p)$ has density

$$L(\theta) \exp\left(-\frac{L(\theta)^2}{2}\right)$$

with respect to $\bar{\Theta}_{(p)}^{\mathrm{bin}}(d\theta)$. It is remarkable that the previous density (apparently) does not depend on $p$.

In the previous form, we recognize the finite-dimensional marginals of Aldous' continuum random tree [1], [2]. To give a more explicit description, the discrete skeleton $\bar{\tau}(2e; t_1, \ldots, t_p)$ is distributed uniformly on the set of labelled rooted binary trees with $p$ leaves. This set has $b_p$ elements, with

$$b_p = p! \, 2^{-(p-1)} c_p = 1 \times 3 \times \cdots \times (2p-3).$$

Then, conditionally on the discrete skeleton, the marks $h_v$ are distributed with the density

$$b_p \left( \sum h_v \right) \exp \left( - \frac{\left( \sum h_v \right)^2}{2} \right)$$

(verify that this is a probability density on $\mathbb{R}_+^{2p-1}$ !).

*Bibliographical notes.* The coding of real trees described in subsection 2.2 is taken from [12], although the underlying ideas can be found in earlier papers (see in particular [2] and [25]). A simple approach to Theorem 2.2, based on the four-point condition for real trees, can be found in Lemma 3.1 of [18]. See e.g. Chapter XII of [38] or the last section of [39] for a thorough discussion of the Itô excursion measure. The CRT was introduced and studied by Aldous [1], [2]. The direct approach to finite-dimensional marginals of the CRT which is presented in subsections 2.5 and 2.6 above is taken from [26].

## 3. The Brownian Snake and its Connections with Partial Differential Equations

Our goal in this section is to combine the continuous tree structure studied in the previous section with independent spatial motions: In additional to the genealogical structure, "individuals" move in space according to the law of a certain Markov process $\xi$. This motivates the definition of the path-valued process called the Brownian snake. We study basic properties of the Brownian snake, and we use our previous calculation of marginal distributions of random real trees (subsection 2.5) to give explicit formulas for moment functionals. We then introduce the exit measure of the Brownian snake from a domain, and we derive a key integral equation for the Laplace functional of this random measure. In the case when the spatial motion $\xi$ is Brownian motion in $\mathbb{R}^d$, this integral equation leads to important connections with semilinear partial differential equations, which have been studied recently by several authors.

### 3.1. Combining the branching structure of a real tree with a spatial displacement

We consider a Markov process $(\xi_t, \Pi_x)_{t \geq 0, x \in E}$ with values in a Polish space $E$. We will assume that $\xi$ has continuous sample paths and that there exist an integer $m > 2$ and positive constants $C$ and $\varepsilon$ such that for every $x \in E$ and $t > 0$,

$$\Pi_x \left( \sup_{0 \leq r \leq t} \delta(x, \xi_r)^m \right) \leq C t^{2+\varepsilon}, \tag{27}$$

where $\delta$ denotes the distance on $E$. This assumption is not really necessary, but it will simplify our treatment. It holds for Brownian motion or for solutions of stochastic differential equations with smooth coefficients in $\mathbb{R}^d$ or on a manifold. Later $\xi$ will simply be $d$-dimensional Brownian motion, but for the moment it is preferable to argue in a more general setting.

We denote by $\mathcal{W}$ the set of all finite $E$-valued paths. An element of $\mathcal{W}$ is a continuous mapping $w : [0, \zeta] \to E$, where $\zeta = \zeta_{(w)} \geq 0$ depends on $w$ and is called the lifetime of $w$. The final point of $w$ will be denoted by $\hat{w} = w(\zeta)$. If $x \in E$, the trivial path $w$ such that $\zeta_{(w)} = 0$ and $w(0) = x$ is identified with the point $x$, so that $E$ is embedded in $\mathcal{W}$. The set $\mathcal{W}$ is a Polish space for the distance

$$d(w, w') = |\zeta - \zeta'| + \sup_{t \geq 0} \delta(w(t \wedge \zeta), w'(t \wedge \zeta')).$$

For $x \in E$ we denote by $\mathcal{W}_x$ the set $\mathcal{W}_x = \{w \in \mathcal{W} : w(0) = x\}$.

If $g : \mathbb{R}_+ \to \mathbb{R}_+$ is a continuous function with compact support such that $g(0) = 0$, we saw in subsection 2.2 that $g$ codes a real tree $\mathcal{T}_g$. Our goal is now to combine this branching structure with spatial motions distributed according to the law of the process $\xi$. To this end, we will not make an explicit use of the tree $\mathcal{T}_g$ but rather give our definitions in terms of the coding function. It will be convenient to drop the compact support assumption on $g$, and to consider instead a general function $f \in C(\mathbb{R}_+, \mathbb{R}_+)$.

**Notation.** Let $w : [0, \zeta] \to E$ be an element of $\mathcal{W}$, let $a \in [0, \zeta]$ and $b \geq a$. We denote by $R_{a,b}(w, dw')$ the unique probability measure on $\mathcal{W}$ such that

(i) $\zeta' = b$, $R_{a,b}(w, dw')$ a.s.
(ii) $w'(t) = w(t)$, $\forall t \leq a$, $R_{a,b}(w, dw')$ a.s.
(iii) The law under $R_{a,b}(w, dw')$ of $(w'(a+t), 0 \leq t \leq b-a)$ coincides with the law of $(\xi_t, 0 \leq t \leq b-a)$ under $\Pi_{w(a)}$.

Under $R_{a,b}(w, dw')$, the path $w'$ is the same as $w$ up to time $a$ and then behaves according to the spatial motion $\xi$ up to time $b$.

Let $(W_s, s \geq 0)$ denote the canonical process on the space $C(\mathbb{R}_+, \mathcal{W})$ of continuous functions from $\mathbb{R}_+$ into $\mathcal{W}$. We also denote by $\zeta_s = \zeta_{(W_s)}$ the lifetime of $W_s$.

**Proposition 3.1** *Let $f \in C(\mathbb{R}_+, \mathbb{R}_+)$. Assume that $f$ is locally Hölder continuous with exponent $\eta$ for every $\eta \in (0, \frac{1}{2})$. Let $w \in \mathcal{W}$ with $\zeta_{(w)} = f(0)$. Then, there exists a unique probability measure $\Gamma_w^f$ on $C(\mathbb{R}_+, \mathcal{W})$ such that $W_0 = w$, $\Gamma_w^f$ a.s., and, under $\Gamma_w^f$, the canonical process $(W_s, s \geq 0)$ is (time-inhomogeneous) Markov with transition kernel between times $s$ and $s'$ given by*

$$R_{m_f(s,s'),f(s')}(W_s, dw')$$

*where $m_f(s, s') = \inf_{[s,s']} f(r)$. We have in particular $\zeta_s = f(s)$ for every $s \geq 0$, $\Gamma_w^f$ a.s.*

The intuitive meaning of this construction should be clear : At least when $f(0) = 0$ and $f$ has compact support (so that we can use the theory of subsection 2.2), the vertices $p_f(s)$ and $p_f(s')$ in the tree $\mathcal{T}_f$ have the same ancestors up to generation $m_f(s, s')$. Therefore, the corresponding spatial motions $W_s$ and $W_{s'}$ must be the same up to time $m_f(s, s')$ and then behave independently. This means that the path $W_{s'}$ is obtained from the path $W_s$ through the kernel $R_{m_f(s,s'),f(s')}(W_s, dw')$.

**Proof.** We give the detailed argument only in the case when $f(0) = 0$, which implies that $w = x \in E$. We leave the general case as an exercise for the reader (cf the proof of Proposition IV.5 in [27]).

For each choice of $0 \leq t_1 \leq t_2 \leq \cdots \leq t_p$, we can consider the probability measure $\pi_{t_1,\ldots,t_p}^{x,f}$ on $\mathcal{W}^p$ defined by

$$
\begin{aligned}
&\pi_{t_1,\ldots,t_p}^{x,f}(dw_1 \ldots dw_p) \\
&= R_{0,f(t_1)}(x, dw_1) R_{m_f(t_1,t_2),f(t_2)}(w_1, dw_2) \ldots R_{m_f(t_{p-1},t_p),f(t_p)}(w_{p-1}, dw_p).
\end{aligned}
$$

It is easy to verify that this collection is consistent when $p$ and $t_1, \ldots, t_p$ vary. Hence the Kolmogorov extension theorem yields the existence of a process $(\widetilde{W}_s, s \geq 0)$ with values in $\mathcal{W}$ (in fact in $\mathcal{W}_x$) whose finite-dimensional marginals are the measures $\pi_{t_1,\ldots,t_p}^{x,f}$.

We then verify that $(\widetilde{W}_s, s \geq 0)$ has a continuous modification. Thanks to the classical Kolmogorov lemma, it is enough to show that, for every $T > 0$ there are constants $\beta > 0$ and $C$ such that

$$
E[d(\widetilde{W}_s, \widetilde{W}_{s'})^m] \leq C|s - s'|^{1+\beta}, \tag{28}
$$

for every $s \leq s' \leq T$. (Here $m$ is as in assumption (27).)

Our assumption on $f$ guarantees that for every $\eta \in (0, 1/2)$ there exists a finite constant $C_{\eta,T}$ such that, for every $s, s' \in [0, T]$,

$$
|f(s) - f(s')| \leq C_{\eta,T}|s - s'|^{\eta}. \tag{29}
$$

By our construction, the joint distribution of $(\widetilde{W}_s, \widetilde{W}_{s'})$ is

$$
R_{0,f(s)}(x, dw) R_{m_f(s,s'),f(s')}(w, dw').
$$

This means that $\widetilde{W}_s$ and $\widetilde{W}_{s'}$ are two random paths that coincide up to time $m_f(s, s')$ and then behave independently according to the law of the process $\xi$. Using the definition of the distance $d$, we get for every $s, s' \in [0, T]$, $s \leq s'$,

$$
\begin{aligned}
&E\big(d(\widetilde{W}_s, \widetilde{W}_{s'})^m\big) \\
&\leq c_m\Big(|f(s) - f(s')|^m + 2\,\Pi_x\Big(\Pi_{\xi_{m_f(s,s')}}\big(\sup_{0 \leq t \leq (f(s) \vee f(s')) - m_f(s,s')} \delta(\xi_0, \xi_t)^m\big)\Big)\Big) \\
&\leq c_m\Big(|f(s) - f(s')|^m + 2C|(f(s) \vee f(s')) - m_f(s,s')|^{2+\varepsilon}\Big) \\
&\leq c_m\Big(C_{\eta,T}^m|s - s'|^{m\eta} + 2C\,C_{\eta,T}^{2+\varepsilon}|s - s'|^{(2+\varepsilon)\eta}\Big),
\end{aligned}
$$

where we used assumption (27) in the second inequality, and (29) in the third one. We can choose $\eta$ close enough to $\frac{1}{2}$ so that $m\eta > 1$ and $(2+\varepsilon)\eta > 1$. The bound (28) then follows.

We then define $\Gamma_x^f$ as the law on $C(\mathbb{R}_+, \mathcal{W})$ of the continuous modification of $(\widetilde{W}_s, s \geq 0)$. The fact that under $\Gamma_x^f$ the canonical process is Markov with the given transition kernels is obvious from the choice of the finite-dimensional marginals. The form of the marginals also readily shows that $\zeta_{t_1} = f(t_1), \ldots, \zeta_{t_p} = f(t_p)$, $\Gamma_x^f$ a.s., for every finite collection of times $t_1, \ldots, t_p$. The last assertion of the proposition then follows from a continuity argument. $\qquad\square$

The process $(W_s, s \geq 0)$ under the probability measure $\Gamma_w^f$ is sometimes called the *snake* driven by the function $f$ (with spatial motion $\xi$ and initial value w). From the form of the transition kernels, and an easy continuity argument we have, $\Gamma_w^f$ a.s. for every $0 \leq s \leq s'$,

$$W_s(t) = W_{s'}(t) \qquad \text{for every } t \leq m_f(s, s').$$

We sometimes refer to this property as the *snake property*. Note in particular that if $x = w(0)$ we have $W_s(0) = x$ for every $s \geq 0$, $\Gamma_w^f$ a.s.

### 3.2. The Brownian snake

We now randomize $f$ in the construction of the previous subsection. For every $r \geq 0$, we denote by $P_r(df)$ the law of reflected Brownian motion started at $r$ (the law of $(|B_s|, s \geq 0)$ if $B$ is a linear Brownian motion started at $r$). Then $P_r(df)$ is a probability measure on the set $C(\mathbb{R}_+, \mathbb{R}_+)$. Note that the assumption of Proposition 3.1 holds $P_r(df)$ a.s.

For every $s > 0$, we denote by $\rho_s^r(da\,db)$ the law under $P_r$ of the pair

$$(\inf_{0 \leq u \leq s} f(u), f(s)).$$

The reflection principle (cf subsection 2.4) easily gives the explicit form of $\rho_s^r(da, db)$:

$$\rho_s^r(da, db) = \frac{2(r + b - 2a)}{(2\pi s^3)^{1/2}} \exp\Big(-\frac{(r + b - 2a)^2}{2s}\Big) 1_{(0 < a < b \wedge r)} \, da \, db$$

$$+ 2\,(2\pi s)^{-1/2} \exp\Big(-\frac{(r + b)^2}{2s}\Big) 1_{(0 < b)} \delta_0(da) \, db.$$

**Theorem 3.2** *For every* $w \in \mathcal{W}$, *denote by* $\mathbb{P}_w$ *be the probability measure on* $C(\mathbb{R}_+, \mathcal{W})$ *defined by*

$$\mathbb{P}_w(d\omega) = \int_{C(\mathbb{R}_+, \mathbb{R}_+)} P_{\zeta_{(w)}}(df) \, \Gamma_w^f(d\omega).$$

*The process* $(W_s, \mathbb{P}_w)_{s \geq 0, w \in \mathcal{W}}$ *is a continuous (time-homogeneous) Markov process with values in* $\mathcal{W}$, *with transition kernels*

$$\mathbb{Q}_s(w_1, dw_2) = \int \rho_s^{\zeta_{(w_1)}}(da, db) \, R_{a,b}(w_1, dw_2).$$

*This process is called the Brownian snake with spatial motion $\xi$.*

**Proof.** The continuity is obvious and it is clear that $W_0 = \mathrm{w}$, $\mathbb{P}_\mathrm{w}$ a.s. The semigroup property for the kernels $\mathbb{Q}_s$ is also easy to verify. As for the Markov property, we write

$$\mathbb{E}_\mathrm{w}[F(W_{s_1}, \ldots, W_{s_p})G(W_{s_{p+1}})]$$
$$= \int P_{\zeta_{(\mathrm{w})}}(df) \int_{\mathcal{W}^{p+1}} R_{m_f(0,s_1), f(s_1)}(x, d\mathrm{w}_1) \ldots$$
$$\ldots R_{m_f(s_p, s_{p+1}), f(s_{p+1})}(\mathrm{w}_p, d\mathrm{w}_{p+1})\, F(\mathrm{w}_1, \ldots, \mathrm{w}_p)G(\mathrm{w}_{p+1})$$
$$= \int_{\mathbb{R}_+^{2(p+1)}} \rho_{s_1}^{\zeta_{(\mathrm{w})}}(da_1, db_1)\rho_{s_2-s_1}^{b_1}(da_2, db_2) \ldots \rho_{s_{p+1}-s_p}^{b_p}(da_{p+1}, db_{p+1})$$
$$\int_{\mathcal{W}^{p+1}} R_{a_1, b_1}(\mathrm{w}, d\mathrm{w}_1) \ldots R_{a_{p+1}, b_{p+1}}(\mathrm{w}_p, d\mathrm{w}_{p+1})\, F(\mathrm{w}_1, \ldots, \mathrm{w}_p)G(\mathrm{w}_{p+1})$$
$$= \int_{\mathbb{R}_+^{2(p+1)}} \rho_{s_1}^{\zeta_{(\mathrm{w})}}(da_1, db_1)\rho_{s_2-s_1}^{b_1}(da_2, db_2) \ldots \rho_{s_p-s_{p-1}}^{b_{p-1}}(da_p, db_p)$$
$$\int_{\mathcal{W}^p} R_{a_1, b_1}(\mathrm{w}, d\mathrm{w}_1) \ldots R_{a_p, b_p}(\mathrm{w}_{p-1}, d\mathrm{w}_p)\, F(\mathrm{w}_1, \ldots, \mathrm{w}_p)\mathbb{Q}_{s_{p+1}-s_p}G(\mathrm{w}_p)$$
$$= \mathbb{E}_\mathrm{w}[F(W_{s_1}, \ldots, W_{s_p})\mathbb{Q}_{s_{p+1}-s_p}G(W_{s_p})].$$

This completes the proof. $\square$

Under $\mathbb{P}_\mathrm{w}$, the process $(\zeta_s, s \geq 0)$ is a reflected Brownian motion started at $\zeta_{(\mathrm{w})}$. This property is obvious from the last assertion of Proposition 3.1 and the very definition of $\mathbb{P}_\mathrm{w}$.

The snake property can then be stated in the form

$$W_s(t) = W_{s'}(t) \qquad \text{for all } t \leq \inf_{s \leq r \leq s'} \zeta_r,$$

for every $s < s'$, $\mathbb{P}_\mathrm{w}$ a.s.

In particular, if $\mathrm{w}(0) = x$ we have $W_s \in \mathcal{W}_x$ for every $s \geq 0$, $\mathbb{P}_\mathrm{w}$ a.s.

We now state the strong Markov property of $W$, which is very useful in applications. We denote by $\mathcal{F}_s$ the canonical filtration on $C(\mathbb{R}_+, \mathcal{W})$ ($\mathcal{F}_s$ is the $\sigma$-field generated by $W_r$, $0 \leq r \leq s$) and as usual we take

$$\mathcal{F}_{s+} = \bigcap_{r > s} \mathcal{F}_r\,.$$

**Theorem 3.3** *The process $(W_s, \mathbb{P}_\mathrm{w})_{s \geq 0, \mathrm{w} \in \mathcal{W}}$ is strong Markov with respect to the filtration $(\mathcal{F}_{s+})$.*

**Proof.** Let $T$ be a stopping time of the filtration $(\mathcal{F}_{s+})$ such that $T \leq K$ for some $K < \infty$. Let $F$ be bounded and $\mathcal{F}_{T+}$ measurable, and let $\Psi$ be a bounded measurable function on $\mathcal{W}$. It is enough to prove that for every $s > 0$,

$$\mathbb{E}_\mathrm{w}\big(F\,\Psi(W_{T+s})\big) = \mathbb{E}_\mathrm{w}\big(F\,\mathbb{E}_{W_T}(\Psi(W_s))\big).$$

We may assume that $\Psi$ is continuous. Then,

$$
\begin{aligned}
\mathbb{E}_{\mathrm{w}}\big(F\,\Psi(W_{T+s})\big) &= \lim_{n\to\infty}\sum_{k=0}^{\infty}\mathbb{E}_{\mathrm{w}}\big(\mathbf{1}_{\{\frac{k}{n}\le T<\frac{k+1}{n}\}}F\,\Psi(W_{\frac{k+1}{n}+s})\big)\\
&= \lim_{n\to\infty}\sum_{k=0}^{\infty}\mathbb{E}_{\mathrm{w}}\big(\mathbf{1}_{\{\frac{k}{n}\le T<\frac{k+1}{n}\}}F\,\mathbb{E}_{W_{\frac{k+1}{n}}}(\Psi(W_s))\big).
\end{aligned}
$$

In the first equality, we used the continuity of paths and in the second one the ordinary Markov property, together with the fact that $\mathbf{1}_{\{k/n\le T<(k+1)/n\}}F$ is $\mathcal{F}_{(k+1)/n}$-measurable. At this point, we need an extra argument. We claim that

$$
\lim_{\varepsilon\downarrow 0}\Big(\sup_{t\le K,t\le r\le t+\varepsilon}|\mathbb{E}_{W_r}\big(\Psi(W_s)\big)-\mathbb{E}_{W_t}\big(\Psi(W_s)\big)|\Big)=0,\qquad \mathbb{P}_{\mathrm{w}}\text{ a.s.}\qquad (30)
$$

Clearly, the desired result follows from (30), because on the set $\{k/n\le T<(k+1)/n\}$ we can bound

$$
|\mathbb{E}_{W_{\frac{k+1}{n}}}\big(\Psi(W_s)\big)-\mathbb{E}_{W_T}\big(\Psi(W_s)\big)|\le\sup_{t\le K,t\le r\le t+\frac{1}{n}}|\mathbb{E}_{W_r}\big(\Psi(W_s)\big)-\mathbb{E}_{W_t}\big(\Psi(W_s)\big)|.
$$

To prove (30), we write down explicitly

$$
\mathbb{E}_{W_r}\big(\Psi(W_s)\big)=\int\rho_s^{\zeta_r}(da,db)\int R_{a,b}(W_r,dw')\,\Psi(\mathrm{w}'),
$$

and a similar expression holds for $\mathbb{E}_{W_t}\big(\Psi(W_s)\big)$. Set

$$
c(\varepsilon)=\sup_{t\le K,t\le r\le t+\varepsilon}|\zeta_r-\zeta_t|
$$

and note that $c(\varepsilon)$ tends to 0 as $\varepsilon\to 0$, $\mathbb{P}_{\mathrm{w}}$ a.s. Then observe that if $t\le K$ and $t\le r\le r+\varepsilon$, the paths $W_r$ and $W_t$ coincide at least up to time $(\zeta_t-c(\varepsilon))_+$. Therefore we have

$$
R_{a,b}(W_r,dw')=R_{a,b}(W_t,dw')
$$

for every $a\le(\zeta_t-c(\varepsilon))_+$ and $b\ge a$. The claim (30) follows from this observation and the known explicit form of $\rho_s^{\zeta_r}(da,db)$. $\qquad\square$

**Remark.** The strong Markov property holds for $W$ even if the underlying spatial motion $\xi$ is not strong Markov.

### 3.3. Excursion measures of the Brownian snake

For every $x\in E$, the excursion measure $\mathbb{N}_x$ is the $\sigma$-finite measure on $C(\mathbb{R}_+,\mathcal{W})$ defined by

$$
\mathbb{N}_x(d\omega)=\int_{C(\mathbb{R}_+,\mathbb{R}_+)}n(de)\,\Gamma_x^e(d\omega),\qquad (31)
$$

where $n(de)$ denotes Itô's excursion measure as in Section 2. In particular, $(\zeta_s, s \geq 0)$ is distributed under $\mathbb{N}_x$ according to the Itô measure $n(de)$. As in Section 2, we will use the notation $\sigma$ for the duration of the excursion under $\mathbb{N}_x(d\omega)$:

$$\sigma = \inf\{s > 0 : \zeta_s = 0\}.$$

Note that $W_s \in \mathcal{W}_x$ for every $s \geq 0$, $\mathbb{N}_x$ a.e.

Let us comment on the intuitive interpretation of $\mathbb{N}_x$. As we saw in Section 2, the excursion $(\zeta_s, 0 \leq s \leq \sigma)$ codes a random real tree. If $a$ is a vertex of this real tree, corresponding say to $s \in [0, \sigma]$, the path $W_s$ gives the spatial positions along the line of ancestors of $a$ in the tree, and in particular its terminal point $\widehat{W}_s$ is the position of $a$. The snake property reflects the fact that the ancestors of the vertices corresponding to $s$ and $s'$ are the same up to generation $\inf_{[s,s']} \zeta_r$. To summarize, the paths

$$W_s , \ 0 \leq s \leq \sigma$$

form under $\mathbb{N}_x$ a "tree of paths" of the spatial motion $\xi$, whose genealogical structure is governed by the tree coded by the Brownian excursion $(\zeta_s, 0 \leq s \leq \sigma)$.

**Remark.** We know that the process $(W_s, \mathbb{P}_w)_{s \geq 0, w \in \mathcal{W}}$, where the driving random function is reflected Brownian motion, is a continuous strong Markov process. Furthermore, every point $x \in E$ (viewed as an element of $\mathcal{W}$) is regular for $W$, in the sense that $\mathbb{P}_x(T_{\{x\}} = 0) = 1$ if $T_{\{x\}} = \inf\{s > 0, W_s = x\}$. This last property is immediate from the analogous property for reflected linear Brownian motion. Thus it makes sense to consider the excursion measure of $W$ away from $x$, in the sense of the general Itô excursion theory (see e.g. Blumenthal [4]), and this excursion measure is easily identified with $\mathbb{N}_x$.

In what follows, since we use the notation $(W_s, s \geq 0)$ for the canonical process on $C(\mathbb{R}_+, \mathcal{W})$, it is important to realize that the driving random function of our Brownian snake $(W_s, s \geq 0)$ is either reflected Brownian motion (under the probability measures $\mathbb{P}_w$) or a Brownian excursion (under the excursion measures $\mathbb{N}_x$).

We will make use of the strong Markov property under $\mathbb{N}_x$. To state it, it will be convenient to introduce the following notation. For every $w \in \mathcal{W}$, we denote by $\mathbb{P}_w^*$ the distribution under $\mathbb{P}_w$ of the process

$$W_s^* = W_{s \wedge \sigma}, \quad s \geq 0$$

where $\sigma := \inf\{s > 0 : \zeta_s = 0\}$ as above.

**Theorem 3.4** *Let $T$ be a stopping time of the filtration $(\mathcal{F}_{s+})$. Assume that $0 < T \leq \sigma$, $\mathbb{N}_x$ a.e. Then, if $F$ and $G$ are nonnegative measurable functionals on $C(\mathbb{R}_+, \mathcal{W})$, and if $F$ is $\mathcal{F}_{T_+}$-measurable, we have*

$$\mathbb{N}_x\Big( F\, G(W_{T+s}, s \geq 0)\Big) = \mathbb{N}_x\Big( F\, \mathbb{E}_{W_T}^*(G)\Big).$$

If we interpret $\mathbb{N}_x$ as an excursion measure away from $x$ (as explained in the remark above), the preceding theorem becomes a well-known fact of the theory of Markov processes: See e.g. [4]. Alternatively, it is also easy to give a direct proof of Theorem 3.4 using the same method as in the proof of Theorem 3.3, together with the simple Markov property under the Itô excursion measure (cf Proposition 2.6 (ii)).

We will now use Theorem 2.7 to derive some important formulas under the excursion measure $\mathbb{N}_x$. Let $p \geq 1$ be an integer. Recall from subsection 2.5 the notation $\Theta_{(p)}^{\text{bin}}$ for the set of all marked trees with $p$ leaves, and let $\theta \in \Theta_{(p)}^{\text{bin}}$. For every $x \in E$, we associate with $\theta$ a probability measure on $\mathcal{W}^p$, denoted by $\Pi_x^\theta$, which is defined inductively as follows.

If $p = 1$, then $\theta = (\{\varnothing\}, h)$ for some $h \geq 0$ and we let $\Pi_x^\theta = \Pi_x^h$ be the law of $(\xi_t, 0 \leq t \leq h)$ under $\Pi_x$.

If $p \geq 2$, then we can write in a unique way

$$\theta = \theta' \underset{h}{*} \theta'' \ ,$$

where $\theta' \in \Theta_{(j)}^{\text{bin}}$, $\theta'' \in \Theta_{(p-j)}^{\text{bin}}$, and $j \in \{1, \ldots, p-1\}$ (the notation $\theta' \underset{h}{*} \theta''$, for the concatenation of $\theta'$ and $\theta''$ with root mark $h$, was introduced in the proof of Theorem 1.15) . We then define $\Pi_x^\theta$ by

$$\int \Pi_x^\theta(d\mathrm{w}_1, \ldots, d\mathrm{w}_p) F(\mathrm{w}_1, \ldots, \mathrm{w}_p)$$

$$= \Pi_x \Big( \int\!\!\int \Pi_{\xi_h}^{\theta'}(d\mathrm{w}_1', \ldots, d\mathrm{w}_j') \Pi_{\xi_h}^{\theta''}(d\mathrm{w}_1'', \ldots, d\mathrm{w}_{p-j}'')$$

$$F(\xi_{[0,h]} \odot \mathrm{w}_1', \ldots, \xi_{[0,h]} \odot \mathrm{w}_j', \xi_{[0,h]} \odot \mathrm{w}_1'', \ldots, \xi_{[0,h]} \odot \mathrm{w}_{p-j}'') \Big)$$

where $\xi_{[0,h]} \odot \mathrm{w}$ denotes the concatenation (defined in an obvious way) of the paths $(\xi_t, 0 \leq t \leq h)$ and $(\mathrm{w}(t), 0 \leq t \leq \zeta_{(\mathrm{w})})$.

Informally, $\Pi_x^\theta$ is obtained by running independent copies of $\xi$ along the branches of the tree $\theta$.

Finally, we recall the notation $\theta(f; t_1, \ldots, t_p)$ from subsection 2.5, and we let $\Lambda_p$ be as in Theorem 2.7 the uniform measure on $\Theta_{(p)}^{\text{bin}}$.

**Proposition 3.5** (i) *Let $f \in C(\mathbb{R}_+, \mathbb{R}_+)$ such that $f(0) = 0$, and let $0 \leq t_1 \leq t_2 \leq \cdots \leq t_p$. Then the law under $\Gamma_x^f$ of $(W_{t_1}, \ldots, W_{t_p})$ is $\Pi_x^{\theta(f; t_1, \ldots, t_p)}$.*

(ii) *For any nonnegative Borel measurable function $F$ on $\mathcal{W}^p$,*

$$\mathbb{N}_x \Big( \int_{\{0 \leq s_1 \leq \cdots \leq s_p \leq \sigma\}} ds_1 \ldots ds_p \, F(W_{s_1}, \ldots, W_{s_p}) \Big) = 2^{p-1} \int \Lambda_p(d\theta) \, \Pi_x^\theta(F) \ .$$

**Proof.** Assertion (i) follows easily from the definition of $\Gamma_x^f$ and the construction of the trees $\theta(f; t_1, \ldots, t_p)$. A precise argument can be given using induction on $p$, but we leave details to the reader. To get (ii), we write

$$\mathbb{N}_x \Big( \int_{\{0 \le t_1 \le \cdots \le t_p \le \sigma\}} dt_1 \ldots dt_p \, F\big(W_{t_1}, \ldots, W_{t_p}\big) \Big)$$

$$= \int n(de) \int_{\{0 \le t_1 \le \cdots \le t_p \le \sigma\}} dt_1 \ldots dt_p \, \Gamma_x^e \Big( F\big(W_{t_1}, \ldots, W_{t_p}\big) \Big)$$

$$= \int n(de) \int_{\{0 \le t_1 \le \cdots \le t_p \le \sigma\}} dt_1 \ldots dt_p \, \Pi_x^{\theta(e; t_1, \ldots, t_p)}(F)$$

$$= 2^{p-1} \int \Lambda_p(d\theta) \, \Pi_x^\theta(F).$$

The first equality is the definition of $\mathbb{N}_x$, the second one is part (i) of the proposition, and the last one is Theorem 2.7. $\qquad\square$

The cases $p = 1$ and $p = 2$ of Proposition 3.5 (ii) are of special interest. Let us rewrite the corresponding formulas in a particular case. Recall that we denote by $\hat{w}$ the terminal point of w. For any nonnegative Borel measurable function $g$ on $E$, we have

$$\mathbb{N}_x \Big( \int_0^\sigma ds \, g(\hat{W}_s) \Big) = \Pi_x \Big( \int_0^\infty dt \, g(\xi_t) \Big),$$

and

$$\mathbb{N}_x \Big( \Big( \int_0^\sigma ds \, g(\hat{W}_s) \Big)^2 \Big) = 4 \, \Pi_x \Big( \int_0^\infty dt \, \Big( \Pi_{\xi_t} \Big( \int_0^\infty dr \, g(\xi_r) \Big) \Big)^2 \Big).$$

**Remark**. In addition to (31), we could also consider the associated normalized excursion measure

$$\mathbb{N}_x^{(1)}(d\omega) = \int_{C(\mathbb{R}_+, \mathbb{R}_+)} n_{(1)}(de) \, \Gamma_x^e(d\omega), \tag{32}$$

where $n_{(1)}$ is as in Section 2 the law of the normalized Brownian excursion. Let $\mathcal{Z}$ be the random probability measure on $E$ defined under $\mathbb{N}_x^{(1)}$ by

$$\langle \mathcal{Z}, g \rangle = \int_0^1 ds \, g(\hat{W}_s).$$

In the case when $\xi$ is Brownian motion in $E = \mathbb{R}^d$ and $x = 0$, the random measure $\mathcal{Z}$ is called ISE (Aldous [3]) for Integrated Super-Brownian Excursion. ISE and its variants play an important role in various asymptotics for statistical mechanics models (see e.g. [7], [21]). In such applications, the explicit formula for the moments of the random measure $\mathcal{Z}$ is often useful. This formula is proved by the same method as Proposition 3.5 (ii), using Theorem 2.11 instead of Theorem 2.7. Precisely, we have the following result.

**Proposition 3.6** *For any nonnegative Borel measurable function $g$ on $E$,*

$$\mathbb{N}_x^{(1)}(\langle \mathcal{Z}, g \rangle^p)$$
$$= p! 2^{p-1} \int \Lambda_p(d\theta)\, L(\theta) \exp(-2L(\theta)^2) \int \Pi_x^\theta(dw_1 \ldots dw_p)\, g(\hat{w}_1) \cdots g(\hat{w}_p).$$

We conclude this subsection with an important technical lemma. We fix $w \in \mathcal{W}_x$ with $\zeta_{(w)} > 0$ and now consider the Brownian snake under $\mathbb{P}_w$, that is the driving random function of the snake is reflected Brownian motion (and no longer a Brownian excursion as in the beginning of this section). We again use the notation

$$\sigma = \inf\{s > 0 : \zeta_s = 0\}$$

and denote by $(\alpha_i, \beta_i)$, $i \in I$ the excursion intervals of $\zeta_s - \inf_{[0,s]} \zeta_r$ before time $\sigma$. In other words, $(\alpha_i, \beta_i)$, $i \in I$ are the connected components of the open set $[0, \sigma] \cap \{s \geq 0, \zeta_s > \inf_{[0,s]} \zeta_r\}$. Then, for every $i \in I$ we define $W^i \in C(\mathbb{R}_+, \mathcal{W})$ by setting for every $s \geq 0$,

$$W_s^i(t) = W_{(\alpha_i + s) \wedge \beta_i}(\zeta_{\alpha_i} + t), \qquad 0 \leq t \leq \zeta_s^i := \zeta_{(\alpha_i + s) \wedge \beta_i} - \zeta_{\alpha_i}.$$

From the snake property we have in fact $W^i \in C(\mathbb{R}_+, \mathcal{W}_{w(\zeta_{\alpha_i})})$.

**Lemma 3.7** *The point measure*

$$\sum_{i \in I} \delta_{(\zeta_{\alpha_i}, W^i)}$$

*is under $\mathbb{P}_w$ a Poisson point measure on $\mathbb{R}_+ \times C(\mathbb{R}_+, \mathcal{W})$ with intensity*

$$2\, 1_{[0, \zeta_{(w)}]}(t) dt\, \mathbb{N}_{w(t)}(d\omega)\,.$$

**Proof.** A well-known theorem of Lévy (already used in Section 1) states that, if $(\beta_t, t \geq 0)$ is a linear Brownian motion started at $a$, the process $\beta_t - \inf_{[0,t]} \beta_r$ is a reflected Brownian motion whose local time at $0$ is $t \to 2\big(a - \inf_{[0,t]} \beta_r\big)$. From this and excursion theory, it follows that the point measure

$$\sum_{i \in I} \delta_{(\zeta_{\alpha_i}, \zeta^i)}$$

is under $\mathbb{P}_w$ a Poisson point measure with intensity

$$2\, 1_{[0, \zeta_w]}(t) dt\, n(de)\,.$$

It remains to combine this result with the spatial displacements.

To this end, fix a function $f \in C(\mathbb{R}_+, \mathbb{R}_+)$ such that $f(0) = \zeta_{(w)}$, $\sigma(f) = \inf\{t > 0 : f(t) = 0\} < \infty$ and $f$ is locally Hölder with exponent $\frac{1}{2} - \gamma$ for every $\gamma > 0$. Recall the notation $\Gamma_w^f$ from subsection 3.1 above. Denote by $e_j, j \in J$ the excursions of $f(s) - \inf_{[0,s]} f(r)$ away from $0$ before time $\sigma(f)$, by $(a_j, b_j)$, $j \in J$ the corresponding time intervals, and define for every $j \in J$

$$W_s^j(t) = W_{(a_j + s) \wedge b_j}\big(f(a_j) + t\big)\,, \quad 0 \leq t \leq f\big((a_j + s) \wedge b_j\big) - f(a_j)\,,$$

From the definition of $\Gamma_{\mathrm{w}}^f$, it is easily verified that the processes $W^j$, $j \in J$ are independent under $\Gamma_{\mathrm{w}}^f$, with respective distributions $\Gamma_{\mathrm{w}(f(a_j))}^{e_j}$.

Let $F$ be a bounded nonnegative measurable function on $\mathbb{R}_+ \times C(\mathbb{R}_+, \mathcal{W})$, such that $F(t, \omega) = 0$ if $\sup \zeta_s(\omega) \leq \gamma$, for some $\gamma > 0$. Recall the notation $P_r(df)$ for the law of reflected Brownian motion started at $r$. By using the last observation and then the beginning of the proof, we get

$$
\begin{aligned}
\mathbb{E}_{\mathrm{w}}\Big(\exp -\sum_{i \in I} F(\zeta_{\alpha_i}, W^i)\Big) &= \int P_{\zeta_{(\mathrm{w})}}(df)\Gamma_{\mathrm{w}}^f\Big(\exp -\sum_{j \in J} F\big(f(a_j), W^j\big)\Big) \\
&= \int P_{\zeta_{(\mathrm{w})}}(df) \prod_{j \in J} \Gamma_{\mathrm{w}(f(a_j))}^{e_j}\Big(e^{-F(f(a_j), \cdot)}\Big) \\
&= \exp -2 \int_0^{\zeta_{(\mathrm{w})}} dt \int n(de)\Gamma_{\mathrm{w}(t)}^e\big(1 - e^{-F(t, \cdot)}\big) \\
&= \exp -2 \int_0^{\zeta_{(\mathrm{w})}} dt\, \mathbb{N}_{\mathrm{w}(t)}\big(1 - e^{-F(t, \cdot)}\big)\ .
\end{aligned}
$$

The third equality is the exponential formula for Poisson measures, and the last one is the definition of $\mathbb{N}_x$. This completes the proof. $\qquad\square$

### *3.4. The exit measure*

Let $D$ be an open set in $E$ and fix $x \in D$. For every $\mathrm{w} \in \mathcal{W}_x$ set

$$\tau(\mathrm{w}) = \inf\big\{t \in [0, \zeta_{(\mathrm{w})}], \mathrm{w}(t) \notin D\big\}\ ,$$

where $\inf \varnothing = +\infty$. Define

$$\mathcal{E}^D = \big\{W_s\big(\tau(W_s)\big); s \geq 0, \tau(W_s) < \infty\big\}\ ,$$

so that $\mathcal{E}^D$ is the set of all exit points from $D$ of the paths $W_s$, for those paths that do exit $D$. Our goal is to construct $\mathbb{N}_x$ a.e. a random measure that is in some sense uniformly spread over $\mathcal{E}^D$. To avoid trivial cases, we first assume that

$$\Pi_x(\exists t \geq 0, \xi_t \notin D) > 0\ . \tag{33}$$

We start by constructing a continuous increasing process that increases only on the set $\{s \geq 0 : \tau(W_s) = \zeta_s\}$.

**Proposition 3.8** *The formula*

$$L_s^D = \lim_{\varepsilon \downarrow 0} \frac{1}{\varepsilon} \int_0^s dr\, 1_{\{\tau(W_r) < \zeta_r < \tau(W_r) + \varepsilon\}}$$

*defines a continuous increasing process* $(L_s^D, s \geq 0)$, $\mathbb{N}_x$ *a.e. or* $\mathbb{P}_{\mathrm{w}}$ *a.s. for any* $\mathrm{w} \in \mathcal{W}_x$.

**Proof.** Since $\mathbb{N}_x$ can be viewed as the excursion measure of $W$ away from $x$, it is enough to prove that the given statement holds under $\mathbb{P}_x$. Indeed, if follows from the construction of the Itô measure that, for every $h > 0$, $\mathbb{N}_x(\cdot | \sup \zeta_s > h)$ is the law under $\mathbb{P}_x$ of the first excursion of $W$ away from $x$ with "height" greater than $h$, and so the result under $\mathbb{N}_x$ can easily be derived from the case of $\mathbb{P}_x$.

We use the following lemma, where $\mathrm{w} \in \mathcal{W}_x$ is fixed.

**Lemma 3.9** *Set* $\gamma_s = \left( \zeta_s - \tau(W_s) \right)^+$ *and* $\sigma_s = \inf\{v \geq 0 : \int_0^v dr \, 1_{\{\gamma_r > 0\}} > s\}$. *Then* $\sigma_s < \infty$ *for every* $s \geq 0$, $\mathbb{P}_\mathrm{w}$ *a.s., and the process* $\Gamma_s = \gamma_{\sigma_s}$ *is under* $\mathbb{P}_\mathrm{w}$ *a reflected Brownian motion started at* $(\zeta_\mathrm{w} - \tau(\mathrm{w}))^+$.

Proposition 3.8 easily follows from Lemma 3.9: Denote by $(\ell_s, s \geq 0)$ the local time at $0$ of $\Gamma$. Then, $\mathbb{P}_x$ a.s. for every $s \geq 0$,

$$\ell_s = \lim_{\varepsilon \to 0} \frac{1}{\varepsilon} \int_0^s dr \, 1_{\{0 < \Gamma_r < \varepsilon\}}.$$

Set $A_s = \int_0^s dr \, 1_{\{\gamma_r > 0\}}$ and $L_s^D = \ell_{A_s}$. We get

$$L_s^D = \lim_{\varepsilon \downarrow 0} \frac{1}{\varepsilon} \int_0^{A_s} dr \, 1_{\{0 < \Gamma_r < \varepsilon\}} = \lim_{\varepsilon \downarrow 0} \frac{1}{\varepsilon} \int_0^s dr \, 1_{\{0 < \gamma_r < \varepsilon\}},$$

which is the desired result. $\qquad\qquad\qquad\qquad\qquad\qquad\qquad\qquad\square$

**Proof of Lemma 3.9.** For every $\varepsilon > 0$, introduce the stopping times

$$S_1^\varepsilon = \inf\{s \geq 0 : \zeta_s \geq \tau(W_s) + \varepsilon\} \qquad T_1^\varepsilon = \inf\{s \geq S_1^\varepsilon : \zeta_s \leq \tau(W_s)\}$$

$$S_{n+1}^\varepsilon = \inf\{s \geq T_n^\varepsilon : \zeta_s \geq \tau(W_s) + \varepsilon\} \quad T_{n+1}^\varepsilon = \inf\{s \geq S_{n+1}^\varepsilon : \zeta_s \leq \tau(W_s)\}$$

We first verify that the stopping times $S_n^\varepsilon$ and $T_n^\varepsilon$ are finite $\mathbb{P}_\mathrm{w}$ a.s. By applying the strong Markov property at $\inf\{s \geq 0, \zeta_s = 0\}$, it is enough to consider the case when $\mathrm{w} = x$. Still another application of the strong Markov property shows that it is enough to verify that $S_1^\varepsilon < \infty$ a.s. To this end, observe that $\mathbb{P}_x\left(\zeta_1 \geq \tau(W_1) + \varepsilon\right) > 0$ (by (33) and because, conditionally on $\zeta_1$, $W_1$ is a path of $\xi$ with length $\zeta_1$) and apply the strong Markov property at $\inf\{s \geq 1, \zeta_s = 0\}$.

From the snake property and the continuity of $s \to \zeta_s$, one easily gets that the mapping $s \to \gamma_s$ is also continuous. It follows that $\gamma_{S_1^\varepsilon} = \varepsilon \vee (\zeta_{(\mathrm{w})} - \tau(\mathrm{w}))$ and $\gamma_{S_n^\varepsilon} = \varepsilon$ for $n \geq 2$.

We then claim that, for every $n \geq 1$, we have

$$T_n^\varepsilon = \inf\{s \geq S_n^\varepsilon : \zeta_s = \tau(W_{S_n^\varepsilon})\}.$$

Indeed the snake property implies that for

$$S_n^\varepsilon \leq r \leq \inf\{s \geq S_n^\varepsilon : \zeta_s = \tau(W_{S_n^\varepsilon})\},$$

the paths $W_r$ and $W_{S_n^\varepsilon}$ coincide for $t \le \tau(W_{S_n^\varepsilon})$, so that $\tau(W_r) = \tau(W_{S_n^\varepsilon})$. This argument also shows that $\gamma_r = \zeta_r - \tau(W_{S_n^\varepsilon})$ for $S_n^\varepsilon \le r \le T_n^\varepsilon$.

From the previous observations and the strong Markov property of the Brownian snake, we see that the processes

$$\left(\gamma_{(S_n^\varepsilon + r) \wedge T_n^\varepsilon}, r \ge 0\right), \quad n = 1, 2, \dots$$

are independent and distributed according to the law of a linear Brownian motion started at $\varepsilon$ (at $\varepsilon \vee (\zeta_{(\mathrm{w})} - \tau(\mathrm{w}))$ for $n = 1$) and stopped when it hits 0. Hence, if

$$\sigma_r^\varepsilon = \inf\left\{ s : \int_0^s \sum_{n=1}^\infty 1_{[S_n^\varepsilon, T_n^\varepsilon]}(u)du > r \right\},$$

the process $(\gamma_{\sigma_r^\varepsilon}, r \ge 0)$ is obtained by pasting together a linear Brownian motion started at $\varepsilon \vee (\zeta_{(\mathrm{w})} - \tau(\mathrm{w}))$ and stopped when it hits 0, with a sequence of independent copies of the same process started at $\varepsilon$. A simple coupling argument shows that $(\gamma_{\sigma_r^\varepsilon}, r \ge 0)$ converges in distribution as $\varepsilon \to 0$ to reflected Brownian motion started at $(\zeta_{(\mathrm{w})} - \tau(\mathrm{w}))^+$. The lemma follows since it is clear that $\sigma_r^\varepsilon \downarrow \sigma_r$ a.s. for every $r \ge 0$. $\qquad \square$

**Definition.** *The exit measure $\mathcal{Z}^D$ from $D$ is defined under $\mathbb{N}_x$ by the formula*

$$\langle \mathcal{Z}^D, g \rangle = \int_0^\sigma dL_s^D g(\hat{W}_s) \ .$$

From Proposition 3.8 it is easy to obtain that $L_s^D$ increases only on the (closed) set $\{s \in [0, \sigma] : \zeta_s = \tau(W_s)\}$. It follows that $\mathcal{Z}^D$ is ($\mathbb{N}_x$ a.e.) supported on $\mathcal{E}^D$.

Let us consider the case when (33) does not hold. Then a first moment calculation using the case $p = 1$ of Proposition 3.5 (ii) shows that

$$\int_0^\infty ds \, 1_{\{\tau(W_s) < \infty\}} = 0 \ , \ \mathbb{N}_x \ \text{ a.e.}$$

Therefore the result of Proposition 3.8 still holds under $\mathbb{N}_x$ with $L_s^D = 0$ for every $s \ge 0$. Consequently, we take $\mathcal{Z}^D = 0$ in that case.

We will need a first moment formula for $L^D$. With a slight abuse of notation, we also denote by $\tau$ the first exit time from $D$ for $\xi$.

**Proposition 3.10** *Let $\Pi_x^D$ denote the law of $(\xi_r, 0 \le r \le \tau)$ under the subprobability measure $\Pi_x(\cdot \cap \{\tau < \infty\})$ ($\Pi_x^D$ is viewed as a measure on $\mathcal{W}_x$). Then, for every bounded nonnegative measurable function $G$ on $\mathcal{W}_x$,*

$$\mathbb{N}_x\left(\int_0^\sigma dL_s^D G(W_s)\right) = \Pi_x^D(G) \ .$$

*In particular, for every bounded nonnegative measurable function $g$ on $E$,*

$$\mathbb{N}_x\big(\langle \mathcal{Z}^D, g\rangle\big) = \Pi_x\big(1_{\{\tau < \infty\}}g(\xi_\tau)\big).$$

**Proof.** We may assume that $G$ is continuous and bounded, and $G(w) = 0$ if $\zeta_{(w)} \leq K^{-1}$ or $\zeta_{(w)} \geq K$, for some $K > 0$. By Proposition 3.8,

$$\int_0^\sigma dL_s^D\, G(W_s) = \lim_{\varepsilon \to 0} \frac{1}{\varepsilon}\int_0^\sigma ds\, 1_{\{\tau(W_s) < \zeta_s < \tau(W_s)+\varepsilon\}}\, G(W_s) \qquad (34)$$

$\mathbb{N}_x$ a.e. If we can justify the fact that the convergence (34) also holds in $L^1(\mathbb{N}_x)$, we will get from the case $p = 1$ of Proposition 3.5 (ii):

$$
\begin{aligned}
\mathbb{N}_x\Big(\int_0^\sigma dL_s^D G(W_s)\Big) &= \lim_{\varepsilon \to 0} \frac{1}{\varepsilon}\int_0^\infty dh\, \Pi_x\Big(1_{\{\tau < h < \tau+\varepsilon\}}G(\xi_r, 0 \leq r \leq h)\Big) \\
&= \lim_{\varepsilon \to 0} \Pi_x\Big(1_{\{\tau < \infty\}}\, \varepsilon^{-1}\int_\tau^{\tau+\varepsilon} dh\, G(\xi_r, 0 \leq r \leq h)\Big) \\
&= \Pi_x\Big(1_{\{\tau < \infty\}}G(\xi_r, 0 \leq r \leq \tau)\Big).
\end{aligned}
$$

It remains to justify the convergence in $L^1(\mathbb{N}_x)$. Because of our assumption on $G$ we may deal with the finite measure $\mathbb{N}_x\big(\cdot \cap \{\sup \zeta_s > K^{-1}\}\big)$ and so it is enough to prove that

$$\sup_{0 < \varepsilon < 1}\, \mathbb{N}_x\bigg(\Big(\frac{1}{\varepsilon}\int_0^\sigma ds\, 1_{\{\tau(W_s) < \zeta_s < \tau(W_s)+\varepsilon\}}G(W_s)\Big)^2\bigg)$$

is finite. This easily follows from the case $p = 2$ of Proposition 3.5 (ii), using now the fact that $G(w) = 0$ if $\zeta_{(w)} \geq K$. $\qquad \square$

Let us give an important remark. Without any additional effort, the previous construction applies to the more general case of a space-time open set $D \subset \mathbb{R}_+ \times E$, such that $(0, x) \in D$. In this setting, $\mathcal{Z}^D$ is a random measure on $\partial D \subset \mathbb{R}_+ \times E$ such that for $g \in C_{b+}(\partial D)$

$$\langle \mathcal{Z}^D, g\rangle = \lim_{\varepsilon \to 0} \frac{1}{\varepsilon}\int_0^\sigma ds\, 1_{\{\tau(W_s) < \zeta_s < \tau(W_s)+\varepsilon\}}g(\zeta_s, \hat{W}_s)$$

where $\tau(w) = \inf\big\{t \geq 0,\, (t, w(t)) \notin D\big\}$. To see that this more general case is in fact contained in the previous construction, simply replace $\xi$ by the space-time process $\xi'_t = (t, \xi_t)$, which also satisfies assumption (27), and note that the Brownian snake with spatial motion $\xi'$ is related to the Brownian snake with spatial motion $\xi$ in a trivial manner.

We will now derive an integral equation for the Laplace functional of the exit measure. This result is the key to the connections with partial differential equations that will be investigated later.

**Theorem 3.11** *Let $g$ be a nonnegative bounded measurable function on $E$. For every $x \in E$, set*

$$u(x) = \mathbb{N}_x\big(1 - \exp-\langle \mathcal{Z}^D, g \rangle\big) , \quad x \in D .$$

*The function $u$ solves the integral equation*

$$u(x) + 2\Pi_x\Big(\int_0^\tau u(\xi_s)^2 ds\Big) = \Pi_x\big(1_{\{\tau < \infty\}} g(\xi_\tau)\big) . \tag{35}$$

Our proof of Theorem 3.11 is based on Lemma 3.7. Another more computational proof would rely on calculations of moments of the exit measure from Proposition 3.5 above.

**Proof.** For every $r > 0$ set $\eta_r^D = \inf\{s \geq 0 : L_s^D > r\}$, with the usual convention $\inf \varnothing = \infty$. By the definition of $\mathcal{Z}^D$, we have

$$
\begin{aligned}
u(x) &= \mathbb{N}_x\Big(1 - \exp-\int_0^\sigma dL_s^D g(\hat{W}_s)\Big) \\
&= \mathbb{N}_x\Big(\int_0^\sigma dL_s^D g(\hat{W}_s) \exp\big(-\int_s^\sigma dL_r^D g(\hat{W}_r)\big)\Big) \\
&= \mathbb{N}_x\Big(\int_0^\infty dr \, 1_{\{\eta_r^D < \infty\}} g(\hat{W}_{\eta_r^D}) \, \exp-\int_{\eta_r^D}^\sigma dL_s^D g(\hat{W}_s)\big)\Big) \\
&= \mathbb{N}_x\Big(\int_0^\infty dr \, 1_{\{\eta_r^D < \infty\}} g(\hat{W}_{\eta_r^D}) \, \mathbb{E}_{W_{\eta_r^D}}\big(\exp-\int_0^\sigma dL_s^D g(\hat{W}_s)\big)\Big) \\
&= \mathbb{N}_x\Big(\int_0^\sigma dL_s^D g(\hat{W}_s) \mathbb{E}_{W_s}\big(\exp-\int_0^\sigma dL_r^D g(\hat{W}_r)\big)\Big).
\end{aligned}
$$

The second equality is the simple identity $1 - \exp(-A_t) = \int_0^t dA_s \, \exp(-(A_t - A_s))$ valid for any continuous nondecreasing function $A$. The third equality is the change of variables $s = \eta_r^D$ and the fourth one follows from the strong Markov property under $\mathbb{N}_x$ (cf Theorem 3.4) at the stopping time $\eta_r^D$.

Let $w \in \mathcal{W}_x$ be such that $\zeta_{(w)} = \tau_{(w)}$. From Lemma 3.7, we have

$$
\begin{aligned}
\mathbb{E}_w\Big(\exp-\int_0^\sigma dL_r^D g(\hat{W}_r)\Big) &= \mathbb{E}_w\Big(\exp-\sum_{i \in I}\int_{\alpha_i}^{\beta_i} dL_r^D g(\hat{W}_r)\Big) \\
&= \exp\Big(-2\int_0^{\zeta_{(w)}} dt \, \mathbb{N}_{w(t)}\big(1 - \exp-\int_0^\sigma dL_r^D g(\hat{W}_r)\big)\Big) \\
&= \exp\Big(-2\int_0^{\zeta_{(w)}} dt \, u(w(t))\Big) .
\end{aligned}
$$

Hence,

$$
\begin{aligned}
u(x) &= \mathbb{N}_x\Big(\int_0^\sigma dL_s^D g(\hat{W}_s) \exp\big(-2\int_0^{\zeta_s} dt \, u(W_s(t))\big)\Big) \\
&= \Pi_x\Big(1_{\{\tau < \infty\}} g(\xi_\tau) \exp\big(-2\int_0^\tau dt \, u(\xi_t)\big)\Big)
\end{aligned}
$$

by Proposition 3.10. The proof is now easily completed by the usual Feynman-Kac argument:

$$
\begin{aligned}
u(x) &= \Pi_x\big(1_{\{\tau < \infty\}}g(\xi_\tau)\big) - \Pi_x\Big(1_{\{\tau < \infty\}}g(\xi_\tau)\big(1 - \exp{-2\int_0^\tau dt\ u(\xi_t)}\big)\Big) \\
&= \Pi_x\big(1_{\{\tau < \infty\}}g(\xi_\tau)\big) - 2\Pi_x\Big(1_{\{\tau < \infty\}}g(\xi_\tau)\int_0^\tau dt\ u(\xi_t)\exp\big(-2\int_t^\tau dr\ u(\xi_r)\big)\Big) \\
&= \Pi_x\big(1_{\{\tau < \infty\}}g(\xi_\tau)\big) - 2\Pi_x\Big(\int_0^\tau dt u(\xi_t)\Pi_{\xi_t}\big(1_{\{\tau < \infty\}}g(\xi_\tau)\exp\big(-2\int_0^\tau dru(\xi_r)\big)\big)\Big) \\
&= \Pi_x\big(1_{\{\tau < \infty\}}g(\xi_\tau)\big) - 2\Pi_x\Big(\int_0^\tau dt\ u(\xi_t)^2\Big)\ .
\end{aligned}
$$

### 3.5. The probabilistic solution of the nonlinear Dirichlet problem

In this subsection, we assume that $\xi$ is Brownian motion in $\mathbb{R}^d$. The results however could easily be extended to an elliptic diffusion process in $\mathbb{R}^d$ or on a manifold.

We say that $y \in \partial D$ is regular for $D^c$ if

$$
\inf\{t > 0 : \xi_t \notin D\} = 0 \qquad , \quad \Pi_y \text{ a.s.}
$$

The open set $D$ is called regular if every point $y \in \partial D$ is regular for $D^c$. We say that a real-valued function $u$ defined on $D$ solves $\Delta u = 4u^2$ in $D$ if $u$ is of class $C^2$ on $D$ and the equality $\Delta u = 4u^2$ holds pointwise on $D$.

**Theorem 3.12** *Let $D$ be a domain in $\mathbb{R}^d$ and let $g$ be a bounded nonnegative measurable function on $\partial D$. For every $x \in D$, set $u(x) = \mathbb{N}_x(1 - \exp{-\langle \mathcal{Z}^D, g\rangle})$. Then $u$ solves $\Delta u = 4u^2$ in $D$. If in addition $D$ is regular and $g$ is continuous, then $u$ solves the problem*

$$
\begin{aligned}
\Delta u &= 4u^2 &&in\ D \\
u_{|\partial D} &= g
\end{aligned}
\tag{36}
$$

*where the notation $u_{|\partial D} = g$ means that for every $y \in \partial D$,*

$$
\lim_{D \ni x \to y} u(x) = g(y)\ .
$$

**Proof.** First observe that, by (35),

$$
u(x) \le \Pi_x\big(1_{\{\tau < \infty\}}g(\xi_\tau)\big) \le \sup_{y \in \partial D} g(y)\ ,
$$

so that $u$ is bounded in $D$. Let $B$ be a ball whose closure is contained in $D$, and denote by $\tau_B$ the first exit time from $B$. From (35) and the strong Markov property at time $\tau_B$ we get for $x \in B$

$$
\begin{aligned}
u(x) &+ 2\Pi_x\Big(\int_0^{\tau_B} u(\xi_s)^2 ds\Big) + 2\Pi_x\Big(\Pi_{\xi_{\tau_B}}\Big(\int_0^\tau u(\xi_s)^2 ds\Big)\Big) \\
&= \Pi_x\big(\Pi_{\xi_{\tau_B}}(1_{\{\tau < \infty\}}g(\xi_\tau))\big).
\end{aligned}
$$

By combining this with formula (35) applied with $x = \xi_{\tau_B}$, we arrive at

$$u(x) + 2\Pi_x\Big(\int_0^{\tau_B} u(\xi_s)^2 ds\Big) = \Pi_x\big(u(\xi_{\tau_B})\big) \ . \tag{37}$$

The function $h(x) = \Pi_x\big(u(\xi_{\tau_B})\big)$ is harmonic in $B$, so that $h$ is of class $C^2$ and $\Delta h = 0$ in $B$. Set

$$f(x) := \Pi_x\Big(\int_0^{\tau_B} u(\xi)^2 ds\Big) = \int_B dy\, G_B(x, y) u(y)^2$$

where $G_B$ is the Green function of Brownian motion in $B$. Since $u$ is measurable and bounded, Theorem 6.6 of [37] shows that $f$ is continuously differentiable in $B$, and so is $u$ since $u = h - 2f$. Then again by Theorem 6.6 of [37], the previous formula for $f$ implies that $f$ is of class $C^2$ in $B$ and $-\frac{1}{2}\Delta f = u^2$ in $B$, which leads to the desired equation for $u$.

For the second part of the theorem, suppose first that $D$ is bounded, and let $y \in \partial D$ be regular for $D^c$. Then, if $g$ is continuous at $y$, it is well known that

$$\lim_{D \ni x \to y} \Pi_x\big(g(\xi_\tau)\big) = g(y) \ .$$

On the other hand, we have also

$$\limsup_{D \ni x \to y} \Pi_x\Big(\int_0^\tau u(\xi_s)^2 ds\Big) \leq \big(\sup_{x \in D} u(x)\big)^2 \limsup_{D \ni x \to y} E_x(\tau) = 0 \ .$$

Thus (35) implies that

$$\lim_{D \ni x \to y} u(x) = g(y) \ .$$

When $D$ is unbounded, a similar argument applies after replacing $D$ by $D \cap B$, where $B$ is now a large ball: Argue as in the derivation of (37) to verify that for $x \in D \cap B$,

$$u(x) + 2\Pi_x\Big(\int_0^{\tau_{D \cap B}} u(\xi_s)^2 ds\Big) = \Pi_x\big(1_{\{\tau \leq \tau_B\}} g(\xi_\tau)\big) + \Pi_x\big(1_{\{\tau_B < \tau\}} u(\xi_{\tau_B})\big)$$

and then follow the same route as in the bounded case. $\qquad\square$

The nonnegative solution of the problem (36) is always unique. When $D$ is bounded, this is a consequence of the following analytic lemma.

**Lemma 3.13** (Comparison principle) *Let $h : \mathbb{R}_+ \to \mathbb{R}_+$ be a monotone increasing function. Let $D$ be a bounded domain in $\mathbb{R}^d$ and let $u, v$ be two nonnegative functions of class $C^2$ on $D$ such that $\Delta u \geq h(u)$ and $\Delta v \leq h(v)$. Suppose that for every $y \in \partial D$,*

$$\limsup_{D \ni x \to y} \big(u(x) - v(x)\big) \leq 0 \ .$$

*Then $u \leq v$.*

**Proof.** Set $f = u - v$ and $D' = \{x \in D, f(x) > 0\}$. If $D'$ is not empty, we have

$$\Delta f(x) \geq h\big(u(x)\big) - h\big(v(x)\big) \geq 0$$

for every $x \in D'$. Furthermore, it follows from the assumption and the definition of $D'$ that

$$\limsup_{D' \ni x \to z} f(x) \leq 0$$

for every $z \in \partial D'$. Then the classical maximum principle implies that $f \leq 0$ on $D'$, which is a contradiction. $\qquad\square$

**Corollary 3.14** (Mean value property) *Let $D$ be a domain in $\mathbb{R}^d$ and let $U$ be a bounded regular subdomain of $D$ whose closure is contained in $D$. Then, if $u$ is a nonnegative solution of $\Delta u = 4\,u^2$ in $D$, we have for every $x \in U$*

$$u(x) = \mathbb{N}_x\big(1 - \exp -\langle \mathcal{Z}^U, u \rangle\big).$$

**Proof.** For every $x \in U$, set

$$v(x) = \mathbb{N}_x\big(1 - \exp -\langle \mathcal{Z}^U, u \rangle\big).$$

By Theorem 3.12, $v$ solves $\Delta v = 4v^2$ in $U$ with boundary value $v_{|\partial U} = u_{|\partial U}$. By Lemma 3.13, we must have $v(x) = u(x)$ for every $x \in U$. $\qquad\square$

The last proposition of this subsection provides some useful properties of nonnegative solutions of $\Delta u = 4u^2$ in a domain. For $x \in \mathbb{R}^d$ and $\varepsilon > 0$, we denote by $B(x, \varepsilon)$ the open ball of radius $\varepsilon$ centered at $x$. We also denote by

$$\mathcal{R} = \{\hat{W}_s, 0 \leq s \leq \sigma\}$$

the range of the Brownian snake.

**Proposition 3.15** (i) *There exists a positive constant $c_d$ such that for every $x \in \mathbb{R}^d$ and $\varepsilon > 0$,*

$$\mathbb{N}_x\big(\mathcal{R} \cap B(x, \varepsilon)^c \neq \varnothing\big) = c_d\, \varepsilon^{-2}\,.$$

(ii) *Let $u$ be a nonnegative solution of $\Delta u = 4u^2$ in the domain $D$. Then for every $x \in D$,*

$$u(x) \leq c_d\, \mathrm{dist}(x, \partial D)^{-2}\,.$$

(iii) *The set of all nonnegative solutions of $\Delta u = 4\,u^2$ in $D$ is closed under pointwise convergence.*

**Proof.** (i) By translation invariance we may assume that $x = 0$. We then use a scaling argument. For $\lambda > 0$, the law under $n(de)$ of $e_\lambda(s) = \lambda^{-1}e(\lambda^2 s)$ is $\lambda^{-1}n$ (exercise !). It easily follows that the law under $\mathbb{N}_0$ of $W_s^{(\varepsilon)}(t) = \varepsilon^{-1}W_{\varepsilon^4 s}(\varepsilon^2 t)$ is $\varepsilon^{-2}\mathbb{N}_0$. Then, with an obvious notation,

$$\mathbb{N}_0\big(\mathcal{R} \cap B(0, \varepsilon)^c \neq \varnothing\big) \;\;=\;\; \mathbb{N}_0\big(\mathcal{R}^{(\varepsilon)} \cap B(0, 1)^c \neq \varnothing\big)$$

$$=\;\; \varepsilon^{-2}\mathbb{N}_0\big(\mathcal{R} \cap B(0, 1)^c \neq \varnothing\big)\,.$$

It remains to verify that $\mathbb{N}_0\big(\mathcal{R} \cap B(0,1)^c \neq \varnothing\big) < \infty$. If this were not true, excursion theory would imply that $\mathbb{P}_0$ a.s., infinitely many excursions of the Brownian snake exit the ball $B(0,1)$ before time 1. Clearly this would contradict the continuity of $s \to W_s$ under $\mathbb{P}_0$.

(ii) Let $x \in D$ and $r > 0$ be such that $\bar{B}(x,r) \subset D$. By Corollary 3.14, we have for every $y \in B(x,r)$

$$u(y) = \mathbb{N}_y\big(1 - \exp -\langle \mathcal{Z}^{B(x,r)}, u\rangle\big) \ .$$

In particular,

$$u(x) \leq \mathbb{N}_x\big(\mathcal{Z}^{B(x,r)} \neq 0\big) \leq \mathbb{N}_x\big(\mathcal{R} \cap B(x,r)^c \neq \varnothing\big) = c_d \, r^{-2} \ .$$

In the second inequality we used the fact that $\mathcal{Z}^{B(x,r)}$ is supported on $\mathcal{E}^{B(x,r)} \subset \mathcal{R} \cap B(x,r)^c$.

(iii) Let $(u_n)$ be a sequence of nonnegative solutions of $\Delta u = 4\,u^2$ in $D$ such that $u_n(x) \longrightarrow u(x)$ as $n \to \infty$ for every $x \in D$. Let $U$ be an open ball whose closure is contained in $D$. By Corollary 3.14, for every $n \geq 1$ and $x \in U$,

$$u_n(x) = \mathbb{N}_x\big(1 - \exp -\langle \mathcal{Z}^U, u_n\rangle\big).$$

Note that $\mathbb{N}_x(\mathcal{Z}^U \neq 0) < \infty$ (by (i)) and that the functions $u_n$ are uniformly bounded on $\partial U$ (by (ii)). Hence we can pass to the limit in the previous formula and get $u(x) = \mathbb{N}_x\big(1 - \exp -\langle \mathcal{Z}^U, u\rangle\big)$ for $x \in U$. The desired result then follows from Theorem 3.12. □

Let us conclude this subsection with the following remark. Theorem 3.11 could be applied as well to treat parabolic problems for the operator $\Delta u - 4u^2$. To this end we need only replace the Brownian motion $\xi$ by the space-time process $(t, \xi_t)$. If we make this replacement and let $D \subset \mathbb{R}_+ \times \mathbb{R}^d$ be a space-time domain, then for every bounded nonnegative measurable function $g$ on $\partial D$, the formula

$$u(t,x) = \mathbb{N}_{t,x}\big(1 - \exp -\langle \mathcal{Z}^D, g\rangle\big)$$

gives a solution of

$$\frac{\partial u}{\partial t} + \frac{1}{2}\Delta u - 2u^2 = 0$$

in $D$. Furthermore, $u$ has boundary condition $g$ under suitable conditions on $D$ and $g$. The proof proceeds from the integral equation (35) as for Theorem 3.11.

### *3.6. Solutions with boundary blow-up*

**Proposition 3.16** *Let $D$ be a bounded regular domain. Then $u_1(x) = \mathbb{N}_x(\mathcal{Z}^D \neq 0)$, $x \in D$ is the minimal nonnegative solution of the problem*

$$\begin{aligned}\Delta u &= 4u^2 \qquad in \ D \\ u_{|\partial D} &= +\infty.\end{aligned} \tag{38}$$

**Proof.** First note that $u_1(x) < \infty$ by Proposition 3.15 (i). For every $n \geq 1$, set $v_n(x) = \mathbb{N}_x(1 - \exp -n\langle \mathcal{Z}^D, 1 \rangle)$, $x \in D$. By Theorem 3.12, $v_n$ solves (36) with $g = n$. By Proposition 3.15 (iii), $u_1 = \lim \uparrow v_n$ also solves $\Delta u = 4u^2$ in $D$. The condition $u_1|_{\partial D} = \infty$ is clear since $u_1 \geq v_n$ and $v_n|_{\partial D} = n$. Finally if

$v$ is another nonnegative solution of the problem (38), the comparison principle (Lemma 3.13) implies that $v \geq v_n$ for every $n$ and so $v \geq u_1$. □

**Proposition 3.17** *Let $D$ be any open set in $\mathbb{R}^d$ and $u_2(x) = \mathbb{N}_x(\mathcal{R} \cap D^c \neq \varnothing)$ for $x \in D$. Then $u_2$ is the maximal nonnegative solution of $\Delta u = 4u^2$ in $D$ (in the sense that $u \leq u_2$ for any other nonnegative solution $u$ in $D$).*

**Proof** First note that $\mathcal{R}$ is connected $\mathbb{N}_x$ a.e. as the range of the continuous mapping $s \to \widehat{W}_s$. It follows that we may deal separately with each connected component of $D$, and thus assume that $D$ is a domain. Then we can easily construct a sequence $(D_n)$ of bounded regular subdomains of $D$, such that $D = \lim \uparrow D_n$ and $\bar{D}_n \subset D_{n+1}$ for every $n$. Set

$$v_n(x) = \mathbb{N}_x(\mathcal{Z}^{D_n} \neq 0) \;,\; \tilde{v}_n(x) = \mathbb{N}_x(\mathcal{R} \cap D_n^c \neq \varnothing)$$

for $x \in D_n$. By the support property of the exit measure, it is clear that $v_n \leq \tilde{v}_n$. We also claim that $\tilde{v}_{n+1}(x) \leq v_n(x)$ for $x \in D_n$. To verify this, observe that on the event $\{\mathcal{R} \cap D_{n+1}^c \neq \varnothing\}$ there exists a path $W_s$ that hits $D_{n+1}^c$. For this path $W_s$, we must have $\tau_{D_n}(W_s) < \zeta_s$ (here $\tau_{D_n}$ stands for the exit time from $D_n$), and it follows from the properties of the Brownian snake that

$$A_\sigma^n := \int_0^\sigma dr \; 1_{\{\tau_{D_n}(W_r) < \zeta_r\}} > 0 \;,$$

$\mathbb{N}_x$ a.e. on $\{\mathcal{R} \cap D_{n+1}^c \neq \varnothing\}$. However, from the construction of the exit measure in subsection 3.4 above, $\langle \mathcal{Z}^{D_n}, 1 \rangle$ is obtained as the local time at level $0$ and at time $A_\sigma^n$ of a reflected Brownian motion started at $0$. Since the local time at $0$ of a reflected Brownian motion started at $0$ immediately becomes (strictly) positive, it follows that $\{\mathcal{R} \cap D_{n+1}^c \neq \varnothing\} \subset \{\mathcal{Z}^{D_n} \neq 0\}$, $\mathbb{N}_x$ a.e., which gives the inequality $\tilde{v}_{n+1}(x) \leq v_n(x)$.

We have then for $x \in D$

$$u_2(x) = \lim_{n \to \infty} \downarrow \tilde{v}_n(x) = \lim_{n \to \infty} \downarrow v_n(x) \;. \tag{39}$$

This follows easily from the fact that the event $\{\mathcal{R} \cap D^c \neq \varnothing\}$ is equal $\mathbb{N}_x$ a.e. to the intersection of the events $\{\mathcal{R} \cap D_n^c \neq \varnothing\}$. By Proposition 3.16, $v_n$ solves $\Delta u = 4u^2$ in $D_n$. It then follows from (39) and Proposition 3.15 (iii) that $u_2$ solves $\Delta u = 4u^2$ in $D$. Finally, if $u$ is another nonnegative solution in $D$, the comparison principle implies that $u \leq v_n$ in $D_n$ and it follows that $u \leq u_2$. □

**Example.** Let us apply the previous proposition to compute $\mathbb{N}_x(0 \in \mathcal{R})$ for $x \neq 0$. By rotational invariance and the same scaling argument as in the proof of Proposition 3.15 (i), we get $\mathbb{N}_x(0 \in \mathcal{R}) = C|x|^{-2}$ with a nonnegative constant $C$. On the other hand, by Proposition 3.17, we know that $u(x) = \mathbb{N}_x(0 \in \mathcal{R})$ solves $\Delta u = 4u^2$ in $\mathbb{R}^d \backslash \{0\}$. A short calculation, using the expression of the Laplacian for a radial function, shows that the only possible values of $C$ are $C = 0$ and $C = 2 - \frac{d}{2}$. Since $u$ is the maximal solution, we conclude that if $d \leq 3$,

$$\mathbb{N}_x(0 \in \mathcal{R}) = \left(2 - \frac{d}{2}\right)|x|^{-2}$$

whereas $\mathbb{N}_x(0 \in \mathcal{R}) = 0$ if $d \geq 4$. In particular, points are polar (in the sense that they are not hit by the range) if and only if $d \geq 4$.

Let us conclude with some remarks. First note that, if $D$ is bounded and regular (the boundedness is superfluous here), the function $u_2$ of Proposition 2 also satisfies $u_2|_{\partial D} = +\infty$. This is obvious since $u_2 \geq u_1$. We may ask the following two questions.

**1.** If $D$ is regular, is it true that $u_1 = u_2$? (uniqueness of the solution with boundary blow-up)

**2.** For a general domain $D$, when is it true that $u_2|_{\partial D} = +\infty$? (existence of a solution with boundary blow-up)

A complete answer to question **2** is provided in [8] (see also [27]). A general answer to **1** is still an open problem (see [27] and the references therein for partial results).

*Bibliographical notes.* Much of this section is taken from [27], where additional references about the Brownian snake can be found. The connections with partial differential equations that are discussed in subsections 3.5 and 3.6 were originally formulated by Dynkin [13] in the language of superprocesses (see Perkins [35] for a recent account of the theory of superprocesses). These connections are still the subject of an active research: See Dynkin's books [14], [15]. Mselati's thesis [33] gives an application of the Brownian snake to the classification and probabilistic representation of the solutions of $\Delta u = u^2$ in a smooth domain.

# References

[1] D. Aldous. The continuum random tree I. *Ann. Probab.* **19** (1991), 1-28. MR1085326

[2] D. Aldous. The continuum random tree III. *Ann. Probab.* **21** (1993), 248-289. MR1207226

[3] D. Aldous. Tree-based models for random distribution of mass. *J. Stat. Phys.* **73** (1993), 625-641. MR1251658

[4] R.M. Blumenthal. *Excursions of Markov Processes.* Birkhäuser 1992.

[5] M. Bramson, J.T. Cox, J.F. Le Gall. Super-Brownian limits of voter model clusters. *Ann. Probab.* **29** (2001), 1001-1032. MR1872733

[6] P. Chassaing, G. Schaeffer. Random planar lattices and integrated superBrownian excursion. *Probab. Th. Rel. Fields* **128** (2004), 161-212. MR2031225

[7] E. Derbez, G. Slade. The scaling limit of lattice trees in high dimensions. *Comm. Math. Phys.* **198** (1998), 69-104. MR1620301

[8] J.S. Dhersin, J.F. Le Gall. Wiener's test for super-Brownian motion and for the Brownian snake. *Probab. Th. Rel. Fields* **108** (1997), 103-129. MR1452552

[9] M. Drmota, B. Gittenberger. On the profile of random trees. *Random Struct. Alg.* **10** (1997), 421-451. MR1608230

[10] T. Duquesne. A limit theorem for the contour process of conditioned Galton-Watson trees. *Ann. Probab.* **31** (2003), 996–1027. MR1964956

[11] T. Duquesne, J.F. Le Gall. Random Trees, Lévy Processes and Spatial Branching Processes. *Astérisque* **281** (2002)

[12] T. Duquesne, J.F. Le Gall. Probabilistic and fractal aspects of Lévy trees. *Probab. Th. Rel. Fields* **131** (2005), 553-603. MR2147221

[13] E.B. Dynkin. A probabilistic approach to one class of nonlinear differential equations. *Probab. Th. Rel. Fields* **90** (1991), 89-115. MR1109476

[14] E.B. Dynkin. *Diffusions, Superdiffusions and Partial Differential Equations.* American Math. Society 2002. MR1883198

[15] E.B. Dynkin. *Superdiffusions and Positive Solutions of Nonlinear Partial Differential Equations.* American Math. Society 2004. MR2089791

[16] S.N. Ethier, T. Kurtz. *Markov Processes, Characterization and Convergence.* Wiley 1986.

[17] S.N. Evans, J. Pitman, A. Winter. Rayleigh processes, real trees, and root growth with re-grafting. *Probab. Th. Rel. Fields*, to appear.

[18] S.N. Evans, A. Winter. Subtree prune and re-graft: A reversible tree-valued Markov process. *Ann. Probab.*, to appear. MR1159560

[19] M. Gromov. *Metric Structures for Riemannian and Non-Riemannian Spaces.* Birkhäuser 1999.

[20] B. Haas, G. Miermont. The genealogy of self-similar fragmentations with negative index as a continuum random tree. *Electron. J. Probab.* **9** (2004), 57-97. MR2041829

[21] T. Hara, G. Slade. The scaling limit of the incipient infinite cluster in high-dimensional percolation. II. Integrated super-Brownian excursion. Probabilistic techniques in equilibrium and nonequilibrium statistical physics. *J. Math. Phys.* **41** (2000), 1244-1293. MR1757958

[22] R. van der Hofstad, A. Sakai. Gaussian scaling for the critical spread-out contact process above the upper critical dimension. *Electron. J. Probab.* **9** (2004), 710-769. MR2110017

[23] R. van der Hofstad, G. Slade. Convergence of critical oriented percolation to super-Brownian motion above $4 + 1$ dimensions. *Ann. Inst. H. Poincaré Probab. Statist.* **39** (2003), 413-485. MR1978987

[24] W.D. Kaigh. An invariance principle for random walk conditioned by a late

return to zero. *Ann. Probab.* **4** (1976), 115-121. MR415706

[25] J.F. Le Gall. Brownian excursions, trees and measure-valued processes. *Ann. Probab.* **19** (1991), 1399-1439. MR1127710

[26] J.F. Le Gall. The uniform random tree in a Brownian excursion. *Probab. Theory Rel. Fields* **95** (1993), 25-46. MR1207305

[27] J.F. Le Gall. *Spatial Branching Processes, Random Snakes and Partial Differential Equations.* Birkhäuser 1999.

[28] J.F. Le Gall. *Random trees and spatial branching processes.* Maphysto Lecture Notes Series (University of Aarhus), vol. 9 (80 pp.), 2000. Available at `http://www.maphysto.dk/publications/MPS-LN/2000/9.pdf`

[29] J.F. Le Gall, Y. Le Jan. Branching processes in Lévy processes: The exploration process. *Ann. Probab.* **26** (1998), 213-252. MR1617047

[30] J.F. Le Gall, M. Weill. Conditioned Brownian trees. *Ann. Inst. H. Poincaré Probab. Statist.*, to appear.

[31] J.F. Marckert, A. Mokkadem. The depth first processes of Galton-Watson trees converge to the same Brownian excursion. *Ann. Probab.* **31** (2003), 1655-1678. MR1989446

[32] J.-F. Marckert, G. Miermont. Invariance principles for labeled mobiles and bipartite planar maps. Preprint (2005).

[33] B. Mselati. Classification and probabilistic representation of the positive solutions of a semilinear elliptic equation. *Memoirs Amer. Math. Soc.* **798** (2004)

[34] J. Neveu. Arbres et processus de Galton-Watson. *Ann. Inst. Henri Poincaré Probab. Stat.* **22** (1986), 199-207. MR850756

[35] E.A. Perkins. Dawson-Watanabe superprocesses and measure-valued diffusions. Ecole d'été de probabilités de Saint-Flour 1999. In: *Lecture Notes Math.* **1781**. Springer, 2002. MR1915445

[36] J.W. Pitman. Combinatorial stochastic processes. Ecole d'été de probabilités de Saint-Flour 2002. *Lecture Notes Math.*, Springer. To appear.

[37] S.C. Port, C.J. Stone. *Brownian Motion and Classical Potential Theory.* Academic, New York, 1978. MR492329

[38] D. Revuz, M. Yor. *Continuous Martingales and Brownian Motion.* Springer 1991.

[39] L.C.G. Rogers, D. Williams. *Diffusions, Markov Processes and Martingales. Vol.2: Itô Calculus.* Wiley 1987. MR921238

[40] F. Spitzer. *Principles of Random Walk.* Van Nostrand 1963.

[41] R.P. Stanley. *Enumerative Combinatorics, Vol.2.* Cambridge University Press 1999